\documentclass[11pt]{amsart}
\usepackage{amsmath,amssymb,amsthm,latexsym,cite}
\usepackage[english]{babel}
\usepackage{xcolor}
\usepackage{enumitem,bm,booktabs,longtable,etoolbox}
\usepackage[left=2.1cm,right=2.1cm,top=2.71cm,bottom=2.71cm]{geometry}
\usepackage[colorlinks=true,urlcolor=black,citecolor=black,linkcolor=black,
linktocpage,pdfpagelabels,bookmarksnumbered,bookmarksopen]{hyperref}
\makeatletter
\patchcmd{\H@old@sect}{\@tocwrite{#1}{#8}}{\@tocwrite{#1}{#7}}{}{}
\patchcmd{\H@old@sect}{\@tocwrite{#1}{#8}}{\@tocwrite{#1}{#7}}{}{}
\def\@seccntformat#1{\protect\textup{\protect\@secnumfont
\csname the#1\endcsname\protect\enspace}}
\makeatother
\newtheorem{theorem}{Theorem}[section]
\newtheorem{definition}[theorem]{Definition}
\newtheorem{proposition}[theorem]{Proposition}
\newtheorem{lemma}[theorem]{Lemma}
\theoremstyle{remark}
\newtheorem{remark}[theorem]{Remark}
\theoremstyle{plain}
\newtheorem{corollary}[theorem]{Corollary}

\newcommand{\R}{\mathbb R}

\numberwithin{equation}{section}
\title[Hardy--Sobolev--Maz'ya inequalities]
{Hardy--Sobolev--Maz'ya inequalities on the lattice graph \(\mathbb Z^N\)}
\author[C. Ji]{Chao Ji}

\address[C. Ji]{\newline\indent
	School of Mathematics
	\newline\indent
	East China University of Science and Technology
	\newline\indent
	Shanghai 200237, PR China }
\email{\href{mailto:jichao@ecust.edu.cn}{jichao@ecust.edu.cn}}

\author[K. Sheng]{Kai Sheng}

\address[K. Sheng]{\newline\indent
	School of Mathematics
	\newline\indent
	East China University of Science and Technology
	\newline\indent
	Shanghai 200237, PR China}

\email{\href{mailto:shengkai2001@outlook.com}{shengkai2001@outlook.com}}

\date{October 20, 2026}
\begin{document}
\emergencystretch=4em
\begin{abstract}
In this paper, we study the following
Hardy--Sobolev--Maz'ya inequality
on the lattice graph \(\mathbb Z^N\):
\[
  \sum_{j=1}^N\sum_{z\in\mathbb Z^N}|u(z)-u(z-e_j)|^2
  -C_H(n,m)\sum_{z\in\mathbb Z^N\setminus\Sigma}
  \frac{|u(z)|^2}{d_\Sigma(z)^2}
  \geq\overline\kappa_{m,n}\|u\|_{\ell^{2^*}(\mathbb Z^N)}^2,
  \quad u\in C_c(\mathbb Z^N\setminus\Sigma),
\]
where \(n,m\geq1\), \(N=n+m\geq3\),
\(\Sigma=\mathbb Z^n\times\{0\}\), \(2^*=2N/(N-2)\), \(e_j\) is the \(j\)-th coordinate vector, \(d_\Sigma(z):=\inf_{\zeta\in\Sigma}|z-\zeta|\), and \(C_H(n,m)\) is the
best constant for the Hardy inequality on \(\mathbb Z^N\setminus\Sigma\).
 Using the supersolution method and multilinear interpolation, we first prove that the best constant
\(\overline\kappa_{m,n}\) is positive if and only if either
\(1\leq m\leq8\) and \(n\geq1\), or \(m=9\) and \(n\geq14\). Then, we show that
\[
  \overline\kappa_{m,n}=\overline S_N
  \quad\text{if }m=2\text{ or }(m,n)=(1,2),
\]
where \(\overline S_N\) is the best constant in the critical
Sobolev inequality on \(D^{1,2}(\mathbb Z^N)\). Moreover,
for \(m\neq2\) and \(N\geq4\), we prove
\(0\leq\overline\kappa_{m,n}<\overline S_N\).  Finally, using the equality
\(\overline\kappa_{2,n}=\overline S_{n+2}\), we prove that \(\overline\kappa_{2,n}\) is
not attained, while for \((m,n)=(1,2)\), we establish nonattainment of
\(\overline\kappa_{1,2}\) by developing a duality argument based on a comparison of Green functions.
\end{abstract}
\maketitle

\section{Introduction}

In this paper, we study the following
Hardy--Sobolev--Maz'ya inequality on \(\mathbb Z^N\):
\begin{equation}
  \label{eq:discrete-HSM}
  \sum_{z\in\mathbb Z^N}|Du(z)|^2
  -C_H(n,m)\sum_{z\in\mathbb Z^N\setminus\Sigma}
  \frac{|u(z)|^2}{d_\Sigma(z)^2}
  \geq\overline\kappa_{m,n}\|u\|_{\ell^{2^*}(\mathbb Z^N)}^2,
  \quad u\in C_c(\mathbb Z^N\setminus\Sigma),
\end{equation}
where \(N=n+m\geq3\), \(n=1,\ldots,N-1\),
\(z=(x,y)\in\mathbb Z^n\times\mathbb Z^m\),
\(\Sigma=\mathbb Z^n\times\{0\}\), \(2^*=2N/(N-2)\), \(d_\Sigma(z):=\inf_{\zeta\in\Sigma}|z-\zeta|\), and
\(C_c(\mathbb Z^N\setminus\Sigma)\) is the space of finitely supported
functions on \(\mathbb Z^N\setminus\Sigma\).
For \(j=1,\ldots,N\), let \(e_j\) be the \(j\)-th coordinate vector.
To define \(Du(z):=(D_1u(z),\ldots,D_Nu(z))\), where
\(D_ju(z)=u(z)-u(z-e_j)\), for every \(z\in\mathbb Z^N\),
we extend each \(u\in C_c(\mathbb Z^N\setminus\Sigma)\) by
\(u(z)=0\) for \(z\in\Sigma\).

We denote by \(\overline\kappa_{m,n}\) the best constant in \eqref{eq:discrete-HSM}, and by \(C_H(n,m)\) the best constant in the following
 Hardy inequality:
\begin{equation}
  \label{eq:hardy-sigma}
  \sum_{z\in\mathbb Z^N}|Du(z)|^2
  \geq C_H(n,m)\sum_{z\in\mathbb Z^N\setminus\Sigma}
  \frac{|u(z)|^2}{d_\Sigma(z)^2},
  \qquad u\in C_c(\mathbb Z^N\setminus\Sigma).
\end{equation}

\par

Inequality \eqref{eq:discrete-HSM} arises from the study of remainder terms in the sharp Hardy inequality.
There is an extensive literature on such remainder terms, including terms
involving \(L^2\) norms, logarithmic corrections and Sobolev-type remainder terms.
On a bounded domain \(\Omega\subset\mathbb R^N\) containing the origin,
\(N\geq3\), Br\'ezis and V\'azquez \cite[Theorem~4.1]{BrezisVazquez1997}
obtained a remainder \(C_\Omega\|u\|_{L^2(\Omega)}^2\) in the Hardy inequality on \(H_0^1(\Omega)\), where \(C_\Omega>0\). When \(\Omega\) is a ball, they determined the optimal coefficient
and asked whether the resulting inequality admits further remainder terms
\cite[Section~8, Problem~2]{BrezisVazquez1997}. Motivated by this question,
Filippas and Tertikas \cite[Theorems~B and~D]{FilippasTertikas2002} studied when further improvements are possible. For a logarithmic weight, they obtained the remainder
\[
  \frac14\int_\Omega
  \frac{|u(x)|^2}{|x|^2[1+\log(D/|x|)]^2}\,dx,
  \qquad D\geq\sup_{x\in\Omega}|x|,
\]
and proved that the coefficient \(1/4\) is optimal.
Barbatis, Filippas and Tertikas \cite[Theorem~A]{BarbatisFilippasTertikas2003}
subsequently extended logarithmic improvements to \(L^p\) Hardy inequalities.
Another direction concerns Sobolev-type remainder terms.
Filippas, Maz'ya and Tertikas \cite[Theorem~1.1]{FilippasMazyaTertikas2007}
proved
\[
  \int_\Omega|\nabla u|^2\,dx
  -\frac14\int_\Omega\frac{|u|^2}{d(x)^2}\,dx
  \geq C_\Omega\|u\|_{L^{2^*}(\Omega)}^2,
  \qquad u\in C_c^\infty(\Omega),
\]
where \(\Omega\subset\mathbb R^N\) (\(N\geq3\)) is a bounded convex domain
with \(C^2\) boundary, \(d(x)=\operatorname{dist}(x,\partial\Omega)\),  \(C_\Omega>0\) and \(2^*=2N/(N-2)\).

In addition, the Sobolev remainder in the Hardy inequality on the Euclidean
upper half-space has a geometric interpretation in terms of the Poincar\'e
inequality on hyperbolic space. For \(m=1\) and \(N=n+1\geq3\), the upper
half-space \(\mathbb R^n\times(0,\infty)\) equipped with the metric
\(g=y^{-2}(dx^2+dy^2)\) is a model of hyperbolic space \(\mathbb H^N\).
The substitution \(v=y^{(N-2)/2}u\) gives the Poincar\'e--Sobolev inequality
\cite[Appendix~B]{ManciniSandeep2008}
\[
  \int_{\mathbb H^N}|\nabla_{\mathbb H}v|^2\,dV_{\mathbb H}
  -\frac{(N-1)^2}{4}\int_{\mathbb H^N}|v|^2\,dV_{\mathbb H}
  \geq C\|v\|_{L^{2^*}(\mathbb H^N)}^2,
  \qquad v\in C_c^\infty(\mathbb H^N),
\]
where \(C>0\), and \(\nabla_{\mathbb H}\) and \(dV_{\mathbb H}\) denote
the hyperbolic gradient and volume element, respectively.
Here \((N-1)^2/4\) is the best constant in the Poincar\'e inequality
on \(\mathbb H^N\) \cite[(1.1)]{ManciniSandeep2008}.
Moreover,
\(\|v\|_{L^{2^*}(\mathbb H^N)}=\|u\|_{L^{2^*}(\mathbb R^n\times(0,\infty))}\).
Thus the Poincar\'e--Sobolev inequality on \(\mathbb H^N\) and the
Hardy--Sobolev--Maz'ya inequality on \(\mathbb R^n\times(0,\infty)\)
have the same best constant.
For related results on Poincar\'e--Sobolev and Hardy--Sobolev--Maz'ya
inequalities on hyperbolic spaces, see \cite{BFL2008,LuYang2019,ManciniSandeep2008}
and the references therein. In \eqref{eq:discrete-HSM}, we are interested in an \(\ell^{2^*}\)
remainder term of the form $\overline\kappa_{m,n}\|u\|_{\ell^{2^*}(\mathbb Z^N)}^2$. Two natural questions arise: whether $\overline\kappa_{m,n}>0$
and whether this best constant is attained.
For further results on remainder terms in Hardy inequalities, we refer to \cite{BrezisVazquez1997,BarbatisFilippasTertikas2003,FilippasTertikas2002,FilippasMazyaTertikas2007,FilippasTertikasTidblom2009,FrankLoss2012} and the references therein.

The question of attainment is also closely related to the existence of solutions to the nonlinear Schr\"odinger equation with a Hardy potential
and a critical nonlinearity on \(\mathbb Z^N\)
(see, for instance, \cite{Smets2005,HeJi2026}). Indeed, if $\overline\kappa_{m,n}>0$ is attained by some $u$
satisfying $\|u\|_{\ell^{2^*}(\mathbb Z^N)}=1$,
then $u$ solves the corresponding Euler--Lagrange equation
\[
\Delta_{\mathrm d}u(z)
-C_H(n,m)\frac{u(z)}{d_\Sigma(z)^2}
=\overline\kappa_{m,n}|u(z)|^{2^*-2}u(z),
\qquad z\in\mathbb Z^N\setminus\Sigma,
\]
with $u|_\Sigma=0$, where $\Delta_{\mathrm d}$ denotes
the discrete Laplacian.
Thus, the study of the best constant in \eqref{eq:discrete-HSM}
is naturally connected to the existence of solutions to this nonlinear Schr\"odinger equation on lattice graphs.

\subsection{Literature review}

The Hardy--Sobolev--Maz'ya inequality was first studied by Maz'ya
\cite[Corollary~3]{Mazya1985}.
For \(n,m\geq1\) and \(N=n+m\geq3\), this inequality takes the form
\[
  \int_{\mathbb R^N}|\nabla u|^2\,dz
  -\frac{(m-2)^2}{4}\int_{\mathbb R^N}\frac{|u|^2}{|y|^2}\,dz
  \geq\kappa_{m,n}\|u\|_{L^{2^*}(\mathbb R^N)}^2,
\]
where $u\in C_c^\infty\bigl(\mathbb R^N\setminus(\mathbb R^n\times\{0\})\bigr)$, \(z=(x,y)\in\mathbb R^n\times\mathbb R^m\),
\(2^*=2N/(N-2)\), and \(\kappa_{m,n}\geq0\) denotes the best constant.
 By applying a weighted Sobolev inequality to \(|y|^{(m-2)/2}u\),
Maz'ya \cite[Corollary~3]{Mazya1985} proved
\(\kappa_{m,n}>0\) for \(n,m\geq1\) with \(n+m\geq3\).
Once \(\kappa_{m,n}>0\), it is natural to study its exact value and attainment.
By translating test functions away from \(\mathbb R^n\times\{0\}\)
in the \(y\)-variables, it is easy to see that
\[
  \kappa_{m,n}\leq S_N,
\]
where \(S_N\) is the best constant in the critical Sobolev inequality
on \(\mathbb R^N\). This raises the question of when \(\kappa_{m,n}<S_N\)
and when \(\kappa_{m,n}=S_N\).
Tertikas and Tintarev \cite[Lemma~2.2]{TertikasTintarev2007}
proved  \(\kappa_{m,n}<S_N\) when \(m\geq3\) and \(n\geq1\),
or \(m=1\) and \(n\geq3\).
Their proof uses a Sobolev extremal function \cite{Talenti1976}, together with scaling
and a cutoff when \(m=1\).
However, when \((m,n)=(1,2)\), the Hardy term and the increase in the
Sobolev quotient due to the cutoff have the same order, so this argument
does not give \(\kappa_{1,2}<S_3\).
To determine \(\kappa_{1,2}\), Benguria, Frank and Loss
\cite[Theorem~1.1]{BFL2008} used a duality argument based on a comparison
of Green functions and proved  \(\kappa_{1,2}=S_3\).

 The next question
is whether \(\kappa_{m,n}\) is attained.
The main difficulty in proving attainment is the lack of compactness.
To overcome this difficulty, Tertikas and Tintarev combined the strict
inequality \(\kappa_{m,n}<S_N\) with a concentration--compactness
argument to prove the existence of extremal functions in the above
ranges \cite[Theorem~1.2]{TertikasTintarev2007}.
For \((m,n)=(1,2)\), Benguria, Frank and Loss
\cite[Theorem~1.1]{BFL2008} proved by a duality argument that \(\kappa_{1,2}\) is not attained. See also
\cite{FilippasMazyaTertikas2007,FilippasTertikasTidblom2009,FrankLoss2012,PinchoverTintarev2009}
and the references therein for further results on Hardy--Sobolev--Maz'ya inequalities.

In recent years, functional inequalities on lattice graphs have attracted considerable attention (see \cite{HuaLi2021,Gupta2023,HajaiejHanHua2026,HanLi2026,HuangYe2026} and the references therein).
As mentioned above, comparison with the sharp Sobolev constant plays an important role in the study of the Hardy--Sobolev--Maz'ya inequality.
We therefore first recall some results on the best constants in Sobolev inequalities on \(\mathbb Z^N\) and their attainment.
Hua and Li \cite[Theorem~1]{HuaLi2021} developed a discrete concentration--compactness principle and proved the existence of extremal functions for \(N\geq3\) and \(q>2^*\).
A key step in their proof is to exclude vanishing by means of an interpolation estimate.
However, this estimate requires \(q>2^*\) and does not apply to the critical case \(q=2^*\).
To overcome this difficulty, He and Ji \cite[Theorem~1.8]{HeJi2026} subsequently established \(\overline S_N<S_N\) for \(N\geq3\), where \(\overline S_N\) is the best constant in the critical Sobolev inequality on \(\mathbb Z^N\).
Combining this comparison with singular limit analysis on lattice graphs, they ruled out vanishing and proved that \(\overline S_N\) is attained.
For related results on Sobolev, Hardy--Littlewood--Sobolev, and Caffarelli--Kohn--Nirenberg inequalities on lattice graphs, see \cite{HajaiejHanHua2026,HanLi2026,HuangLiYin2015,HuaLiMunch2025} and the references therein.

To study the positivity of \(\overline\kappa_{m,n}\), we recall results
on the value of \(c_m\) and its growth as \(m\to\infty\), where \(c_m\)
is the best constant in the following Hardy inequality
\[
  \sum_{j=1}^m\sum_{y\in\mathbb Z^m}|a(y)-a(y-e_j)|^2
  \geq c_m\sum_{y\in\mathbb Z^m\setminus\{0\}}\frac{|a(y)|^2}{|y|^2},
  \qquad a\in C_c(\mathbb Z^m\setminus\{0\}),
\]
with \(a(0)=0\). For each fixed dimension \(m\), a natural question is whether
\(c_m=(m-2)^2/4\), where $(m-2)^2/4$ is the best constant in the corresponding Hardy
inequality on \(\mathbb R^m\).
For \(m\geq3\), Lizama \cite{Lizama2026} proved
\(c_m=(m-2)^2/4\) for \(m=3,4\) by constructing positive supersolutions
and applying the ground-state transformation.
He also  proved \(c_m<(m-2)^2/4\) for \(m\geq9\).
However,  the supersolution construction in \cite{Lizama2026} does not give the lower bound
\(c_m\geq(m-2)^2/4\) for \(5\leq m\leq8\).
Alpay \cite[Theorem~2.1]{AlpayLow2026} subsequently extended the
reciprocal-edge method to prove \(c_m=(m-2)^2/4\) for \(5\leq m\leq8\).

Concerning the growth of \(c_m\), Gupta \cite[Corollary~1.5]{Gupta2023}
used Fourier methods and weighted Hardy inequalities on a torus to prove
\[
  Am\leq c_m\leq Bm
  \qquad\text{for all sufficiently large }m,
\]
where \(A,B>0\) are independent of \(m\).
He also asked for the value of \(\lim_{m\to\infty}c_m/m\)
\cite[Remark~1.9]{Gupta2023}.
Subsequently, Huang and Ye
\cite[Theorem~1.1]{HuangYe2026} answered this question by combining Fourier methods with weighted
Hardy--Rellich identities and probabilistic concentration estimates,
and proved
\begin{equation}
  \label{eq:hardy-constant-asymptotic}
  \lim_{m\to\infty}\frac{c_m}{m}=2.
\end{equation}
Independently, Gupta \cite[Theorem~1]{GuptaSharp2026}
proved \eqref{eq:hardy-constant-asymptotic} by using a ground-state
representation and the Poincar\'e inequality on the torus.

Motivated by \cite{Mazya1985,TertikasTintarev2007,HuangYe2026}, we study the Hardy--Sobolev--Maz'ya inequality~\eqref{eq:discrete-HSM} on \(\mathbb Z^N\).
Our main questions are whether \(\overline\kappa_{m,n}>0\), how \(\overline\kappa_{m,n}\) compares with \(\overline S_N\), and whether \(\overline\kappa_{m,n}\) is attained. Compared with the Euclidean setting, two main difficulties arise in proving positivity and comparing the best constants.
First, the standard supersolution method for proving positivity requires a positive supersolution with suitable decay at infinity.
In the Euclidean setting, such a supersolution can be constructed explicitly.
However, for the discrete problem with \(m\geq3\), \cite[Theorem~4.2]{KellerPinchoverPogorzelski2020} guarantees the existence of a positive supersolution but provides neither an explicit expression nor the required decay estimates.
Second, scaling plays an important role in constructing test functions to establish \(\kappa_{1,n}<S_N\) for \(n\geq3\) in the Euclidean setting \cite[proof of Lemma~2.2]{TertikasTintarev2007}.
However, the corresponding scaling argument is not available on \(\mathbb Z^N\).
These difficulties prevent us from directly applying the methods developed for the Euclidean problem.

\raggedbottom
\subsection{Main results}

 We first study the positivity of \(\overline\kappa_{m,n}\).

\begin{theorem}
\label{thm:main}\samepage
Let \(n,m\geq1\) and \(N=n+m\geq3\). Then
\[
  \overline\kappa_{m,n}>0
  \quad\Longleftrightarrow\quad
  1\leq m\leq8,\ n\geq1\quad\text{or}\quad m=9,\ n\geq14,
\]
where
\[
  \overline\kappa_{m,n}
  =
  \inf_{\substack{u\in C_c(\mathbb Z^N\setminus\Sigma)\\u\neq0}}
  \frac{\displaystyle\sum_{z\in\mathbb Z^N}|Du(z)|^2
  -C_H(n,m)\sum_{z\in\mathbb Z^N\setminus\Sigma}
  \frac{|u(z)|^2}{d_\Sigma(z)^2}}
  {\displaystyle\|u\|_{\ell^{2^*}(\mathbb Z^N)}^2}.
\]
\end{theorem}
\par

To study the value of \(\overline\kappa_{m,n}\), let \(\overline S_N\)
denote the best constant in the critical Sobolev inequality on \(\mathbb Z^N\)
(see \cite{HuaLi2021})
\[
  \overline S_N:=\inf_{u\in C_c(\mathbb Z^N)\setminus\{0\}}
  \frac{\displaystyle\sum_{z\in\mathbb Z^N}|Du(z)|^2}
  {\displaystyle\|u\|_{\ell^{2^*}(\mathbb Z^N)}^2}>0.
\]
By translating test functions away from \(\Sigma\) in the \(y\)-variables,
it is easy to see that
\[
  \overline\kappa_{m,n}\leq\overline S_N
\]
(see Lemma~\ref{lem:general-sobolev-upper}).
This raises the question of when \(\overline\kappa_{m,n}<\overline S_N\)
and when \(\overline\kappa_{m,n}=\overline S_N\).
The comparison also plays a key role in the study of attainment (for more details, see Corollary~\ref{cor:nonattainment-m2}).

Our second result reads as follows.

\begin{theorem}
\label{thm:sobolev-comparison}
Let \(n,m\geq1\) and \(N=n+m\geq3\). Then
\[
  \begin{cases}
    \overline\kappa_{m,n}=\overline S_N,
    &m=2\text{ or }(m,n)=(1,2),\\[2pt]
    0\leq\overline\kappa_{m,n}<\overline S_N,
    &m\neq2\text{ and }N\geq4.
  \end{cases}
\]
\end{theorem}

 We next study whether \(\overline\kappa_{m,n}\) is attained.
For \(\overline\kappa_{m,n}>0\), let \(\mathcal H_{m,n}\) be the
completion of \(C_c(\mathbb Z^N\setminus\Sigma)\) under the norm
\[
  \|u\|_{\mathcal H_{m,n}}
  :=\left(
  \sum_{z\in\mathbb Z^N}|Du(z)|^2
  -C_H(n,m)\sum_{z\in\mathbb Z^N\setminus\Sigma}
  \frac{|u(z)|^2}{d_\Sigma(z)^2}
  \right)^{1/2},
  \qquad u\in C_c(\mathbb Z^N\setminus\Sigma).
\]
With this choice of norm, every \(u\in\mathcal H_{m,n}\) belongs to \(\ell^{2^*}(\mathbb Z^N)\)
and satisfies \(u(z)=0\) for every \(z\in\Sigma\)
(see Lemma~\ref{lem:Hmn-embedding}).

As a byproduct of Theorem~\ref{thm:sobolev-comparison}, we first show
that \(\overline\kappa_{2,n}\) is not attained in \(\mathcal H_{2,n}\)
for \(n\geq1\).

\begin{corollary}
\label{cor:nonattainment-m2}
Let \(m=2\) and \(n\geq1\). Then
\(\overline\kappa_{2,n}=\overline S_{n+2}>0\), and  $\overline\kappa_{2,n}$
is not attained in \(\mathcal H_{2,n}\).
\end{corollary}

 We also prove that \(\overline\kappa_{1,2}\) is not attained.

\begin{theorem}
\label{thm:attainment}
Let \((m,n)=(1,2)\). Then
\(\overline\kappa_{1,2}=\overline S_3>0\), and  $\overline\kappa_{1,2}$
is not attained in \(\mathcal H_{1,2}\).
\end{theorem}

\begin{remark}
It remains to determine whether \(\overline\kappa_{m,n}\) is attained in
\(\mathcal H_{m,n}\) when \(0<\overline\kappa_{m,n}<\overline S_N\).
A main difficulty is to rule out \(\|u_j\|_{\ell^\infty}\to0\) for minimizing
sequences with \(\|u_j\|_{\ell^{2^*}}=1\).
For \(m=1\) and \(n\geq3\), a possible approach is to use the ground-state
transformation and interpolation to prove that such a sequence satisfies
\[
  \liminf_{j\to\infty}\|u_j\|_{\mathcal H_{1,n}}^2\geq\kappa_{1,n},
\]
where \(\kappa_{1,n}\) is the best constant in the Hardy--Sobolev--Maz'ya
inequality on \(\mathbb R^{n+1}_+\).
This would exclude vanishing if \(\overline\kappa_{1,n}<\kappa_{1,n}\).

\end{remark}

To the best of our knowledge, the Hardy--Sobolev--Maz'ya inequality~\eqref{eq:discrete-HSM} on \(\mathbb Z^N\) has not been established.
Theorem~\ref{thm:main} gives the first complete classification of the
dimensions for which its best constant \(\overline\kappa_{m,n}\) is positive.
It is worth pointing out that Theorem~\ref{thm:main} differs from Maz'ya's
result for the Hardy--Sobolev--Maz'ya inequality
\cite[Corollary~3, p.~97]{Mazya1985}.
In the continuous case, the best constant for the Hardy--Sobolev--Maz'ya
inequality is positive for every
\(n,m\geq1\) with \(n+m\geq3\).
In contrast, we prove \(\overline\kappa_{m,n}=0\) when
\(m=9\) and \(1\leq n\leq13\), or when \(m\geq10\) and \(n\geq1\).
Moreover, Theorem~\ref{thm:sobolev-comparison} can be regarded as a discrete
analogue of \cite[Lemma~2.2]{TertikasTintarev2007} and
\cite[Theorem~1.1]{BFL2008}.
 When \(m=2\), we have  \(C_H(n,2)=0\) (see Proposition~\ref{prop:m2}), and hence
inequality~\eqref{eq:discrete-HSM} reduces to the critical Sobolev inequality on
\(C_c(\mathbb Z^{n+2}\setminus\Sigma)\).
For the critical Sobolev inequality on \(D^{1,2}(\mathbb Z^{n+2})\),
the best constant is attained \cite[Theorem~1.8]{HeJi2026}.
However, Corollary~\ref{cor:nonattainment-m2} shows that the best
constant in \eqref{eq:discrete-HSM} is not attained in \(\mathcal H_{2,n}\). Thus, although the best constants coincide, their attainment properties differ because the admissible spaces are different. Finally, Theorem~\ref{thm:attainment} is a discrete counterpart of
\cite[Theorem~1.1]{BFL2008}.
\par

\subsection{Difficulties and ideas of the proof}\label{sec:proof-ideas}

Our study of \eqref{eq:discrete-HSM}  focuses on the following three problems:
\begin{enumerate}[label=(\Alph*),leftmargin=2.4em,itemsep=2pt]
\item determining the dimensions \(n\) and \(m\) for which
\(\overline\kappa_{m,n}>0\);
\item identifying the cases in which \(\overline\kappa_{m,n}<\overline S_N\)
or \(\overline\kappa_{m,n}=\overline S_N\);
\item studying the attainment of \(\overline\kappa_{m,n}\) in
\(\mathcal H_{m,n}\) when \(\overline\kappa_{m,n}>0\).
\end{enumerate}

With regard to (A), we first prove \(C_H(n,m)=c_m\).
We use the identity
\(d_\Sigma(z)=|y|\) for \(z=(x,y)\) to reduce the computation of
\(C_H(n,m)\) to the Hardy inequality on \(\mathbb Z^m\setminus\{0\}\).
More precisely, Proposition~\ref{prop:hardy-reduction} shows that
\begin{equation}
  C_H(n,m)=c_m:=
  \inf_{\substack{a\in C_c(\mathbb Z^m)\setminus\{0\}\\a(0)=0}}
  \frac{\displaystyle\sum_{j=1}^m\sum_{y\in\mathbb Z^m}
  |a(y)-a(y-e_j)|^2}
  {\displaystyle\sum_{y\in\mathbb Z^m\setminus\{0\}}
  \dfrac{|a(y)|^2}{|y|^2}}.
\end{equation}
In the sequel, we use \(c_m\) in place of \(C_H(n,m)\).

The key step is to prove the following criterion:
\[
  \overline\kappa_{m,n}>0
  \quad\Longleftrightarrow\quad
  c_m>\left(\frac{m-2}{2}\right)^2
  -\left(\frac{n}{n+m}\right)^2.
\]
Combining this criterion with the known values of \(c_m\) for
\(1\leq m\leq8\) and suitable estimates for \(m\geq9\), we obtain
the dimension classification in Theorem~\ref{thm:main}.

 To prove the sufficiency of this criterion, we use a positive
supersolution of the corresponding equation and the ground-state
transformation. See, for example,
\cite[Proposition~4.8]{KellerPinchoverPogorzelski2020},
\cite[Theorem~1]{KellerPinchoverPogorzelski2018} and
\cite[Theorem~1]{KrejcirikStampach2022}.
More precisely, the left-hand side of \eqref{eq:discrete-HSM}
can be written as
\[
  \begin{aligned}
    &\sum_{z\in\mathbb Z^N}|Du(z)|^2
    -c_m\sum_{z\in\mathbb Z^N\setminus\Sigma}
    \frac{|u(z)|^2}{d_\Sigma(z)^2}\\
    &=\sum_{z\in\mathbb Z^N\setminus\Sigma}
    u(z)\left(\Delta_{\mathrm d}u(z)
    -\frac{c_m}{d_\Sigma(z)^2}u(z)\right),
    \qquad u\in C_c(\mathbb Z^N\setminus\Sigma),
  \end{aligned}
\]
where \(\Delta_{\mathrm d}u(z):=
\sum_{j=1}^N(2u(z)-u(z+e_j)-u(z-e_j))\).
We therefore seek a function \(h\) with \(h>0\) on
\(\mathbb Z^N\setminus\Sigma\) and \(h|_\Sigma=0\) satisfying
\[
  \Delta_{\mathrm d}h(z)-\frac{c_m}{d_\Sigma(z)^2}h(z)\geq0,
  \qquad z\in\mathbb Z^N\setminus\Sigma,
\]
and set \(u=hv\).

In order to apply the supersolution method to prove this criterion, we consider the cases \(m=2\) and \(m\ne2\) separately.
When \(m=2\), we have \(c_2=0\), so the positivity of
\(\overline\kappa_{m,n}\) follows directly from the critical
Sobolev inequality on \(\mathbb Z^N\).
For \(m\ne2\), we use the ground-state transformation and \(Q_1\) interpolation
\cite[Section~2.2, item~(1)]{OrtnerShapeev2012}
(see Subsection~\ref{sec:Q1-interpolation}).
 For \(m=1\), the function \(h(k)=\sqrt{k}\), with \(h(0)=0\), satisfies
\[
   \Delta_{\mathrm d}h(k)-\frac{h(k)}{4k^2}\geq0,
  \qquad k\geq1,
\]
 as in \cite[proof of Theorem~1]{KellerPinchoverPogorzelski2018}. We set
\[
  u(x,\sigma k)=\sqrt{k}\,v_\sigma(x,k),
  \qquad x\in\mathbb Z^n,\quad k\geq1,\quad\sigma\in\{-1,1\}.
\]
Combining the growth of \(h\) with \(Q_1\) interpolation and the weighted Sobolev inequality in Lemma~\ref{lem:weighted-sobolev}, we obtain \(\overline\kappa_{1,n}>0\).

The main difficulty in proving Theorem~\ref{thm:main} lies in
obtaining the required decay estimates when \(m\geq3\).
 When \(m=1\), the interpolation argument in
Lemma~\ref{lem:weighted-Q1-common}
uses the growth of the explicit supersolution \(h(k)=\sqrt{k}\).
 For \(m\geq3\),
\cite[Theorem~4.2]{KellerPinchoverPogorzelski2020} guarantees the existence
of a positive supersolution \(h_m\) but does not provide the bounds needed
 to apply the same interpolation argument, namely
\[
  c|y|^{-\alpha_m}\leq h_m(y)\leq C|y|^{-\alpha_m}
  \qquad\text{for all sufficiently large }|y|,
\]
 where, for \(m\ne2\),
\[
  \alpha_m:=\frac{m-2}{2}
  +\sqrt{\left(\frac{m-2}{2}\right)^2-c_m}.
\]
To overcome this difficulty, we combine ideas from the exhaustion
argument in \cite[Lemma~3.8]{HaeselerKeller2010} and the comparison
arguments using subsolutions and supersolutions in
\cite[Theorem~8 and Section~2.4]{JexStampach2025} to construct a
positive supersolution with the required decay estimates
(see Lemmas~\ref{lem:capacitary-supersolution} and~\ref{lem:hm-decay}).
When the condition on \(c_m\) fails, we construct a sequence of test functions \(v_R\) with \(\|v_R\|_{\ell^{2^*}}=1\) for which the left-hand side of \eqref{eq:discrete-HSM} tends to zero as \(R\to\infty\).

Concerning (B), the case \(m=2\) follows from \(c_2=0\) and
Lemma~\ref{lem:general-sobolev-upper}.
To prove  \(\overline\kappa_{m,n}<\overline S_N\)
when \(m\ne2\) and \(N\geq4\), we construct test functions from an extremal
function for \(\overline S_N\).
In the continuous setting, Tertikas and Tintarev
\cite[proof of Lemma~2.2]{TertikasTintarev2007} used scaling and a cutoff
of an extremal function for the Sobolev inequality on \(\mathbb R^N\). However, this argument cannot be
applied directly to our problem, since scaling is not available on
\(\mathbb Z^N\). To overcome this difficulty, we adapt their cutoff
construction to an extremal function for \(\overline S_N\) and
combine it with translation in the \(y\)-direction.

For \((m,n)=(1,2)\), we adapt the duality argument of Benguria,
Frank and Loss \cite[proof of Theorem~1.1]{BFL2008} to prove
\(\overline\kappa_{1,2}=\overline S_3\).
The key step is to compare the Green function of
\(\Delta_{\mathrm d}-1/(4k^2)\) on \(\mathbb Z^2\times\mathbb N\)
with that of \(\Delta_{\mathrm d}\) on \(\mathbb Z^3\), as stated in
Proposition~\ref{prop:green-order}.

As for (C), we prove nonattainment for \(m=2\) and \((m,n)=(1,2)\).
Since every nonnegative extremal function for \(\overline S_{n+2}\)
is strictly positive, we obtain Corollary~\ref{cor:nonattainment-m2}.
For \((m,n)=(1,2)\), the strict comparison of Green functions in
Proposition~\ref{prop:green-order}, together with the duality argument above,
gives Theorem~\ref{thm:attainment}.
\par

The rest of the paper is organized as follows. In
 Section~\ref{sec:preliminaries}, we introduce the functional setting and notation,
prove \(C_H(n,m)=c_m\), and establish basic estimates for \(Q_1\) interpolants.
Section~\ref{sec:positivity-criterion} is devoted to the proof of
Theorem~\ref{thm:main}.
In Section~\ref{sec:sobolev-comparison}, we
prove Theorem~\ref{thm:sobolev-comparison}, including the Green function
comparison for \((m,n)=(1,2)\).
Section~\ref{sec:attainment} contains the proofs of
Corollary~\ref{cor:nonattainment-m2} and Theorem~\ref{thm:attainment}. Finally, in Appendix~\ref{sec:nine-estimates}, we establish the estimates for \(c_9\)
needed to complete the proof of Theorem~\ref{thm:main}.

Throughout the paper, we use the following conventions:
\begin{itemize}[leftmargin=1.6em,itemsep=0pt,topsep=3pt]
\item \(C\) and \(C_i\) denote positive constants which may vary
from line to line.
\item Dependence on the fixed dimensions is omitted unless it
is relevant to the estimate.
\item \(\mathbf1_A\) denotes the characteristic function of a set \(A\).
\end{itemize}

\raggedbottom
\section{Preliminaries}\label{sec:preliminaries}

In this section, we introduce some notation and prove several preliminary results.

\subsection[The lattice graph and basic notation]{The lattice graph and basic notation}
For \(R>0\), we write \(B_R:=\{z\in\mathbb Z^d:|z|\leq R\}\),
where the dimension \(d\) is determined by the context.
For \(d\geq1\), let \(G=(V,E)\) be the \(d\)-dimensional lattice graph, where
\[
  V:=\mathbb Z^d,\qquad
  E:=\left\{(z,w)\in\mathbb Z^d\times\mathbb Z^d:
  \sum_{j=1}^d|z_j-w_j|=1\right\}.
\]
Two vertices \(z,w\in\mathbb Z^d\) are called neighbors, denoted by
\(z\sim w\), if \((z,w)\in E\).  The combinatorial distance between \(z\) and \(w\) is defined by
\[
  d(z,w):=\inf\{k:z=z_0\sim z_1\sim\cdots\sim z_k=w\}.
\]

 Let
\(C(\mathbb Z^d)\) be the space of real-valued functions on
\(\mathbb Z^d\), and let \(C_c(\mathbb Z^d)\) be the subspace of
functions with finite support.
For \(u\in C(\mathbb Z^d)\), the \(\ell^p\)-norm of \(u\) is defined by
\[
  \|u\|_{\ell^p(\mathbb Z^d)}:=
  \begin{cases}
    \left(\displaystyle\sum_{z\in\mathbb Z^d}|u(z)|^p\right)^{1/p},
    &1\leq p<\infty,\\
    \displaystyle\sup_{z\in\mathbb Z^d}|u(z)|,&p=\infty
  \end{cases}
  .
\]
 When no confusion can arise, we write \(\|u\|_p\) instead of
\(\|u\|_{\ell^p(\mathbb Z^d)}\).
For \(u\in C(\mathbb Z^d)\), define
\[
  D_ju(z):=u(z)-u(z-e_j)\quad(1\leq j\leq d),\qquad
  Du(z):=(D_1u(z),\ldots,D_du(z)).
\]
For \(1\leq p<\infty\), let \(D^{1,p}(\mathbb Z^d)\) be the completion of
\(C_c(\mathbb Z^d)\) under the norm
\begin{equation}
  \label{eq:dirichlet-energy-convention}
  \|u\|_{D^{1,p}(\mathbb Z^d)}:=
  \left(\frac12\sum_{z\in\mathbb Z^d}\sum_{w\sim z}|u(w)-u(z)|^p\right)^{1/p}.
\end{equation}
With this convention, we define the discrete Laplacian  of \(u\in C(\mathbb Z^d)\) by
\[
  \Delta_{\mathrm d}u(z)=\sum_{j=1}^d
  \bigl(2u(z)-u(z+e_j)-u(z-e_j)\bigr).
\]
Before applying \(\Delta_{\mathrm d}\), we extend functions on
\(\mathbb Z^N\setminus\Sigma\) or \(\mathbb Z^m\setminus\{0\}\)
by zero on \(\Sigma\) or at \(0\), respectively.

To study the attainment of \(\overline\kappa_{m,n}\) when \(\overline\kappa_{m,n}>0\),
we define \(\mathcal H_{m,n}\) as the completion of \(C_c(\mathbb Z^N\setminus\Sigma)\)
under the norm
\[
  \|u\|_{\mathcal H_{m,n}}
  :=\left(\sum_{z\in\mathbb Z^N}|Du(z)|^2
  -C_H(n,m)\sum_{z\in\mathbb Z^N\setminus\Sigma}
  \frac{|u(z)|^2}{d_\Sigma(z)^2}\right)^{1/2}.
\]
The next lemma shows that this Hilbert space embeds continuously into
\(\ell^{2^*}(\mathbb Z^N)\).

\begin{lemma}
\label{lem:Hmn-embedding}
Suppose that \(\overline\kappa_{m,n}>0\). For \(U\in\mathcal H_{m,n}\),
choose \(u_j\in C_c(\mathbb Z^N\setminus\Sigma)\) with
\(u_j\to U\) in \(\mathcal H_{m,n}\), and define
\[
  (JU)(z):=
  \begin{cases}
    \displaystyle\lim_{j\to\infty}u_j(z), & z\in\mathbb Z^N\setminus\Sigma,\\
    0, & z\in\Sigma.
  \end{cases}
\]
Then \(J:\mathcal H_{m,n}\to\ell^{2^*}(\mathbb Z^N)\) is a continuous embedding.
\end{lemma}
\begin{proof}
We first show that \(J:\mathcal H_{m,n}\to\ell^{2^*}(\mathbb Z^N)\) is a
well-defined continuous linear map. Let \(U\in\mathcal H_{m,n}\) and choose
\(u_j\in C_c(\mathbb Z^N\setminus\Sigma)\) such that \(u_j\to U\) in \(\mathcal H_{m,n}\).
By \eqref{eq:discrete-HSM}, we have
\[
  \|u_j-u_k\|_{\ell^{2^*}}
  \leq\overline\kappa_{m,n}^{-1/2}\|u_j-u_k\|_{\mathcal H_{m,n}}
  \longrightarrow 0\qquad\text{as }j,k\to\infty.
\]
By the completeness of \(\ell^{2^*}(\mathbb Z^N)\), there exists
\(u\in\ell^{2^*}(\mathbb Z^N)\) such that \(u_j\to u\) in \(\ell^{2^*}(\mathbb Z^N)\).
For each \(z\in\mathbb Z^N\), we have
\[
  |u_j(z)-u(z)|\leq\|u_j-u\|_{\ell^{2^*}}\longrightarrow 0.
\]
In particular, \(u|_\Sigma=0\), since \(u_j|_\Sigma=0\) for every \(j\).
Thus \(JU=u\) by the definition of \(J\).
If \(v_j\in C_c(\mathbb Z^N\setminus\Sigma)\) also converges to \(U\) in
\(\mathcal H_{m,n}\), then \eqref{eq:discrete-HSM} gives
\[
  \|u_j-v_j\|_{\ell^{2^*}}
  \leq\overline\kappa_{m,n}^{-1/2}\|u_j-v_j\|_{\mathcal H_{m,n}}\longrightarrow 0.
\]
Hence \(JU\) is independent of the choice of \((u_j)\).
Linearity follows from the definition of \(J\).
Applying \eqref{eq:discrete-HSM} to \(u_j\) and letting \(j\to\infty\), we obtain
\[
  \|JU\|_{\ell^{2^*}}\leq\overline\kappa_{m,n}^{-1/2}\|U\|_{\mathcal H_{m,n}},
\]
so \(J\) is continuous.

Next, we prove that \(J\) is injective. Suppose that \(JU=0\).
Then \(u_j\to0\) in \(\ell^{2^*}(\mathbb Z^N)\), and hence \(u_j(z)\to0\)
for every \(z\in\mathbb Z^N\). For each \(\varphi\in C_c(\mathbb Z^N\setminus\Sigma)\),
summation by parts gives
\[
  \langle U,\varphi\rangle_{\mathcal H_{m,n}}
  =\lim_{j\to\infty}\sum_{z\in\mathbb Z^N\setminus\Sigma}
  u_j(z)\left(\Delta_{\mathrm d}\varphi(z)
  -\frac{C_H(n,m)}{d_\Sigma(z)^2}\varphi(z)\right)=0,
\]
since \(\varphi\) and \(\Delta_{\mathrm d}\varphi\) have finite support.
Taking \(\varphi=u_k\) and letting \(k\to\infty\), we obtain
\(\|U\|_{\mathcal H_{m,n}}^2=0\) and hence \(U=0\). This completes the proof.
\end{proof}

\subsection[The best constant for the Hardy inequality on ZN minus Sigma]{The best constant for the Hardy inequality on \(\mathbb
Z^N\setminus\Sigma\)}
 To determine when \(\overline\kappa_{m,n}>0\), we first study the best
constant \(C_H(n,m)\) in the Hardy inequality~\eqref{eq:hardy-sigma}. More precisely, we prove
\(C_H(n,m)=c_m\), where \(c_m\) is the best constant in the Hardy inequality on
\(\mathbb Z^m\setminus\{0\}\).

   In order to prove \(C_H(n,m)=c_m\), we first introduce the following cutoff functions:
\[
  \phi_R(x)=
  \begin{cases}
    1,&\|x\|_\infty\leq R,\\[1mm]
    2-\dfrac{\|x\|_\infty}{R},&R<\|x\|_\infty<2R,\\[1mm]
    0,&\|x\|_\infty\geq2R,
  \end{cases}
  \qquad x\in\mathbb Z^n,\quad R\in\mathbb N,\ R\geq1,
\]
where \(\|x\|_\infty=\max_{1\leq j\leq n}|x_j|\).

\begin{lemma}
\label{lem:standard-x-cutoff}
Let \(n\geq1\), \(R\geq1\), and \(1\leq q<\infty\).
Then \(\phi_R\)  satisfies
\[
  \sum_{x\in\mathbb Z^n}|D\phi_R(x)|^2\leq C(n)R^{n-2},
  \qquad
  c(n,q)R^n\leq\sum_{x\in\mathbb Z^n}|\phi_R(x)|^q
  \leq C(n,q)R^n.
\]

\end{lemma}
\begin{proof}
By the definition of \(\phi_R\), we have \(\phi_R(x)=1\)
for \(\|x\|_\infty\leq R\), \(0\leq\phi_R\leq1\), and
\(\phi_R(x)=0\) for \(\|x\|_\infty\geq2R\). Hence,
\[
  (2R+1)^n\leq\sum_{x\in\mathbb Z^n}|\phi_R(x)|^q
  \leq(4R+1)^n,\qquad 1\leq q<\infty.
\]
Since \(R\geq1\), this proves
\begin{equation}
  \label{eq:x-cutoff-Lq}
  c(n)R^n\leq\sum_{x\in\mathbb Z^n}|\phi_R(x)|^q\leq C(n)R^n.
\end{equation}

By the definition of \(\phi_R\), we also have \(|\phi_R(x)-\phi_R(x-e_j)|\leq R^{-1}\)
for \(x\in\mathbb Z^n\) and \(1\leq j\leq n\).
Since \(\phi_R(x)=\phi_R(x-e_j)=0\) for \(\|x\|_\infty>2R\), we have
\begin{equation}
  \label{eq:x-cutoff-gradient}
  \begin{aligned}
    \sum_{x\in\mathbb Z^n}|D\phi_R(x)|^2
    &=\sum_{j=1}^n\sum_{\|x\|_\infty\leq2R}
    |\phi_R(x)-\phi_R(x-e_j)|^2\\
    &\leq\frac{n(4R+1)^n}{R^2}\leq C(n)R^{n-2}.
  \end{aligned}
\end{equation}

By \eqref{eq:x-cutoff-gradient} and \eqref{eq:x-cutoff-Lq},
with \(q=2\), we obtain
\[
  \frac{\displaystyle\sum_{x\in\mathbb Z^n}|D\phi_R(x)|^2}
  {\displaystyle\sum_{x\in\mathbb Z^n}|\phi_R(x)|^2}
  \leq\frac{C(n)}{R^2}\rightarrow0
  \qquad\text{as }R\rightarrow\infty.
\]
This completes the proof.
\end{proof}

With the estimates in Lemma~\ref{lem:standard-x-cutoff}, we  now prove  \(C_H(n,m)=c_m\).

\begin{proposition}
\label{prop:hardy-reduction}
If \(n,m\geq1\) and \(N=n+m\geq3\), then \(C_H(n,m)=c_m\).
\end{proposition}
\begin{proof}
We first prove  \(C_H(n,m)\geq c_m\). Take an arbitrary
\(u\in C_c(\mathbb Z^N\setminus\Sigma)\setminus\{0\}\).
For each \(x\in\mathbb Z^n\), we have \(u(x,\cdot)\in C_c(\mathbb Z^m)\)
and \(u(x,0)=0\). By the definition of \(c_m\), we obtain
\[
  \begin{aligned}
    \sum_{z\in\mathbb Z^N}|Du(z)|^2
    &=\sum_{j=1}^n\sum_{x\in\mathbb Z^n}\sum_{y\in\mathbb Z^m}
    |u(x,y)-u(x-e_j,y)|^2\\
    &\quad+\sum_{j=1}^m\sum_{x\in\mathbb Z^n}\sum_{y\in\mathbb Z^m}
    |u(x,y)-u(x,y-e_j)|^2\\
    &\geq c_m\sum_{x\in\mathbb Z^n}
    \sum_{y\in\mathbb Z^m\setminus\{0\}}\frac{|u(x,y)|^2}{|y|^2}.
  \end{aligned}
\]
 Taking the infimum over all such \(u\), we obtain \(C_H(n,m)\geq c_m\).

For the reverse inequality, fix
\(0\neq a\in C_c(\mathbb Z^m)\)
with \(a(0)=0\), and set \(u_R(x,y)=\phi_R(x)a(y)\).
Then \(u_R\in C_c(\mathbb Z^N\setminus\Sigma)\setminus\{0\}\). By the definition of \(u_R\), we have
\[
  \begin{aligned}
    D_j u_R(x,y)&=a(y)D_j\phi_R(x),&&1\leq j\leq n,\\
    D_{n+j}u_R(x,y)&=\phi_R(x)D_j a(y),&&1\leq j\leq m.
  \end{aligned}
\]
Consequently,
\begin{equation}
  \label{eq:product-energy-hardy}
  \begin{aligned}
    \sum_{z\in\mathbb Z^N}|Du_R(z)|^2
    &=\left(\sum_{y\in\mathbb Z^m}|a(y)|^2\right)
    \left(\sum_{x\in\mathbb Z^n}|D\phi_R(x)|^2\right)\\
    &\quad+\left(\sum_{x\in\mathbb Z^n}|\phi_R(x)|^2\right)
    \left(\sum_{y\in\mathbb Z^m}|Da(y)|^2\right),\\
    \sum_{x\in\mathbb Z^n}\sum_{y\in\mathbb Z^m\setminus\{0\}}
    \frac{|u_R(x,y)|^2}{|y|^2}
    &=\left(\sum_{x\in\mathbb Z^n}|\phi_R(x)|^2\right)
    \left(\sum_{y\in\mathbb Z^m\setminus\{0\}}
    \frac{|a(y)|^2}{|y|^2}\right).
  \end{aligned}
\end{equation}
Hence,
\[
  \begin{aligned}
    C_H(n,m)
    &\leq
    \frac{\displaystyle\sum_{z\in\mathbb Z^N}|Du_R(z)|^2}
    {\displaystyle\sum_{x\in\mathbb Z^n}
    \sum_{y\in\mathbb Z^m\setminus\{0\}}
    \frac{|u_R(x,y)|^2}{|y|^2}}\\
    &=
    \frac{\displaystyle\sum_{y\in\mathbb Z^m}|Da(y)|^2}
    {\displaystyle\sum_{y\in\mathbb Z^m\setminus\{0\}}
    \frac{|a(y)|^2}{|y|^2}}
    +
    \frac{\displaystyle\sum_{y\in\mathbb Z^m}|a(y)|^2}
    {\displaystyle\sum_{y\in\mathbb Z^m\setminus\{0\}}
    \frac{|a(y)|^2}{|y|^2}}
    \frac{\displaystyle\sum_{x\in\mathbb Z^n}|D\phi_R(x)|^2}
    {\displaystyle\sum_{x\in\mathbb Z^n}|\phi_R(x)|^2}.
  \end{aligned}
\]
  Letting \(R\to\infty\) with \(a\) fixed, we deduce from
Lemma~\ref{lem:standard-x-cutoff} that
\[
  C_H(n,m)\leq
  \frac{\displaystyle\sum_{y\in\mathbb Z^m}|Da(y)|^2}
  {\displaystyle\sum_{y\in\mathbb Z^m\setminus\{0\}}
  \frac{|a(y)|^2}{|y|^2}}.
\]
 Taking the infimum over all such \(a\), we obtain \(C_H(n,m)\leq c_m\).
This completes the proof.
\end{proof}

\subsection[Q1 interpolation]{\(Q_1\) interpolation}
\label{sec:Q1-interpolation}

As mentioned in the Introduction, we use \(Q_1\) interpolation to prove the
weighted Sobolev inequalities in Section~\ref{sec:positivity-criterion}.
We first recall its definition; see
\cite[Section~2.1, formula~(1), and Section~2.2, item~(1)]{OrtnerShapeev2012}.

\begin{definition}
Let \(d\geq1\) and let \(f:\mathbb Z^d\rightarrow\mathbb R\). For
\(z\in\mathbb Z^d\), set \(Q_z:=z+[0,1]^d\). For
\(\nu\in\{0,1\}^d\) and \(\theta\in[0,1]^d\), set
\[
  \lambda_\nu(\theta)
  :=
  \prod_{j=1}^d
  \theta_j^{\nu_j}(1-\theta_j)^{1-\nu_j}.
\]
On \(Q_z\), define the \(Q_1\) interpolant of \(f\) by
\begin{equation}
  \label{eq:Q1-interpolant}
  If(z+\theta):=\sum_{\nu\in\{0,1\}^d}
  \lambda_\nu(\theta)f(z+\nu).
\end{equation}
\end{definition}

 We next establish some basic estimates for the interpolant \(If\).

\begin{lemma}
\label{lem:Q1-basic-estimates}
Let \(d\geq1\) and \(f\in C_c(\mathbb Z^d)\). For
\(i\in\{1,\ldots,d\}\) and \(\nu\in\{0,1\}^d\) with \(\nu_i=0\), set
\[
  \lambda_\nu^{(i)}(\theta)
  :=
  \prod_{\substack{1\leq j\leq d\\ j\neq i}}
  \theta_j^{\nu_j}(1-\theta_j)^{1-\nu_j},
  \qquad
  \theta\in[0,1]^d.
\]
Then the following estimates hold.
\begin{enumerate}[label=\textup{(\roman*)}]
\item
For every \(z\in\mathbb Z^d\) and \(i\in\{1,\ldots,d\}\),
\begin{equation}
  \label{eq:Q1-local-Jensen}
  |\partial_iIf(z+\theta)|^2
  \leq
  \sum_{\nu_i=0}
  \lambda_\nu^{(i)}(\theta)
  |f(z+\nu+e_i)-f(z+\nu)|^2
\end{equation}
for almost every \(\theta\in[0,1]^d\).

Here and below, \(\sum_{\nu_i=0}\) denotes the sum over all
\(\nu\in\{0,1\}^d\) such that \(\nu_i=0\).

\item
\[
  \int_{\mathbb R^d}|\nabla If(x)|^2\,dx
  \leq
  \sum_{i=1}^d
  \sum_{z\in\mathbb Z^d}
  |f(z+e_i)-f(z)|^2.
\]

\item
For every \(q>0\) and \(z\in\mathbb Z^d\),
\begin{equation}
  \max_{\nu\in\{0,1\}^d}|f(z+\nu)|^q
  \leq
  C(d,q)
  \int_{Q_z}|If(x)|^q\,dx.
\end{equation}
\end{enumerate}
\end{lemma}

\begin{proof}
We first prove \textup{(i)}. Differentiating
\eqref{eq:Q1-interpolant} with respect to \(\theta_i\), we obtain
\begin{equation}
  \label{eq:Q1-local-derivative}
  \partial_iIf(z+\theta)=\sum_{\nu_i=0}
  \lambda_\nu^{(i)}(\theta)
  \bigl(f(z+\nu+e_i)-f(z+\nu)\bigr).
\end{equation}
 Since \(\lambda_\nu^{(i)}\geq0\) and
\(\sum_{\nu_i=0}\lambda_\nu^{(i)}=1\), Jensen's inequality yields \textup{(i)}.

 To prove \textup{(ii)}, we first observe that
\[
  \int_{[0,1]^d}\lambda_\nu^{(i)}(\theta)\,d\theta
  =\prod_{j\neq i}\int_0^1\theta_j^{\nu_j}
  (1-\theta_j)^{1-\nu_j}\,d\theta_j
  =2^{-(d-1)}.
\]

Integrating \eqref{eq:Q1-local-Jensen} over \(Q_z\) and then
summing over \(z\in\mathbb Z^d\), we deduce that
\[
  \begin{aligned}
    \int_{\mathbb R^d}|\partial_iIf(x)|^2\,dx
    &\leq2^{-(d-1)}\sum_{\nu_i=0}\sum_{z\in\mathbb Z^d}
    |f(z+\nu+e_i)-f(z+\nu)|^2\\
    &=\sum_{w\in\mathbb Z^d}|f(w+e_i)-f(w)|^2.
  \end{aligned}
\]
Summing over \(i\), we obtain \textup{(ii)}.

Finally, to prove \textup{(iii)}, let \(P(\theta)=If(z+\theta)\) and
\(M=\max_{\nu\in\{0,1\}^d}|P(\nu)|\), and choose a vertex
\(\nu_0\) with \(|P(\nu_0)|=M\). Since
\(|\partial_iP|\leq2M\) by \eqref{eq:Q1-local-derivative}, for
 every \(\theta\in[0,1]^d\) satisfying
\(|\theta_i-(\nu_0)_i|\leq(4d)^{-1}\) for all \(i=1,\ldots,d\), we have
\[
  |P(\theta)-P(\nu_0)|
  \leq2M\sum_{i=1}^d|\theta_i-(\nu_0)_i|\leq M/2.
\]

Set \(A:=\{\theta\in[0,1]^d:|\theta_i-(\nu_0)_i|\leq(4d)^{-1},\ 1\leq i\leq d\}\).
Then \(|A|=(4d)^{-d}\) and \(|P(\theta)|\geq M/2\) for \(\theta\in A\). Hence,
\[
  \max_{\nu\in\{0,1\}^d}|f(z+\nu)|^q=M^q\leq2^q(4d)^d
  \int_{[0,1]^d}|P(\theta)|^q\,d\theta,
\]
which proves \textup{(iii)} for every \(q>0\).
This completes the proof.
\end{proof}

\section{Proof of Theorem~\ref{thm:main}}
\label{sec:positivity-criterion}

In this section, we prove Theorem~\ref{thm:main}.
To determine the dimensions for which \(\overline\kappa_{m,n}>0\),
we first prove  the following criterion for \(n,m\geq1\) and \(N=n+m\geq3\):
\begin{equation}
  \label{eq:positivity-criterion}
  \overline\kappa_{m,n}>0
  \quad\Longleftrightarrow\quad
  c_m>\left(\frac{m-2}{2}\right)^2
  -\left(\frac{n}{n+m}\right)^2.
\end{equation}
 Combining \eqref{eq:positivity-criterion} with \(c_m=(m-2)^2/4\) for \(1\leq m\leq8\) (see
\cite{AlpayLow2026,Gupta2024,HuangYe2024,Lizama2026}) and with the estimates
for \(m\geq9\), we obtain Theorem~\ref{thm:main}.

To prove \eqref{eq:positivity-criterion}, we consider the case
\(m=2\) and the case \(m\ne2\) separately.

The following lemma will be used in both cases.
\begin{lemma}
\label{lem:general-sobolev-upper}
For \(n,m\geq1\) and \(N=n+m\geq3\), we have
\begin{equation}
  \label{eq:general-sobolev-upper}
  0\leq\overline\kappa_{m,n}\leq\overline S_N,
  \qquad n,m\geq1,\quad N=n+m\geq3.
\end{equation}
\end{lemma}
\begin{proof}
The Hardy inequality gives \(\overline\kappa_{m,n}\geq0\).
For the upper bound, take \(0\neq f\in C_c(\mathbb Z^N)\), choose an integer
\(L>\max\{|y|:(x,y)\in\operatorname{supp}f\}\), and set
\(u_L(x,y)=f(x,y-Le_1)\), where \(e_1\in\mathbb Z^m\).
Then \(u_L\in C_c(\mathbb Z^N\setminus\Sigma)\setminus\{0\}\), and
\[
  \sum_{z\in\mathbb Z^N}|Du_L(z)|^2=\sum_{z\in\mathbb Z^N}|Df(z)|^2,
  \qquad
  \sum_{z\in\mathbb Z^N}|u_L(z)|^{2^*}=\sum_{z\in\mathbb Z^N}|f(z)|^{2^*}.
\]
By \(c_m\geq0\) and the definition of \(\overline\kappa_{m,n}\), we have
\[
  \overline\kappa_{m,n}
  \leq\frac{\displaystyle\sum_{z\in\mathbb Z^N}|Du_L(z)|^2
  -c_m\sum_{x\in\mathbb Z^n}\sum_{y\in\mathbb Z^m\setminus\{0\}}
  |u_L(x,y)|^2/|y|^2}
  {\displaystyle\left(\sum_{z\in\mathbb Z^N}|u_L(z)|^{2^*}\right)^{2/2^*}}
  \leq\frac{\displaystyle\sum_{z\in\mathbb Z^N}|Df(z)|^2}
  {\displaystyle\left(\sum_{z\in\mathbb Z^N}|f(z)|^{2^*}\right)^{2/2^*}}.
\]
Taking the infimum over \(f\in C_c(\mathbb Z^N)\setminus\{0\}\), we obtain
\eqref{eq:general-sobolev-upper}.
This completes the proof.
\end{proof}

\subsection[The case m=2]{The case \(m=2\)}
Notice that when \(m=2\), we have \(c_2=0\), so
\eqref{eq:discrete-HSM} reduces to the
critical Sobolev inequality on
\(C_c(\mathbb Z^{n+2}\setminus\Sigma)\).
In this case, we prove  \(\overline\kappa_{2,n}=\overline S_{n+2}\).

\begin{proposition}
\label{prop:m2}
Let \(m=2\) and \(n\geq1\). Then
\(\overline\kappa_{2,n}=\overline S_N>0\).
\end{proposition}
\begin{proof}
By \cite[Remark~1.2]{Gupta2024} and
\eqref{eq:dirichlet-energy-convention}, we have \(c_2=0\).
For \(0\neq u\in C_c(\mathbb Z^N\setminus\Sigma)\),
by the critical Sobolev inequality on \(\mathbb Z^N\)~\cite[Theorem~1.8]{HeJi2026}, we obtain
\[
  \frac{\displaystyle\sum_{z\in\mathbb Z^N}|Du(z)|^2}
  {\displaystyle\left(\sum_{z\in\mathbb Z^N}|u(z)|^{2^*}\right)^{2/2^*}}
  \geq\overline S_N>0.
\]
Since \(u\in C_c(\mathbb Z^N\setminus\Sigma)\setminus\{0\}\)
is arbitrary, we obtain \(\overline\kappa_{2,n}\geq\overline S_N\).

 Combining this with \eqref{eq:general-sobolev-upper}, we obtain
 \(\overline\kappa_{2,n}=\overline S_{n+2}>0\).
This completes the proof.
\end{proof}

\subsection[The case m not equal to 2]{The case  \(m\ne2\)}

 To prove \eqref{eq:positivity-criterion} for \(m\ne2\), we first construct
positive supersolutions satisfying
\[
  c|y|^{-\alpha_m}\leq h_m(y)\leq C|y|^{-\alpha_m},
  \qquad y\in\mathbb Z^m\setminus\{0\},
\]
where \(\alpha_m\) is defined in Section~\ref{sec:proof-ideas}.
To construct such a supersolution, we consider the operator
\(\mathcal L_m\) defined, for \(f\in C(\mathbb Z^m)\) with \(f(0)=0\), by
\[
  (\mathcal L_mf)(y):=\sum_{j=1}^m
  \bigl(2f(y)-f(y+e_j)-f(y-e_j)\bigr)-\frac{c_m}{|y|^2}f(y),
  \qquad y\in\mathbb Z^m\setminus\{0\}.
\]

For \(m\geq3\), we construct \(h_m\) by solving minimization problems on finite
subsets of \(\mathbb Z^m\setminus\{0\}\).
By \cite[Theorem~3.1]{KapitanskiLaptev2016},
there exists \(C_m>0\) such that
\[
  \sum_{y\in\mathbb Z^m\setminus\{0\}}\frac{|u(y)|^2}{|y|^2}
  \leq C_m\sum_{y\in\mathbb Z^m}|Du(y)|^2,
  \qquad u\in C_c(\mathbb Z^m),\quad u(0)=0.
\]

For this purpose, for \(f\in C_c(\mathbb Z^m)\) with \(f(0)=0\), define
\[
  q_m(f):=\sum_{y\in\mathbb Z^m}|Df(y)|^2
  -c_m\sum_{y\in\mathbb Z^m\setminus\{0\}}\frac{|f(y)|^2}{|y|^2}.
\]
For \(f\in C(\mathbb Z^m)\) and \(\varphi\in C_c(\mathbb Z^m)\)
with \(f(0)=\varphi(0)=0\), define
\[
  q_m(f,\varphi):=\sum_{y\in\mathbb Z^m}Df(y)\cdot D\varphi(y)
  -c_m\sum_{y\in\mathbb Z^m\setminus\{0\}}\frac{f(y)\varphi(y)}{|y|^2}.
\]
Then
\begin{equation}
  \label{eq:qm-green-identity}
  q_m(f,\varphi)=\sum_{y\in{\mathbb Z^m\setminus\{0\}}}
  (\mathcal L_mf)(y)\varphi(y).
\end{equation}
By the definition of \(c_m\), we have \(q_m(f)\geq0\).
To prove that the minimization problem
\[
  \inf_{\substack{
  f\in C_c(\mathbb Z^m\setminus\{0\}),\ f(e_1)=1\\
  \operatorname{supp}f\subset\{y\in\mathbb Z^m:0<|y|<R\}}}q_m(f)
\]
has a unique solution for each integer \(R\geq3\), we first prove
\(q_m(f)>0\) for every \(0\neq f\in C_c(\mathbb Z^m\setminus\{0\})\).

\begin{lemma}
\label{lem:qm-strict-positivity}
Let \(m\geq3\). Then
\(q_m(f)>0\) for every
\(0\neq f\in C_c(\mathbb Z^m)\) with \(f(0)=0\).
\end{lemma}
\begin{proof}
We argue by contradiction. Suppose that \(q_m(f)=0\) for some \(0\ne f\in C_c(\mathbb Z^m)\) with \(f(0)=0\). Since \(0\leq q_m(|f|)\leq q_m(f)=0\), we may assume that \(f\geq0\).
Since \(\operatorname{supp}f\) is finite and \(\mathbb Z^m\setminus\{0\}\) is connected, there exist \(y,z\ne0\) such that \(y\sim z\), \(f(y)=0\), and \(f(z)>0\). By \eqref{eq:qm-green-identity}, we have
\[
 q_m(f,\delta_y)=-\sum_{\substack{\eta\in\mathbb Z^m\setminus\{0\}\\\eta\sim y}}f(\eta)<0.
\]
Therefore, for all sufficiently small \(t>0\),
\[
 q_m(f+t\delta_y)=2tq_m(f,\delta_y)+t^2q_m(\delta_y)<0,
\]
which contradicts \(q_m(f+t\delta_y)\geq0\). This completes the proof.
\end{proof}

To construct \(h_m\) on finite sets and compare it with subsolutions and
supersolutions, we need the following comparison principle.
For a finite set \(D\subset\mathbb Z^m\setminus\{0\}\), set
\[
  \partial_{\mathrm V}D:=\{z\in\mathbb Z^m\setminus(D\cup\{0\}):
  z\sim y\text{ for some }y\in D\}.
\]

\begin{lemma}
\label{lem:finite-domain-comparison}
Assume that \(m\geq3\) and \(D\subset\mathbb Z^m\setminus\{0\}\) is finite.
If \(u,v\in C(\mathbb Z^m)\) satisfy \(u(0)=v(0)=0\), \(\mathcal L_mu\leq\mathcal L_mv\) in \(D\) and
\(u\leq v\) on \(\partial_{\mathrm V}D\), then \(u\leq v\) in \(D\).
\end{lemma}
\begin{proof}
Set \(w=u-v\) and define
\[
  w_+(y):=\begin{cases}
    \max\{w(y),0\},&y\in D,\\
    0,&\text{else}.
  \end{cases}
\] By \(w\leq0\) on
\(\partial_{\mathrm V}D\) and
\((a-b)(a_+-b_+)\geq|a_+-b_+|^2\), we have
\[
  q_m(w_+)\leq q_m(w,w_+)
  =\sum_{y\in D}(\mathcal L_mw)(y)w_+(y)\leq0.
\]
By Lemma~\ref{lem:qm-strict-positivity}, we obtain \(w_+=0\).
This completes the proof.
\end{proof}

To construct \(h_m\), we first study the following minimization problem:
\[
  \inf_{f\in\mathcal A_R}q_m(f),\qquad
  \mathcal A_R:=\{f\in C_c(\mathbb Z^m):
  \operatorname{supp}f\subset B_R^*,\ f(o)=1\},
\]
where \(o=e_1\), \(B_R^*:=\{y\in\mathbb Z^m:0<|y|<R\}\), and
\(R\geq3\) is an integer. The following lemma gives its unique minimizer
\(h_R\) and proves
\[
  0<h_m(y):=\lim_{\substack{R\rightarrow\infty\\R\in\mathbb N}}h_R(y)<\infty,
  \qquad y\in\mathbb Z^m\setminus\{0\}.
\]

\begin{lemma}
\label{lem:capacitary-supersolution}
Let \(m\geq3\).
\begin{enumerate}[label=\textup{(\roman*)}]
\item For every integer \(R\geq3\), there exists a unique \(h_R\in\mathcal A_R\) such that
\[
  q_m(h_R)=\inf_{f\in\mathcal A_R}q_m(f)=:\mu_R.
\]
Moreover,
\[
  h_R>0\text{ in }B_R^*,\qquad \mu_R>0,\qquad
  \mathcal L_mh_R=\mu_R\delta_o\text{ in }B_R^*.
\]
\item For integers \(3\leq R<S\), we have \(h_R\leq h_S\) on
\(\mathbb Z^m\) and \(0<\mu_S\leq\mu_R\).
\item For every \(y\in{\mathbb Z^m\setminus\{0\}}\),
\[
  0<h_m(y):=\lim_{\substack{R\rightarrow\infty\\R\in\mathbb N}}h_R(y)<\infty.
\]
With \(h_m(0)=0\), we have
\[
  \mu_R\rightarrow\mu_m\geq0,\qquad h_m(o)=1,\qquad
  \mathcal L_mh_m=\mu_m\delta_o\quad\text{on }\mathbb Z^m\setminus\{0\}.
\]
\end{enumerate}
\end{lemma}
\begin{proof}
We first prove \textup{(i)}. By Lemma~\ref{lem:qm-strict-positivity},
\(q_m\) is positive definite on
\[
  \mathcal V_R:=\{f\in C_c(\mathbb Z^m):\operatorname{supp}f\subset B_R^*\}.
\]
Since \(\dim\mathcal V_R=\#B_R^*<\infty\), there exists \(a_R>0\) such that
\(q_m(f)\geq a_R\sum_{y\in B_R^*}|f(y)|^2\).
Since \(\mathcal A_R\) is nonempty, closed and convex, and \(q_m\) is strictly convex,
the preceding lower bound gives a unique minimizer \(h_R\) on \(\mathcal A_R\). Since \(|h_R|\in\mathcal A_R\) and
\(q_m(|h_R|)\leq q_m(h_R)\), uniqueness gives \(h_R\geq0\).
For
\(\varphi\) supported in \(B_R^*\setminus\{o\}\),
\(h_R+t\varphi\in\mathcal A_R\), and hence
\[
  0=\left.\frac d{dt}q_m(h_R+t\varphi)\right|_{t=0}
  =2q_m(h_R,\varphi)
  =2\sum_{y\in B_R^*}(\mathcal L_mh_R)(y)\varphi(y).
\]
Taking \(\varphi=\delta_y\) gives
\begin{equation}
  \label{eq:hR-harmonic}
  \mathcal L_mh_R=0\text{ in }B_R^*\setminus\{o\},
  \qquad h_R(o)=1,\qquad h_R=0\text{ on }\mathbb Z^m\setminus B_R^*.
\end{equation}
To prove \(h_R>0\), we argue by contradiction. Suppose that
\(h_R(y)=0\) at some \(y\in B_R^*\). Then
\eqref{eq:hR-harmonic} and \(h_R\geq0\) give \(h_R=0\) at every neighbor.
Since \(B_R^*\) is connected, repeating this argument along a path to \(o\) contradicts
\(h_R(o)=1\). Hence, \(h_R>0\) in \(B_R^*\). By
\eqref{eq:qm-green-identity}, we have
\[
  \mu_R=q_m(h_R)
  =\sum_{y\in B_R^*}(\mathcal L_mh_R)(y)h_R(y)
  =(\mathcal L_mh_R)(o)>0,
\]
where \(q_m(h_R)>0\) follows from Lemma~\ref{lem:qm-strict-positivity}.

Next, we prove \textup{(ii)}. Let \(S>R\). By \eqref{eq:hR-harmonic}, we have
\[
  \mathcal L_mh_R=\mathcal L_mh_S=0
  \quad\text{in }B_R^*\setminus\{o\}.
\]
Moreover,
\[
  h_R(o)=h_S(o)=1,\qquad
  h_R(y)=0\leq h_S(y)
  \quad\text{for }y\in\mathbb Z^m\setminus B_R^*.
\]
Applying Lemma~\ref{lem:finite-domain-comparison} to
\(B_R^*\setminus\{o\}\), we obtain \(h_R\leq h_S\) on \(\mathbb Z^m\).
Since \(\mathcal A_R\subset\mathcal A_S\), we also have \(0<\mu_S\leq\mu_R\).

Finally, we prove \textup{(iii)}. By \textup{(ii)}, \(h_R(y)\) is nondecreasing in \(R\),
so we first show that
\[
  \sup_{\substack{R\geq3\\R\in\mathbb N}}h_R(y)<\infty
  \qquad\text{for every }y\in\mathbb Z^m\setminus\{0\}.
\]
By \(\mathcal L_mh_R=\mu_R\delta_o\), we have
\[
  \sum_{\substack{z\neq0\\z\sim o}}h_R(z)<2m,
  \qquad\sum_{\substack{z\neq0\\z\sim y}}h_R(z)
  =\left(2m-\frac{c_m}{|y|^2}\right)h_R(y)\leq2mh_R(y)
  \quad(y\in B_R^*\setminus\{o\}).
\]
Fix a path \(o=y_0\sim y_1\sim\cdots\sim y_\ell=y\) in
\({\mathbb Z^m\setminus\{0\}}\). For \(R\) large enough to contain
this path, we deduce that \(h_R(y)\leq(2m)^\ell h_R(o)=(2m)^\ell\).
By \textup{(ii)}, the limit \(h_m(y)\) exists and satisfies
\(h_m(y)\geq h_{R_0}(y)>0\) for \(y\in B_{R_0}^*\).
Letting \(R\rightarrow\infty\) in \(\mathcal L_mh_R=\mu_R\delta_o\), we obtain
\textup{(iii)}. This limit is valid since each vertex has finitely many neighbors.
This completes the proof.
\end{proof}

To  show that \(\alpha_m\) is well defined, we first prove
\(c_m\leq(m-2)^2/4\) for \(m\geq3\). For this purpose, we compare
\(c_m\) with the following constant:
\[
  I_R:=\inf_{\substack{0\neq u\in C_c(\mathbb Z^m)\\
  \operatorname{supp}u\subset\{y\in\mathbb Z^m:|y|>R\}}}
  \frac{\displaystyle\sum_{y\in\mathbb Z^m}|Du(y)|^2}
  {\displaystyle\sum_{y\in\mathbb Z^m\setminus\{0\}}|u(y)|^2/|y|^2},\qquad R>0.
\]

We now prove  \(\sup_{R>0}I_R=(m-2)^2/4\).

\begin{proposition}
\label{prop:euclidean-hardy-comparison}
Let \(m\geq3\). Then
\(c_m\leq\sup_{R>0}I_R=\left(\frac{m-2}{2}\right)^2\).
\end{proposition}
\begin{proof}
By the definitions of \(c_m\) and \(I_R\), we have \(c_m\leq I_R\) for every \(R>0\). Let \(u\) be a nonzero complex-valued function with finite support in
\(\{y\in\mathbb Z^m:|y|>R\}\), and set \(v(y)=|u(y)|\).
Then \(\operatorname{supp}v=\operatorname{supp}u\).
By \(|v(y)-v(z)|\leq|u(y)-u(z)|\) and \(v(y)^2=|u(y)|^2\), we have
\[
  \frac{\displaystyle\sum_{y\in\mathbb Z^m}|Dv(y)|^2}
       {\displaystyle\sum_{y\in\mathbb Z^m\setminus\{0\}}v(y)^2/|y|^2}
  \leq
  \frac{\displaystyle\sum_{y\in\mathbb Z^m}|Du(y)|^2}
       {\displaystyle\sum_{y\in\mathbb Z^m\setminus\{0\}}|u(y)|^2/|y|^2}.
\]
Hence \(I_R=C_R^\infty(m)\), where \(C_R^\infty(m)\) is the infimum
over complex-valued functions in \cite[Section~3.1.3]{LizamaTransition2026}.
For the case \(N=m\), \cite[Proposition~3.3]{LizamaTransition2026} gives
\(\lim_{R\to\infty}I_R=(m-2)^2/4\).
Since \(I_R\) is nondecreasing in \(R\), we obtain \(\sup_{R>0}I_R=(m-2)^2/4\). This completes the proof.
\end{proof}

To estimate \(\mathcal L_m\) on the comparison functions in Lemma~\ref{lem:hm-decay}, we apply Taylor's formula.
For \(F\in C^4(\mathbb R^m\setminus\{0\})\) and
\(y\in\mathbb Z^m\) with \(|y|\geq2\), we have
\begin{equation}
  \label{eq:discrete-Taylor-remainder}
  |\Delta_{\mathrm d}F(y)+\Delta F(y)|
  \leq\frac1{12}\sum_{j=1}^m\max_{|t|\leq1}
  |\partial_j^4F(y+te_j)|.
\end{equation}

By \cite[Theorem~1.3]{HuangYe2024} with \(\alpha=0\), we have \(c_1=1/4\).
Together with Proposition~\ref{prop:euclidean-hardy-comparison}, this gives
\(c_m\leq(m-2)^2/4\) for \(m\ne2\), so \(\alpha_m\) is well defined.
In particular, \(\alpha_1=-1/2\).
To apply the ground-state transformation and \(Q_1\) interpolation,
we next construct \(h_m\) satisfying the following bounds with exponent \(\alpha_m\).

\begin{lemma}
\label{lem:hm-decay}
Let \(m=1\) or \(m\geq3\) . There exist a function
\(h_m\in C(\mathbb Z^m)\) and constants \(0<c\leq C<\infty\) such that
\(h_m(0)=0\), \(h_m(e_1)=1\), and
\[
  h_m(y)>0,\qquad (\mathcal L_mh_m)(y)\geq0,
  \qquad y\in\mathbb Z^m\setminus\{0\}.
\]
Moreover,
\begin{equation}
  \label{eq:global-h-comparability}
  c|y|^{-\alpha_m}\leq h_m(y)\leq C|y|^{-\alpha_m},
  \qquad y\in\mathbb Z^m\setminus\{0\}.
\end{equation}
\end{lemma}
\begin{proof}
For \(m=1\), take \(h_1(y)=\sqrt{|y|}\) for \(y\ne0\), with \(h_1(0)=0\),
as in \cite[proof of Theorem~1]{KellerPinchoverPogorzelski2018}.
For \(k\geq1\), we have
\[
  w_k:=\frac{2h_1(k)-h_1(k+1)-h_1(k-1)}{h_1(k)}
  =2-\sqrt{\frac{k+1}{k}}-\sqrt{\frac{k-1}{k}}
  >\frac1{4k^2}.
\]
Since \(c_1=1/4\) and \(h_1(-k)=h_1(k)\), we obtain
\begin{equation}
  \label{eq:one-dimensional-defect}
  h_1(y)(\mathcal L_1h_1)(y)=\rho_{|y|},\qquad
  \rho_k:=k\left(w_k-\frac1{4k^2}\right)>0,\qquad
  \rho_1=\frac74-\sqrt2>0.
\end{equation}
By \(\alpha_1=-1/2\), \eqref{eq:global-h-comparability} holds with \(c=C=1\).

For \(m\geq3\), take the function \(h_m\) from
Lemma~\ref{lem:capacitary-supersolution}. It satisfies
\(h_m(0)=0\), \(h_m(e_1)=1\), \(h_m(y)>0\) for \(y\ne0\), and
\(\mathcal L_mh_m=\mu_m\delta_{e_1}\geq0\).
It remains to prove \eqref{eq:global-h-comparability}.

To compare \(h_m\) with powers of \(|y|\), let
\(\alpha=\alpha_m\) and \(F_\beta(x)=|x|^{-\beta}\). We have
\[
  \begin{aligned}
    |\partial_j^4F_\beta(y+te_j)|&\leq C|y|^{-\beta-4}\qquad(|t|\leq1),\\
    -\Delta F_\beta(y)&=\beta(m-2-\beta)|y|^{-\beta-2}.
  \end{aligned}
\]
By \eqref{eq:discrete-Taylor-remainder}, we obtain, uniformly in
\(y/|y|\) as \(|y|\to\infty\),
\begin{equation}
  \label{eq:general-indicial-expansion}
  \mathcal L_m(|y|^{-\beta})
  =\bigl[\beta(m-2-\beta)-c_m\bigr]|y|^{-\beta-2}
  +O(|y|^{-\beta-4}).
\end{equation}

By the definition of \(\alpha\), we have
\(\alpha(m-2-\alpha)=c_m\) and
\((\alpha+2)(m-4-\alpha)-c_m=-2(2\alpha+4-m)\leq-4\). By \eqref{eq:general-indicial-expansion}, we can choose \(A>0\), then \(R_0>2\), such that
\[
  \underline\Phi(y):=|y|^{-\alpha}+A|y|^{-\alpha-2},
  \qquad\overline\Phi(y):=|y|^{-\alpha}-A|y|^{-\alpha-2}
\]
satisfy, for \(|y|\geq R_0-1\),
\[
  \mathcal L_m\underline\Phi\leq0\leq\mathcal L_m\overline\Phi,
  \qquad\tfrac12|y|^{-\alpha}\leq\overline\Phi(y)
  \leq\underline\Phi(y)\leq2|y|^{-\alpha}.
\]

We first prove \(h_m(y)\leq C|y|^{-\alpha}\). Set \(D_R=\{y\in\mathbb Z^m:R_0<|y|<R\}\).
Define
\[
  \begin{aligned}
    \Gamma_0&:=\{y\in\mathbb Z^m:|y|\leq R_0,\ \exists z\in D_R\text{ such that }y\sim z\},\\
    \Gamma_R&:=\{y\in\mathbb Z^m:|y|\geq R,\ \exists z\in D_R\text{ such that }y\sim z\}.
  \end{aligned}
\]
Since \(\bigl||y|-|z|\bigr|\leq|y-z|=1\) when \(y\sim z\),
\[
  \Gamma_0\subset\{y:R_0-1\leq|y|\leq R_0\},\qquad
  \Gamma_R\subset\{y:R\leq|y|<R+1\}.
\] By Lemma~\ref{lem:capacitary-supersolution}\textup{(iii)}, \(h_R(y)\leq h_m(y)<\infty\). Since
\(\Gamma_0\subset\{y:R_0-1\leq|y|\leq R_0\}\), we may choose
\[
  C=\max_{\substack{y\in\mathbb Z^m\\R_0-1\leq|y|\leq R_0}}\frac{h_m(y)}{\overline\Phi(y)}<\infty,
\]
so that \(h_R\leq C\overline\Phi\) on \(\Gamma_0\), independently of \(R\). On \(\Gamma_R\), we have
\(h_R=0\leq C\overline\Phi\). Applying
Lemma~\ref{lem:finite-domain-comparison} and letting \(R\rightarrow\infty\),
we obtain
\begin{equation}
  \label{eq:hm-upper-bound}
  h_m(y)\leq C\overline\Phi(y)\leq C|y|^{-\alpha},\qquad |y|>R_0.
\end{equation}

Next, to prove \(h_m(y)\geq c|y|^{-\alpha}\), we modify \(\underline\Phi\) since
\(\underline\Phi>h_R=0\) on \(\Gamma_R\). We first construct
\(\Psi>0\) with \(\mathcal L_m\Psi\geq0\) for large \(|y|\) and
\(\underline\Phi/\Psi\rightarrow0\) uniformly as \(|y|\rightarrow\infty\).

Let \(\nu=(m-2)/2\) and \(F(r)=r^{-\nu}\sqrt{\log r}\) for \(r>1\).
For \(r=|x|>1\), the formula
\(\Delta(F(|x|))=F''(r)+(m-1)F'(r)/r\) gives
\[
  -\Delta(F(|x|))-\frac{c_m}{r^2}F(r)
  =\frac{F(r)}{r^2}
  \left(\nu^2-c_m+\frac{1}{4(\log r)^2}\right).
\]
Set \(\Psi(y)=F(|y|)\) for \(|y|>1\), with \(\Psi(y)=1\) for \(|y|=1\) and \(\Psi(0)=0\).
For \(r=|y|\) sufficiently large and \(1\leq j\leq m\), we have
\[
  \max_{|t|\leq1}\left|\frac{d^4}{dt^4}F(|y+te_j|)\right|
  \leq C r^{-\nu-4}\sqrt{\log r}.
\]
By \eqref{eq:discrete-Taylor-remainder}, we obtain
\[
  \frac{r^2\mathcal L_m\Psi(y)}{\Psi(y)}
  =\nu^2-c_m+\frac{1}{4(\log r)^2}+O(r^{-2}).
\]
Since \(c_m\leq\nu^2\) and \(r^{-2}=o((\log r)^{-2})\), we may increase
\(R_0\) so that \(\mathcal L_m\Psi(y)>0\) for \(|y|\geq R_0\).
Moreover, \(\alpha=\alpha_m\geq\nu\), so
\[
  \frac{\underline\Phi(y)}{\Psi(y)}
  =\frac{r^{\nu-\alpha}(1+Ar^{-2})}{\sqrt{\log r}}
  \longrightarrow0\qquad(r\to\infty).
\]

Set
\(\tau_R=\max_{y\in\Gamma_R}\underline\Phi(y)/\Psi(y)\) and
\(W_R=\underline\Phi-\tau_R\Psi\). Then \(\tau_R\rightarrow0\),
\(W_R\leq0\) on \(\Gamma_R\), and \(\mathcal L_mW_R\leq0\) in \(D_R\).
Choose \(R_1\) with \(\Gamma_0\subset B_{R_1}^*\).
By Lemma~\ref{lem:capacitary-supersolution},
\[
  h_R(y)\geq h_{R_1}(y)\geq\min_{\eta\in\Gamma_0}h_{R_1}(\eta)=:c_0>0,
  \qquad y\in\Gamma_0,\quad R\geq R_1.
\] Since \(W_R\leq\underline\Phi\),
we can choose \(c>0\), independent of \(R\), such that
\(cW_R\leq h_R\) on both boundaries.
The comparison principle gives
\(cW_R\leq h_R\) in \(D_R\). Letting \(R\rightarrow\infty\), we obtain
\(h_m\geq c\underline\Phi\geq c|y|^{-\alpha}\).
Together with \eqref{eq:hm-upper-bound}, we obtain
\(c|y|^{-\alpha}\leq h_m(y)\leq C|y|^{-\alpha}\) for \(|y|\geq R_0\).
Finally, we extend the bounds to all \(y\ne0\).
Since \(h_m(y)>0\) for \(y\ne0\), the finite set
\(\{y\in\mathbb Z^m:0<|y|\leq R_0\}\) satisfies
\[
  0<\min_{0<|y|\leq R_0}h_m(y)|y|^{\alpha_m}
  \leq\max_{0<|y|\leq R_0}h_m(y)|y|^{\alpha_m}<\infty.
\]
Changing \(c\) and \(C\) accordingly, we obtain
\eqref{eq:global-h-comparability} for every \(y\ne0\).
This completes the proof.
\end{proof}

To express the left-hand side of \eqref{eq:discrete-HSM} in terms of
\(v=u/h_m\), we apply the ground-state transformation.
Let \(h:\mathbb Z^m\to[0,\infty)\) satisfy \(h(0)=0\) and
\(h(y)>0\) for \(y\ne0\). For \(v\in C_c(\mathbb Z^m\setminus\{0\})\),
set \((hv)(0)=0\). By \cite[Proposition~4.8]{KellerPinchoverPogorzelski2020},
applied on \(\mathbb Z^m\setminus\{0\}\) with the potential
\(\mathbf1_{\{|y|=1\}}\), we obtain
\begin{equation}
  \label{eq:punctured-ground-state-identity}
  \sum_{y\in\mathbb Z^m}|D(hv)(y)|^2
  =\sum_{\substack{\{y,z\}\subset\mathbb Z^m\setminus\{0\}\\y\sim z}}
  h(y)h(z)|v(y)-v(z)|^2
  +\sum_{y\in\mathbb Z^m\setminus\{0\}}h(y)(\Delta_{\mathrm d}h)(y)|v(y)|^2.
\end{equation}

Taking \(h=h_m\) in \eqref{eq:punctured-ground-state-identity},
we obtain the following lemma.

\begin{lemma}
\label{lem:ground-state-reduction}
Let  \(n,m\geq1\), \(N=n+m\geq3\), and \(m\ne2\).
Assume that \(h_m\) satisfies the conditions of Lemma~\ref{lem:hm-decay}.
For \(u\in C_c(\mathbb Z^N\setminus\Sigma)\), set \(v(x,y)=u(x,y)/h_m(y)\) for \(y\neq0\).
Then
\begin{equation}
  \label{eq:ground-state-defect-full}
  \begin{aligned}
    &\sum_{z\in\mathbb Z^N}|Du(z)|^2
    -c_m\sum_{x\in\mathbb Z^n}\sum_{y\in\mathbb Z^m\setminus\{0\}}
    \frac{|u(x,y)|^2}{|y|^2}\\
    &=\mathcal E_{h_m}(v)
    +\sum_{x\in\mathbb Z^n}\sum_{y\in\mathbb Z^m\setminus\{0\}}
    h_m(y)(\mathcal L_mh_m)(y)|v(x,y)|^2,
  \end{aligned}
\end{equation}
where
\[
  \begin{aligned}
    \mathcal E_{h_m}(v):={}&\sum_{j=1}^n\sum_{x\in\mathbb Z^n}\sum_{y\in\mathbb Z^m\setminus\{0\}}
    h_m(y)^2|D_jv(x,y)|^2\\
    &+\sum_{x\in\mathbb Z^n}
    \sum_{\substack{\{y,z\}\subset\mathbb Z^m\setminus\{0\}\\y\sim z}}
    h_m(y)h_m(z)|v(x,y)-v(x,z)|^2.
  \end{aligned}
\]
Since  \(h_m(y)(\mathcal L_mh_m)(y)\geq0\) for \(y\ne0\), we have
\[\mathcal E_{h_m}(v)\leq
\sum_{z\in\mathbb Z^N}|Du(z)|^2
-c_m\sum_{x\in\mathbb Z^n}\sum_{y\in\mathbb Z^m\setminus\{0\}}
|u(x,y)|^2/|y|^2.\]
Moreover, for
\(1\leq q<\infty\),
\[
  \sum_{z\in\mathbb Z^N}|u(z)|^q
  =\sum_{x\in\mathbb Z^n}\sum_{y\in\mathbb Z^m\setminus\{0\}}h_m(y)^q|v(x,y)|^q.
\]
\end{lemma}
\begin{proof}
For each \(x\in\mathbb Z^n\), apply \eqref{eq:punctured-ground-state-identity}
to \(h=h_m\) and \(v(x,\cdot)\), and subtract
\(c_m\sum_{y\ne0}h_m(y)^2|v(x,y)|^2/|y|^2\).
Since \(h_m\) is independent of \(x\), we have
\(D_ju(x,y)=h_m(y)D_jv(x,y)\) for \(1\leq j\leq n\) and  \(y\ne0\).
Summing  over \(x\), we obtain \eqref{eq:ground-state-defect-full}.

This completes the proof.
\end{proof}

To apply a weighted Sobolev inequality to \(I\widetilde v\),
we next compare the weighted integrals with the sums in
\eqref{eq:ground-state-defect-full}.
Since \(v(x,0)\) is not defined, we set \(\widetilde v(x,0)=v(x,e_1)\).

\begin{lemma}
\label{lem:weighted-Q1-common}
Let  \(n,m\geq1\), \(N=n+m\geq3\), and \(m\ne2\).
Assume that \(h_m\) satisfies the conditions of Lemma~\ref{lem:hm-decay},
that \(h_1(y)=\sqrt{|y|}\) when \(m=1\), and that \(\alpha_m2^*<m\).  For \(u\in C_c(\mathbb Z^N\setminus\Sigma)\),
 set \(v(x,y)=u(x,y)/h_m(y)\) for \(y\ne0\), and define
\begin{equation}
  \label{eq:common-extension}
  \widetilde v(x,y):=\begin{cases}
    v(x,y),&y\neq0,\\
    v(x,e_1),&y=0.
  \end{cases}
\end{equation}
Then there exists \(C=C(m,n)>0\) such that
\begin{equation}
  \label{eq:weighted-Q1-gradient-common}
  \begin{aligned}
    &\int_{\mathbb R^n\times\mathbb R^m}|Y|^{-2\alpha_m}
    |\nabla I\widetilde v(X,Y)|^2\,dX\,dY\\
    &\leq C\left[\sum_{z\in\mathbb Z^N}|Du(z)|^2
    -c_m\sum_{x\in\mathbb Z^n}\sum_{y\in\mathbb Z^m\setminus\{0\}}
    \frac{|u(x,y)|^2}{|y|^2}\right],
  \end{aligned}
\end{equation}
and
\begin{equation}
  \label{eq:weighted-Q1-Lp-common}
  \sum_{x\in\mathbb Z^n}\sum_{y\in{\mathbb Z^m\setminus\{0\}}}h_m(y)^{2^*}|v(x,y)|^{2^*}
  \leq C\int_{\mathbb R^n\times\mathbb R^m}|Y|^{-\alpha_m2^*}
  |I\widetilde v(X,Y)|^{2^*}\,dX\,dY.
\end{equation}
\end{lemma}
\begin{proof}
We first prove \eqref{eq:weighted-Q1-gradient-common}.
Let \(F=I\widetilde v\), \(\alpha=\alpha_m\), and \(R_*=4\sqrt m\).

For \(|y|\geq R_*\) and \((X,Y)\in Q_{(x,y)}\), we have
\(3|y|/4\leq|Y|\leq5|y|/4\).
Let \(z=(\xi,\eta)\) and \(z'=(\xi',\eta')\) belong to
\(Q_{(x,y)}\cap\mathbb Z^N\), with \(z\sim z'\).
By \eqref{eq:global-h-comparability},
\[
  |Y|^{-2\alpha}\leq C h_m(\eta)h_m(\eta').
\]
Here \(\eta=\eta'\) when \(\xi\neq\xi'\), so the product is \(h_m(\eta)^2\).
By \eqref{eq:Q1-local-Jensen},
\[
  \int_{Q_{(x,y)}}|Y|^{-2\alpha}|\nabla F|^2\,dX\,dY
  \leq C\sum_{\substack{\{z,z'\}\subset Q_{(x,y)}\cap\mathbb Z^N\\z\sim z'}}
  h_m(\eta)h_m(\eta')|\widetilde v(z)-\widetilde v(z')|^2.
\]
For each \(z\sim z'\),
\(\#\{(x,y)\in\mathbb Z^N:z,z'\in Q_{(x,y)}\}=2^{N-1}\).
Summing over \(x\in\mathbb Z^n\) and \(y\in\mathbb Z^m\) with \(|y|\geq R_*\),
we obtain
\begin{equation}\label{eq:Q1-far-cubes}
  \sum_{x\in\mathbb Z^n}\sum_{\substack{y\in\mathbb Z^m\\|y|\geq R_*}}
  \int_{Q_{(x,y)}}|Y|^{-2\alpha}|\nabla F|^2\,dX\,dY
  \leq C\mathcal E_{h_m}(v).
\end{equation}

Next, we consider \(|y|<R_*\). Since \(2\alpha<m\),
\(\int_{y+[0,1]^m}|Y|^{-2\alpha}\,dY\leq C\).
By \eqref{eq:Q1-local-Jensen}, we have
\begin{equation}
  \label{eq:Q1-near-cubes}
  \int_{Q_{(x,y)}}|Y|^{-2\alpha}|\nabla F|^2\,dX\,dY
  \leq C\sum_{\substack{\{z,z'\}\subset Q_{(x,y)}\cap\mathbb Z^N\\z\sim z'}}
  |\widetilde v(z)-\widetilde v(z')|^2.
\end{equation}
Set
\[
  a:=\min_{\substack{\eta\in\mathbb Z^m\\0<|\eta|\leq R_*+\sqrt m}}h_m(\eta)>0.
\]
For \(z=(\xi,\eta)\), \(z'=(\xi',\eta')\) with
\(\{z,z'\}\subset Q_{(x,y)}\cap\mathbb Z^N\), \(z\sim z'\), and
\(\eta,\eta'\neq0\), we have
\[
  |\widetilde v(z)-\widetilde v(z')|^2
  \leq a^{-2}h_m(\eta)h_m(\eta')|\widetilde v(z)-\widetilde v(z')|^2.
\]
Hence,
\begin{equation}
  \label{eq:Q1-nonzero-vertices}
  \sum_{x\in\mathbb Z^n}\sum_{|y|<R_*}
  \sum_{\substack{\{z,z'\}\subset Q_{(x,y)}\cap\mathbb Z^N\\z\sim z',\ \eta,\eta'\neq0}}
  |\widetilde v(z)-\widetilde v(z')|^2
  \leq 2^{N-1}a^{-2}\mathcal E_{h_m}(v).
\end{equation}

It remains to estimate
\(|\widetilde v(z)-\widetilde v(z')|^2\) when \(z\sim z'\) and
\(\eta=0\) or \(\eta'=0\).
By \eqref{eq:common-extension}, for \(1\leq i\leq n\),
\begin{equation}
  \label{eq:Q1-zero-horizontal}
  \widetilde v(x,0)-\widetilde v(x-e_i,0)
  =v(x,e_1)-v(x-e_i,e_1),\qquad
  \widetilde v(x,0)-v(x,e_1)=0.
\end{equation}
 For \(m\geq3\), to estimate \(|v(x,e_1)-v(x,\sigma e_j)|^2\) using values at
adjacent vertices in \(\mathbb Z^m\setminus\{0\}\), we connect
\(e_1\) to \(\sigma e_j\) by the following sequences,
where \(\sigma\in\{-1,1\}\), \(1\leq j\leq m\), and \(\sigma e_j\neq e_1\):
\[
  e_1\rightarrow e_1+\sigma e_j\rightarrow\sigma e_j\quad(j\neq1),
  \qquad e_1\rightarrow e_1+e_2\rightarrow e_2\rightarrow-e_1+e_2\rightarrow-e_1.
\]

For each sequence \(y_0=e_1,\ldots,y_\ell=\sigma e_j\), we have
\(\ell\leq4\), \(0<|y_i|\leq\sqrt2\), and \(h_m(y_i)\geq a\).
By the Cauchy--Schwarz
inequality,
\begin{equation}
  \label{eq:Q1-zero-vertical}
  \begin{aligned}
    |\widetilde v(x,0)-v(x,\sigma e_j)|^2
    &=\left|\sum_{i=0}^{\ell-1}\bigl(v(x,y_i)-v(x,y_{i+1})\bigr)\right|^2\\
    &\leq4a^{-2}\sum_{i=0}^{\ell-1}h_m(y_i)h_m(y_{i+1})
    |v(x,y_{i+1})-v(x,y_i)|^2.
  \end{aligned}
\end{equation}
For \(m=1\), \eqref{eq:common-extension} gives
\(\widetilde v(x,0)-v(x,-1)=v(x,1)-v(x,-1)\).
By \eqref{eq:one-dimensional-defect} and \eqref{eq:ground-state-defect-full},
we have
\begin{equation}
  \label{eq:Q1-one-dimensional-bridge}
  \begin{aligned}
    \sum_{x\in\mathbb Z^n}|v(x,1)-v(x,-1)|^2
    &\leq2\sum_{x\in\mathbb Z^n}
    \bigl(|v(x,1)|^2+|v(x,-1)|^2\bigr)\\
    &\leq\frac2{\rho_1}\left[
    \sum_{z\in\mathbb Z^{n+1}}|Du(z)|^2
    -\frac14\sum_{x\in\mathbb Z^n}\sum_{y\in\mathbb Z\setminus\{0\}}
    \frac{|u(x,y)|^2}{y^2}\right].
  \end{aligned}
\end{equation}
By \eqref{eq:Q1-near-cubes}--\eqref{eq:Q1-one-dimensional-bridge} and
\eqref{eq:ground-state-defect-full}, summing over
\(x\in\mathbb Z^n\) and \(y\in\mathbb Z^m\) with \(|y|<R_*\), we obtain
\begin{equation}\label{eq:Q1-near-cubes-summed}
  \begin{aligned}
    &\sum_{x\in\mathbb Z^n}\sum_{|y|<R_*}
    \int_{Q_{(x,y)}}|Y|^{-2\alpha}|\nabla F|^2\,dX\,dY\\
    &\leq C\left[\sum_{z\in\mathbb Z^N}|Du(z)|^2
    -c_m\sum_{x\in\mathbb Z^n}\sum_{y\in\mathbb Z^m\setminus\{0\}}
    \frac{|u(x,y)|^2}{|y|^2}\right].
  \end{aligned}
\end{equation}
Combining \eqref{eq:Q1-far-cubes}, \eqref{eq:Q1-near-cubes-summed}, and
\eqref{eq:ground-state-defect-full}, we obtain
\eqref{eq:weighted-Q1-gradient-common}.

 Finally, to prove \eqref{eq:weighted-Q1-Lp-common} for both signs of \(\alpha\),
for each \(y\ne0\) choose the unit cube
\[
  \widehat Q_y:=\prod_{j=1}^m J_j(y_j),\qquad
  J_j(y_j):=\begin{cases}
    [y_j,y_j+1],&y_j\geq0,\\
    [y_j-1,y_j],&y_j<0.
  \end{cases}
\]
For \(Y\in \widehat Q_y\), we have
\(|y_j|\leq|Y_j|\leq|y_j|+1\), and hence
\[
  |y|\leq|Y|\leq|y|+\sqrt m\leq(1+\sqrt m)|y|.
\]
 By \eqref{eq:global-h-comparability}, we obtain
\(h_m(y)^{2^*}\leq C|Y|^{-\alpha2^*}\).
Since \((x,y)\) is a vertex of \((x+[0,1]^n)\times \widehat Q_y\),
Lemma~\ref{lem:Q1-basic-estimates}\textup{(iii)} gives
\[
  h_m(y)^{2^*}|v(x,y)|^{2^*}
  \leq C\int_{(x+[0,1]^n)\times \widehat Q_y}
  |Y|^{-\alpha2^*}|F(X,Y)|^{2^*}\,dX\,dY.
\]
By the definition of \(J_j\),
\[
  |J_j(a)\cap J_j(b)|=0\quad(a,b\in\mathbb Z,\ a\ne b).
\]
Consequently, for \((x,y)\ne(x',y')\),
\[
  \left|\bigl((x+[0,1]^n)\times\widehat Q_y\bigr)
  \cap\bigl((x'+[0,1]^n)\times\widehat Q_{y'}\bigr)\right|=0.
\]
Summing over \(x\in\mathbb Z^n\) and  \(y\ne0\), we obtain
\eqref{eq:weighted-Q1-Lp-common}.
This completes the proof.
\end{proof}

Combining Lemma~\ref{lem:weighted-Q1-common} with the following inequality
will give the required bound for \(\|u\|_{\ell^{2^*}(\mathbb Z^N)}^2\).
For \(A\subset\mathbb R^d\), let \(C_c^{0,1}(A)\) denote the
real-valued Lipschitz functions with compact support in \(A\).

\begin{lemma}
\label{lem:weighted-sobolev}
Let  \(n,m\geq1\) and \(N=n+m\geq3\). Assume that either
\(m=1\) and \(\alpha=-1/2\), or \(m\geq3\) and \(\alpha<m/2^*\).
There exists \(C=C(m,n,\alpha)>0\) such that
\begin{equation}
  \label{eq:weighted-continuous-common}
  \left(\int_{\mathbb R^n\times\mathbb R^m}|Y|^{-\alpha2^*}
  |F(X,Y)|^{2^*}\,dX\,dY\right)^{2/2^*}
  \leq C\int_{\mathbb R^n\times\mathbb R^m}|Y|^{-2\alpha}
  |\nabla F(X,Y)|^2\,dX\,dY
\end{equation}
for every \(F\in C_c^{0,1}(\mathbb R^{n+m})\).
\end{lemma}
\begin{proof}
To prove \eqref{eq:weighted-continuous-common}, we apply \cite[Theorem~2.2.1 and Remark~2.2.3]{GazziniThesis2008} (see also \cite{GazziniMusina2009}) with
\[
 p=2,\qquad k=m,\qquad a=-2\alpha,\qquad q=2^*,\qquad\Omega=\mathbb R^N.
\]
Since \(n\geq1\) and \(\alpha<m/2^*\) in both cases, we have \(1\leq m<N\) and
\[
 \begin{gathered}
 N-2-2\alpha>\frac{n(N-2)}N>0,\qquad
 -2\alpha>\frac{(2-N)m}N,\\
 \max\left\{2,\frac{2(N-m)}{N-2-2\alpha}\right\}<2^*.
 \end{gathered}
\]
Hence \eqref{eq:weighted-continuous-common} holds for every \(F\in C_c^\infty(\mathbb R^N)\).
To apply \eqref{eq:weighted-continuous-common} to \(I\widetilde v\),
we now consider \(F\in C_c^{0,1}(\mathbb R^{n+m})\).
Choose \(\rho\in C_c^\infty(\mathbb R^{n+m})\) with \(\rho\geq0\),
\(\operatorname{supp}\rho\subset\{|Z|<1\}\), and
\(\int_{\mathbb R^{n+m}}\rho(Z)\,dZ=1\). Set
\[
  F_\varepsilon(Z):=\int_{\mathbb R^{n+m}}
  \varepsilon^{-(n+m)}\rho\bigl((Z-W)/\varepsilon\bigr)F(W)\,dW,
  \qquad0<\varepsilon<1.
\]
Choose \(M>0\) with \(\operatorname{supp}F\subset\{|Z|\leq M\}\).
For \(0<\varepsilon<1\),
\[
  \operatorname{supp}F_\varepsilon\subset\{|Z|\leq M+1\},\qquad
  \|F_\varepsilon\|_\infty\leq\|F\|_\infty,\qquad
  \|\nabla F_\varepsilon\|_\infty\leq\|\nabla F\|_\infty.
\]
Moreover, \(F_\varepsilon\rightarrow F\) and \(\nabla F_\varepsilon\rightarrow\nabla F\)
almost everywhere
as \(\varepsilon\rightarrow0^+\). Since \(\alpha2^*<m\)
and \(2\alpha<m\), by Lebesgue's dominated convergence theorem, we obtain the result.
This completes the proof.
\end{proof}

We now combine the ground-state transformation in Lemma~\ref{lem:ground-state-reduction}
with Lemmas~\ref{lem:weighted-Q1-common} and~\ref{lem:weighted-sobolev}
to prove  \(\overline\kappa_{m,n}>0\).

\begin{proposition}
\label{prop:positivity-sufficiency}
Let  \(n,m\geq1\), \(N=n+m\geq3\), and \(m\ne2\). If
\(c_m>(m-2)^2/4-(n/N)^2\), then \(\overline\kappa_{m,n}>0\).
\end{proposition}
\begin{proof}
Since \(m/2^*=(m-2)/2+n/N\), the definition of \(\alpha_m\) gives
\[
  \alpha_m-\frac m{2^*}
  =\sqrt{\left(\frac{m-2}{2}\right)^2-c_m}-\frac nN<0.
\]
 Hence,  \(\alpha_m<m/2^*\).
For \(u\in C_c(\mathbb Z^N\setminus\Sigma)\), set
\(v(x,y)=u(x,y)/h_m(y)\) for \(y\neq0\), and use the extension
\eqref{eq:common-extension}. Its interpolant belongs to
\(C_c^{0,1}(\mathbb R^{n+m})\). By Lemmas~\ref{lem:ground-state-reduction}--\ref{lem:weighted-sobolev}, we
obtain
\begin{equation}
  \label{eq:HSM-sufficiency-bound}
  \begin{aligned}
    \left(\sum_{z\in\mathbb Z^N}|u(z)|^{2^*}\right)^{2/2^*}
    &=\left(\sum_{x\in\mathbb Z^n}\sum_{y\in\mathbb Z^m\setminus\{0\}}
    h_m(y)^{2^*}|v(x,y)|^{2^*}\right)^{2/2^*}\\
    &\leq C\left(\int_{\mathbb R^{n+m}}|Y|^{-\alpha_m2^*}
    |I\widetilde v|^{2^*}\,dX\,dY\right)^{2/2^*}\\
    &\leq C\int_{\mathbb R^{n+m}}|Y|^{-2\alpha_m}
    |\nabla I\widetilde v|^2\,dX\,dY\\
    &\leq C\left[\sum_{z\in\mathbb Z^N}|Du(z)|^2
    -c_m\sum_{x\in\mathbb Z^n}\sum_{y\in\mathbb Z^m\setminus\{0\}}\frac{|u(x,y)|^2}{|y|^2}\right].
  \end{aligned}
\end{equation}
Hence \(\overline\kappa_{m,n}\geq C^{-1}>0\).
This completes the proof.
\end{proof}

 To prove \(\overline\kappa_{m,n}=0\) for \(m\geq3\), \(n\geq1\), when the condition in
\eqref{eq:positivity-criterion} fails, we will use test functions of the form
\(\phi_R(x)\eta_R(y)h_m(y)\), with \(\eta_R(e_1)=1\).
 For \(m\geq3\), taking \(h=h_m\) in \eqref{eq:punctured-ground-state-identity} and using \(h_m(o)=1\) and
\(\mathcal L_mh_m=\mu_m\delta_o\), we obtain
\begin{equation}
  \label{eq:ground-state-defect-y}
  \begin{aligned}
    &\sum_{y\in\mathbb Z^m}|D(h_mv)(y)|^2
    -c_m\sum_{y\in\mathbb Z^m\setminus\{0\}}\frac{h_m(y)^2|v(y)|^2}{|y|^2}\\
    &=\sum_{\substack{\{y,z\}\subset\mathbb Z^m\setminus\{0\}\\y\sim z}}
    h_m(y)h_m(z)|v(y)-v(z)|^2+\mu_m|v(o)|^2.
  \end{aligned}
\end{equation}
Taking \(v=\eta_R\) in \eqref{eq:ground-state-defect-y} leaves
\(\mu_m|\eta_R(e_1)|^2=\mu_m\). We therefore first prove
\(\mu_m=0\) under the following assumption.

\begin{lemma}
\label{lem:criticality-zero-defect}
Let \(m\geq3\) and  assume that \(c_m<(m-2)^2/4\).
Then  \(\mu_m=0\) and
\(\mathcal L_mh_m=0\) on \(\mathbb Z^m\setminus\{0\}\).
\end{lemma}
\begin{proof}[Proof]
To prove \(\mu_m=0\), it suffices to find
\(f_j\in C_c(\mathbb Z^m\setminus\{0\})\) such that
\(f_j(o)=1\) and \(q_m(f_j)\to0\).
For this purpose, we apply \cite{DasKellerPinchover2026}.
Take \(X=\mathbb Z^m\setminus\{0\}\) and define
\[
  b(y,z):=\begin{cases}1,&y\sim z,\\0,&\text{else},\end{cases}
  \qquad V(y):=\mathbf1_{\{|y|=1\}},\qquad y,z\in X.
\] For \(f\in C_c(X)\), extended by \(f(0)=0\),
\[
  \frac12\sum_{y\in X}\sum_{\substack{z\in X\\z\sim y}}
  |f(z)-f(y)|^2+\sum_{|y|=1}|f(y)|^2
  =\sum_{y\in\mathbb Z^m}|Df(y)|^2.
\]
For \(g(y)=|y|^{-2}\),
by Proposition~\ref{prop:euclidean-hardy-comparison} and the assumption on \(c_m\),
we have
\[
  \lambda_0(g)=c_m<
  \lambda_\infty(g)=\sup_{R>0}I_R=\frac{(m-2)^2}{4}.
\]

Since \(\lambda_0(g)<\lambda_\infty(g)\), by
\cite[Lemma~5.4]{DasKellerPinchover2026}, we have
\[
  \inf_{\substack{f\in C_c(X)\\f(o)\ne0}}
  \frac{q_m(f)}{|f(o)|^2}=0.
\]
By \(q_m(tf)=t^2q_m(f)\), we can choose \(f_j\in C_c(X)\) such that
\[
  f_j(o)=1,\qquad q_m(f_j)\leq\frac1j.
\]
For each \(j\), choose an integer \(R_j\geq\max\{3,j\}\) such that
\(\operatorname{supp}f_j\subset B_{R_j}^*\).
By the definition of \(\mu_{R_j}\), we obtain
\[
  0\leq\mu_m\leq\mu_{R_j}\leq q_m(f_j)\leq\frac1j
  \longrightarrow0\qquad(j\to\infty).
\]
Hence \(\mu_m=0\), and Lemma~\ref{lem:capacitary-supersolution} gives
\(\mathcal L_mh_m=0\) on \(X\).
This completes the proof.
\end{proof}

We now assume that \(m\geq3\) and
\[
  c_m\leq\left(\frac{m-2}{2}\right)^2-\left(\frac nN\right)^2,
  \quad\text{or equivalently}\quad \alpha_m\geq\frac m{2^*}.
\]
Since \(n\geq1\), we have \(c_m<(m-2)^2/4\), and
Lemma~\ref{lem:criticality-zero-defect} gives \(\mu_m=0\).

By \eqref{eq:global-h-comparability}, it suffices to estimate
\(\sum_{0<|y|\leq R}|y|^{-\beta}\) for \(\beta=2\alpha_m\) and
\(\beta=2^*\alpha_m\).
For \(R\geq2\) and \(\beta>0\),
\begin{equation}
  \label{eq:lattice-power-sum-upper}
  \sum_{\substack{y\in\mathbb Z^m\\0<|y|\leq R}}|y|^{-\beta}
  \leq C\begin{cases}
    R^{m-\beta},&\beta<m,\\
    1+\log R,&\beta=m,\\
    1,&\beta>m.
  \end{cases}
\end{equation}
Indeed, the set \(\{y\in\mathbb Z^m:2^j<|y|\leq2^{j+1}\}\) contains at most \(C2^{jm}\) vertices,
so summing \(C2^{j(m-\beta)}\) gives \eqref{eq:lattice-power-sum-upper}.

At \(\beta=m\), we also need a lower bound:
\begin{equation}
  \label{eq:lattice-power-sum-lower}
  \sum_{\substack{y\in\mathbb Z^m\\1\leq|y|\leq R}}|y|^{-m}
  \geq c\log R
\end{equation}
for all sufficiently large \(R\). To see this, the set
\(\{y:\max_i|y_i|=k\}\) has
\((2k+1)^m-(2k-1)^m\geq ck^{m-1}\)
vertices, each with
\(|y|\leq\sqrt m\,k\). Summing their contributions for
\(1\leq k\leq\lfloor R/\sqrt m\rfloor\) gives
\(c\sum_{k\leq R/\sqrt m}k^{-1}\geq c'\log R\).

We now cut off \(h_m\) in the \(y\)-variables. For integers \(R\geq2\), set
\[
  \eta_R(y):=\begin{cases}
    1,&|y|\leq R,\\
    2-|y|/R,&R<|y|<2R,\\
    0,&|y|\geq2R,
  \end{cases}
  \qquad a_R(y):=\eta_R(y)h_m(y).
\]
Then \(a_R\in C_c(\mathbb Z^m)\) and \(a_R(0)=0\).
Define
\[
  \begin{gathered}
    \delta_R:=\sum_{y\in\mathbb Z^m}|Da_R(y)|^2
    -c_m\sum_{y\in\mathbb Z^m\setminus\{0\}}\frac{|a_R(y)|^2}{|y|^2},\\
    A_R:=\sum_{y\in\mathbb Z^m}|a_R(y)|^2,
    \qquad M_R:=\sum_{y\in\mathbb Z^m}|a_R(y)|^{2^*}.
  \end{gathered}
\]
By \eqref{eq:ground-state-defect-y} and \(\mu_m=0\),
\[
  \delta_R=\sum_{\substack{\{y,z\}\subset\mathbb Z^m\setminus\{0\}\\y\sim z}}
  h_m(y)h_m(z)|\eta_R(y)-\eta_R(z)|^2.
\]
If \(y\sim z\) and \(\eta_R(y)\neq\eta_R(z)\), then
\[
  R-1\leq|y|,|z|\leq2R+1,\qquad
  h_m(y)h_m(z)\leq CR^{-2\alpha_m},\qquad
  |\eta_R(y)-\eta_R(z)|\leq R^{-1}.
\]
Since \(\#\{\{y,z\}:y\sim z,\ \eta_R(y)\neq\eta_R(z)\}\leq CR^m\),
\begin{equation}
  \label{eq:y-cutoff-deficit-bound}
  0\leq\delta_R\leq CR^{m-2\alpha_m-2}.
\end{equation}

Also, \eqref{eq:global-h-comparability} and
\eqref{eq:lattice-power-sum-upper} give
\begin{equation}
  \label{eq:aR-l2-bound}
  A_R\leq C\begin{cases}
    R^{m-2\alpha_m},&2\alpha_m<m,\\
    1+\log R,&2\alpha_m=m,\\
    1,&2\alpha_m>m.
  \end{cases}
\end{equation}

Since \(a_R(o)=1\), we have \(M_R\geq1\).
When \(\alpha_m=m/2^*\), the lower bounds in
\eqref{eq:global-h-comparability} and \eqref{eq:lattice-power-sum-lower}
give \(M_R\geq c\log R\).

To prove \(\overline\kappa_{m,n}=0\), we choose the test function
\(u_R(x,y)=\phi_R(x)a_R(y)\), where \(\phi_R\) is as in
Lemma~\ref{lem:standard-x-cutoff}. Then \(u_R\in C_c(\mathbb Z^N\setminus\Sigma)\)
and \(u_R(0,o)=1\). By  \eqref{eq:product-energy-hardy} with \(a=a_R\), and by the definitions
of \(\delta_R\), \(A_R\), and \(M_R\), we obtain
\[
  \begin{aligned}
    \sum_{z\in\mathbb Z^N}|Du_R(z)|^2
    -c_m\sum_{x\in\mathbb Z^n}\sum_{y\in\mathbb Z^m\setminus\{0\}}\frac{|u_R(x,y)|^2}{|y|^2}
    &=A_R\sum_{x\in\mathbb Z^n}|D\phi_R(x)|^2
    +\delta_R\sum_{x\in\mathbb Z^n}|\phi_R(x)|^2,\\
    \sum_{z\in\mathbb Z^N}|u_R(z)|^{2^*}
    &=M_R\sum_{x\in\mathbb Z^n}|\phi_R(x)|^{2^*}.
  \end{aligned}
\]
Using Lemma~\ref{lem:standard-x-cutoff}, together with
\(n-2-2n/2^*=-2m/N\) and \(n-2n/2^*=2n/N\),
gives
\[
  \begin{aligned}
    &\frac{\displaystyle\sum_{z\in\mathbb Z^N}|Du_R(z)|^2
    -c_m\sum_{x\in\mathbb Z^n}\sum_{y\in\mathbb Z^m\setminus\{0\}}|u_R(x,y)|^2/|y|^2}
    {\displaystyle\left(\sum_{z\in\mathbb Z^N}|u_R(z)|^{2^*}\right)^{2/2^*}}\\
    &=\frac{\displaystyle A_R\sum_{x\in\mathbb Z^n}|D\phi_R(x)|^2
    +\delta_R\sum_{x\in\mathbb Z^n}|\phi_R(x)|^2}
    {\displaystyle M_R^{2/2^*}\left(\sum_{x\in\mathbb Z^n}|\phi_R(x)|^{2^*}\right)^{2/2^*}}
    \leq C\frac{A_RR^{-2m/N}+\delta_RR^{2n/N}}{M_R^{2/2^*}}.
  \end{aligned}
\]

Finally, by \eqref{eq:y-cutoff-deficit-bound}, we have
\(\delta_RR^{2n/N}\leq CR^{2(m/2^*-\alpha_m)}\).
If \(\alpha_m>m/2^*\), then \(\delta_RR^{2n/N}\to0\) as \(R\to\infty\).

For sufficiently large \(R\), by
\eqref{eq:aR-l2-bound} and \(M_R\geq1\), we have
\[
  \frac{A_RR^{-2m/N}}{M_R^{2/2^*}}
  \leq A_RR^{-2m/N}
  \leq C\begin{cases}
    R^{2(m/2^*-\alpha_m)},&2\alpha_m<m,\\
    (1+\log R)R^{-2m/N},&2\alpha_m=m,\\
    R^{-2m/N},&2\alpha_m>m
  \end{cases}
  \xrightarrow[R\rightarrow\infty]{}0.
\]
If \(\alpha_m=m/2^*\),
then
\[
  A_RR^{-2m/N}\leq C,\qquad \delta_RR^{2n/N}\leq C,\qquad
  M_R^{2/2^*}\geq c(\log R)^{2/2^*}\rightarrow\infty.
\]
Thus, in both cases,
\begin{equation}
  \label{eq:mge3-necessity-null-sequence}
  \frac{\displaystyle\sum_{z\in\mathbb Z^N}|Du_R(z)|^2
  -c_m\sum_{x\in\mathbb Z^n}\sum_{y\in\mathbb Z^m\setminus\{0\}}|u_R(x,y)|^2/|y|^2}
  {\displaystyle\left(\sum_{z\in\mathbb Z^N}|u_R(z)|^{2^*}\right)^{2/2^*}}
  \rightarrow0\qquad\text{as }R\rightarrow\infty.
\end{equation}

We can now conclude \eqref{eq:positivity-criterion}.
For \(m=2\), Proposition~\ref{prop:m2} gives positivity, and
\(c_2=0>-(n/(n+2))^2\).
For \(m\ne2\), Proposition~\ref{prop:positivity-sufficiency} proves sufficiency.
When \(m=1\) and \(n\geq2\), we have \(c_1=1/4>1/4-(n/(n+1))^2\).
Thus, if the condition fails, then \(m\geq3\), and
\eqref{eq:mge3-necessity-null-sequence} gives \(\overline\kappa_{m,n}\leq0\).
By \eqref{eq:general-sobolev-upper}, we obtain \(\overline\kappa_{m,n}=0\).
This proves \eqref{eq:positivity-criterion}.

To complete the proof of Theorem~\ref{thm:main}, we now estimate \(c_m\)
for \(m\geq9\). Since \(n/(n+m)<1\), the following bound
implies \(\overline\kappa_{m,n}=0\) for every \(m\geq10\) and \(n\geq1\).

\begin{lemma}
\label{lem:high-hardy-gap}
For every \(m\geq10\), \(c_m<(m-2)^2/4-1\).
\end{lemma}
\begin{proof}
We first consider \(m\geq11\). The test functions in
\cite[Theorem~1.1 and Section~5]{Lizama2026} give
\[
  c_m\leq3m-\sqrt{m^2+8m-8}\leq\frac{2m^3}{m^2+2m-2}.
\]

Since \((m-2)^2/4-1=m(m-4)/4\), it suffices to check
\[
  \frac{2m^3}{m^2+2m-2}<\frac{m(m-4)}4
  \quad\Longleftrightarrow\quad m^3-10m^2-10m+8>0.
\]

For \(m\geq11\),
we have
\[
  m^3-10m^2-10m+8=(m-11)^3+23(m-11)^2+133(m-11)+19>0.
\]

It remains to prove \(c_{10}<15\). For this purpose, we use a finitely supported test function.
For \(1\le j\le6\) and \(0\le k\le3\), define
\[
  \Omega_j=\left\{y\in\{-1,0,1\}^{10}:|y|^2=j\right\},
\]
\[
  A_k=\left\{y\in\{-2,-1,0,1,2\}^{10}:
  \#\{i:|y_i|=2\}=1,\quad \#\{i:|y_i|=1\}=k\right\}.
\]
Define \(f\) by the following table.
\begin{center}
\begin{tabular}{crrr}
\toprule
\(K\)&\(f(y)\), \(y\in K\)&\(|K|\)&\(\sum_{z\sim y}f(z)\), \(y\in K\)\\\midrule
\(\Omega_1\)&250&20&1285\\
\(\Omega_2\)&70&180&868\\
\(\Omega_3\)&22&960&331\\
\(\Omega_4\)&8&3360&128\\
\(\Omega_5\)&3&8064&50\\
\(\Omega_6\)&1&13440&18\\
\(A_0\)&25&20&394\\
\(A_1\)&8&360&143\\
\(A_2\)&3&2880&52\\
\(A_3\)&1&13440&17\\\bottomrule
\end{tabular}
\end{center}
For
\[
  y\in\mathbb Z^{10}\setminus
  \left(\bigcup_{j=1}^6\Omega_j\cup\bigcup_{k=0}^3 A_k\right),
\]
set \(f(y)=0\).
By the definitions of \(\Omega_j\) and \(A_k\), we have
\[
  |\Omega_j|=2^j\binom{10}{j},\qquad |A_k|=20\cdot2^k\binom9k.
\]
The table also gives \(\sum_{z\sim y}f(z)\) for \(y\in K\).
By the definition of \(f\) and the values in the table, we obtain
\[
  \begin{aligned}
    \sum_{y\in\mathbb Z^{10}}|Df(y)|^2
    &=\sum_{y\in\mathbb Z^{10}}f(y)
    \left(20f(y)-\sum_{z\sim y}f(z)\right)=28\,920\,640,\\
    \sum_{y\in\mathbb Z^{10}\setminus\{0\}}\frac{|f(y)|^2}{|y|^2}
    &=\frac{9\,651\,841}{5}.
  \end{aligned}
\]
Hence
\[
  \sum_{y\in\mathbb Z^{10}}|Df(y)|^2
  -15\sum_{y\in\mathbb Z^{10}\setminus\{0\}}\frac{|f(y)|^2}{|y|^2}
  =-34\,883<0.
\]
Since \(0\ne f\in C_c(\mathbb Z^{10})\) and \(f(0)=0\),
the definition of \(c_{10}\) gives \(c_{10}<15\).
This completes the proof.
\end{proof}

We are now in a position to prove Theorem~\ref{thm:main}.\par\nopagebreak
\begin{proof}[Proof of Theorem~\ref{thm:main}]
For \(1\leq m\leq8\), we have \(c_m=(m-2)^2/4\).
This follows from \cite[Theorem~1.3]{HuangYe2024} with \(\alpha=0\)
for \(m=1\), \cite[Remark~1.2]{Gupta2024} for \(m=2\),
\cite[Theorem~1.1]{Lizama2026} for \(m=3,4\), and
\cite[Theorem~2.1]{AlpayLow2026} for \(5\leq m\leq8\).
Since \(n\geq1\), \eqref{eq:positivity-criterion} gives
\(\overline\kappa_{m,n}>0\).

For \(m\geq10\), Lemma~\ref{lem:high-hardy-gap} gives
\(c_m<(m-2)^2/4-1\). As \(n/(n+m)<1\),
\eqref{eq:positivity-criterion} gives \(\overline\kappa_{m,n}=0\).

It remains to consider \(m=9\). By Proposition~\ref{prop:nine-bounds},
\[
  \frac{25137}{2116}<\frac{297}{25}\leq c_9<\frac{1440}{121}.
\]
Set \(T_n=49/4-(n/(n+9))^2\). This sequence is strictly decreasing,
with \(T_{13}=1440/121\) and \(T_{14}=25137/2116\).
Thus \(c_9<T_n\) for \(1\leq n\leq13\), and \(c_9>T_n\) for
\(n\geq14\). Applying \eqref{eq:positivity-criterion} proves the
remaining assertions.
This completes the proof.
\end{proof}

\section[Comparison with the best constant for the Sobolev inequality on ZN]{Comparison with the best constant for the Sobolev inequality on \(\mathbb Z^N\)}
\label{sec:sobolev-comparison}

In this section, we prove Theorem~\ref{thm:sobolev-comparison}.
To prove \(\overline\kappa_{m,n}<\overline S_N\) when
\(N\geq4\) and \(c_m>0\), we construct test functions by cutting off
and translating an extremal function for \(\overline S_N\).
Since \(\overline\kappa_{2,n}=\overline S_{n+2}\) has already been proved
in Proposition~\ref{prop:m2}, it remains to consider \((m,n)=(1,2)\).
To prove \(\overline\kappa_{1,2}=\overline S_3\), we first compare
the Green functions in Proposition~\ref{prop:green-order} and then apply a duality argument.
This argument also gives the nonattainment result in Theorem~\ref{thm:attainment}.

\subsection[The case m not equal to 2 and N at least 4]{The case \(m\ne2\) and \(N\ge4\)}
For the construction of the test functions, we first recall the
existence result of He and Ji \cite[Theorem~1.8]{HeJi2026}:
for \(N\geq3\), \(\overline S_N\) is attained by a positive function
in \(D^{1,2}(\mathbb Z^N)\).

For the case \(\varepsilon=1\) in \cite{HeJi2026}, we have
\[
  \|\nabla_1^N u\|_{\ell^2(\mathbb V_1^N)}^2
  =\sum_{j=1}^N\sum_{z\in\mathbb Z^N}|u(z+e_j)-u(z)|^2
  =\sum_{z\in\mathbb Z^N}|Du(z)|^2.
\]
By \cite[Theorem~3.4]{HanLi2026}, with \(p=2\) and \(a=b=0\), we have
\begin{equation}
  \label{eq:full-space-identification}
  D^{1,2}(\mathbb Z^N)
  =\{u\in\ell^{2^*}(\mathbb Z^N):\sum_{z\in\mathbb Z^N}|Du(z)|^2<\infty\}.
\end{equation}
To construct the test functions in Lemma~\ref{lem:strict-sobolev},
for integers \(R\geq2\), define
\[
  \chi_R(z)=\begin{cases}
    1,&|z|\leq R,\\
    2-|z|/R,&R<|z|<2R,\\
    0,&|z|\geq2R,
  \end{cases}\qquad z\in\mathbb Z^N.
\]
By the definition of \(\chi_R\), we have \(\chi_R\in C_c(\mathbb Z^N)\),
\(\chi_R=1\) for \(|z|\leq R\), and \(|D_j\chi_R|\leq R^{-1}\)
for \(1\leq j\leq N\).

In Lemma~\ref{lem:strict-sobolev}, we use \(\chi_{R_k}w\), where
\(w>0\) is an extremal function for \(\overline S_N\).
To choose \(R_k\rightarrow\infty\) satisfying \eqref{eq:cutoff-relative-bound},
we first prove that \(w\in\ell^s(\mathbb Z^N)\) for every \(s>N/(N-2)\).

\begin{lemma}
\label{lem:free-summability}
Let \(N\geq3\), and let \(w>0\) be an extremal function for
\(\overline S_N\) satisfying \(\sum_{z\in\mathbb Z^N}w(z)^{2^*}=1\).
Then \(w\in\ell^s(\mathbb Z^N)\) for every \(s>N/(N-2)\).
In particular, \(w\in\ell^2\) for \(N\geq5\), and
\(w\in\ell^{2+\varepsilon}\) for all \(\varepsilon>0\) when \(N=4\).
\end{lemma}

\begin{proof}
To apply the estimate \eqref{eq:green-hls}, we first prove that \(w=G*(aw)\), where \(a=\overline
S_Nw^{2^*-2}\).
The Euler--Lagrange equation is
\(\Delta_{\mathrm d}w=\overline S_N w^{2^*-1}\) on \(\mathbb Z^N\).
Let \(G\) be the Green function of \(\Delta_{\mathrm d}\), with
\(\Delta_{\mathrm d}G=\delta_0\). By \cite[Theorem~4.3.1]{LawlerLimic2010}, we have
\[
  0<G(z)\leq C(1+|z|)^{2-N},\qquad z\in\mathbb Z^N.
\]
Let \(1<r<N/2\), \(1/s=1/r-2/N\), and \(s'=s/(s-1)\).
For \(f\in\ell^r(\mathbb Z^N)\), set
\[
  F(z):=\sum_{\substack{w\in\mathbb Z^N\\w\neq z}}G(z-w)f(w),
  \qquad (G*f)(z)=G(0)f(z)+F(z).
\]
Since \(1/r+1/s'+(N-2)/N=2\), by
\cite[(1.4)]{HuaLi2021}, for every \(\varphi\in C_c(\mathbb Z^N)\), we have
\[
  \begin{aligned}
    \left|\sum_{z\in\mathbb Z^N}F(z)\varphi(z)\right|
    &\leq C\sum_{\substack{z,w\in\mathbb Z^N\\z\neq w}}
    \frac{|\varphi(z)||f(w)|}{|z-w|^{N-2}}\\
    &\leq C_s\|f\|_{\ell^r(\mathbb Z^N)}
    \|\varphi\|_{\ell^{s'}(\mathbb Z^N)}.
  \end{aligned}
\]
Taking the supremum over \(\varphi\in C_c(\mathbb Z^N)\) with
\(\|\varphi\|_{\ell^{s'}(\mathbb Z^N)}\leq1\), we obtain
\[
  \|F\|_{\ell^s(\mathbb Z^N)}
  =\sup_{\substack{\varphi\in C_c(\mathbb Z^N)\\
  \|\varphi\|_{\ell^{s'}(\mathbb Z^N)}\leq1}}
  \left|\sum_{z\in\mathbb Z^N}F(z)\varphi(z)\right|
  \leq C_s\|f\|_{\ell^r(\mathbb Z^N)}.
\]
Since \(r<s\), we also have
\(\|G(0)f\|_{\ell^s(\mathbb Z^N)}\leq G(0)\|f\|_{\ell^r(\mathbb Z^N)}\).
Hence,
\begin{equation}
  \label{eq:green-hls}
  \|G*f\|_{\ell^s(\mathbb Z^N)}\leq C_s\|f\|_{\ell^r(\mathbb Z^N)},
  \qquad1<r<N/2,\quad\frac1s=\frac1r-\frac2N.
\end{equation}
By \(\sum_{z\in\mathbb Z^N}w(z)^{2^*}=1\) and
\((2^*-2)N/2=2^*\), we have
\[
  \sum_{z\in\mathbb Z^N}|a(z)|^{N/2}
  =\overline S_N^{N/2}\sum_{z\in\mathbb Z^N}w(z)^{2^*}
  =\overline S_N^{N/2}<\infty.
\]
Since \(aw\in\ell^{(2^*)'}(\mathbb Z^N)\), where \((2^*)'=2^*/(2^*-1)\),
\eqref{eq:green-hls} gives \(G*(aw)\in\ell^{2^*}(\mathbb Z^N)\).
Choose \(f_j\in C_c(\mathbb Z^N)\) such that \(f_j\rightarrow aw\)
in \(\ell^{(2^*)'}(\mathbb Z^N)\) as \(j\rightarrow\infty\). By \eqref{eq:green-hls}, we have
\(G*f_j\rightarrow G*(aw)\) in
\(\ell^{2^*}(\mathbb Z^N)\)
and hence at every vertex. Letting \(j\rightarrow\infty\) in
\(\Delta_{\mathrm d}(G*f_j)=f_j\), we obtain \(\Delta_{\mathrm d}(G*(aw))=aw\).
Thus \(w-G*(aw)\) is harmonic and tends to zero at infinity.
By the maximum principle, we deduce that \(w=G*(aw)\).

To prove \(w\in\ell^s(\mathbb Z^N)\), fix \(N/(N-2)<s<2^*\)
and let \(K\subset\mathbb Z^N\) be a finite set to be chosen below.
By \(w=G*(aw)\), we have
\[
  \begin{aligned}
    w(z)
    &=\sum_{\xi\in\mathbb Z^N\setminus K}G(z-\xi)a(\xi)w(\xi)\\
    &\quad+\sum_{\xi\in K}G(z-\xi)a(\xi)w(\xi),
    \qquad z\in\mathbb Z^N.
  \end{aligned}
\]
Define \(Tf=G*(a\mathbf1_{K^c}f)\) and \(1/r_s=1/s+2/N\).
Since \(1<r_s<N/2\), by \eqref{eq:green-hls} and H\"older's inequality,
for \(f\in\ell^s(\mathbb Z^N)\), we have
\[
  \|Tf\|_s\leq C_s\|a\mathbf1_{K^c}f\|_{r_s}
  \leq C_s\|a\mathbf1_{K^c}\|_{N/2}\|f\|_s.
\]

By \eqref{eq:green-hls} and H\"older's inequality with \(s=2^*\), we have
\[
  \|Tf\|_{2^*}\leq C_{2^*}\|a\mathbf1_{K^c}\|_{N/2}\|f\|_{2^*}.
\]
Since \(a\in\ell^{N/2}(\mathbb Z^N)\), we choose \(K\)
such that
\[
  \max\{C_s,C_{2^*}\}\|a\mathbf1_{K^c}\|_{N/2}<1.
\]
Since \(a\mathbf1_Kw\) has finite support, we have
\(g:=G*(a\mathbf1_Kw)\in\ell^s(\mathbb Z^N)\cap\ell^{2^*}(\mathbb Z^N)\).
Since \(\|T\|_{\ell^s\to\ell^s}<1\) and
\(\|T\|_{\ell^{2^*}\to\ell^{2^*}}<1\),
\(\sum_{j\geq0}T^jg\) converges in \(\ell^s(\mathbb Z^N)\)
and \(\ell^{2^*}(\mathbb Z^N)\) to the same pointwise limit. By uniqueness for \((I-T)w=g\)
in \(\ell^{2^*}\), we have
\(w=\sum_{j\geq0}T^jg\). Hence, \(w\in\ell^s\) for \(N/(N-2)<s<2^*\).
For \(s\geq2^*\), the result follows from \(\ell^{2^*}\subset\ell^s\).
This completes the proof.
\end{proof}

 To prove \(\overline\kappa_{m,n}<\overline S_N\) for \(m\ne2\) and \(N\geq4\),
we cut off an extremal function for \(\overline S_N\) and translate  the
resulting function away from \(\Sigma\).

\begin{lemma}
\label{lem:strict-sobolev}
If \(N\geq4\) and \(m\ne2\), then \(\overline\kappa_{m,n}<\overline S_N\).
\end{lemma}
\begin{proof}
Let \(U\) be as in Lemma~\ref{lem:free-summability}.
For \(\chi\in C_c(\mathbb Z^N)\) with \(0\leq\chi\leq1\),
we first estimate
\[
  \sum_{z\in\mathbb Z^N}|D(\chi U)(z)|^2
  -\overline S_N\|\chi U\|_{2^*}^2.
\]
For this purpose, multiplying
\[
  \Delta_{\mathrm d}U(z)=\overline S_NU(z)^{2^*-1},
  \qquad z\in\mathbb Z^N,
\]
by \(\chi(z)^2U(z)\) and summing over \(z\in\mathbb Z^N\), we obtain
\begin{equation}
  \label{eq:cutoff-product}
  \sum_{z\in\mathbb Z^N}|D(\chi U)(z)|^2=\overline S_N\sum_{z\in\mathbb Z^N}\chi(z)^2U(z)^{2^*}
  +\sum_{\substack{\{z,w\}\subset\mathbb Z^N\\z\sim w}}U(z)U(w)|\chi(z)-\chi(w)|^2.
\end{equation}
Since \(\sum_zU(z)^{2^*}=1\) and \(2/2^*<1\), by Jensen's inequality, we have
\begin{equation}
  \label{eq:cutoff-jensen}
  \sum_{z\in\mathbb Z^N}\chi(z)^2U(z)^{2^*}
  =\sum_{z\in\mathbb Z^N}\bigl(\chi(z)^{2^*}\bigr)^{2/2^*}U(z)^{2^*}
  \leq\left(\sum_{z\in\mathbb Z^N}\chi(z)^{2^*}U(z)^{2^*}\right)^{2/2^*}
  =\|\chi U\|_{2^*}^2.
\end{equation}
Next, to apply \eqref{eq:cutoff-product} with \(\chi=\chi_R\), set
\(A(R)=\sum_{z\in B_R}U(z)^2\) and
\[
  D_R:=\sum_{\substack{\{z,w\}\subset\mathbb Z^N\\z\sim w}}
  U(z)U(w)|\chi_R(z)-\chi_R(w)|^2.
\]
For \(R\geq4\), if \(z\sim w\) and \(\chi_R(z)\neq\chi_R(w)\), then
\[
  R-1\leq|z|,|w|\leq2R+1,\qquad
  z,w\in B_{3R}\setminus B_{R/2},\qquad
  |\chi_R(z)-\chi_R(w)|\leq R^{-1}.
\]
By \(2U(z)U(w)\leq U(z)^2+U(w)^2\), we obtain
\[
  D_R\leq CR^{-2}\bigl(A(3R)-A(R/2)\bigr).
\]
For \(N\geq5\), the function \(A\) is bounded.
For \(N=4\), by Lemma~\ref{lem:free-summability} and H\"older's inequality,
we obtain
\[
  A(R)\leq |B_R|^{\varepsilon/(2+\varepsilon)}
  \|U\|_{2+\varepsilon}^2
  \leq C_\varepsilon R^{4\varepsilon/(2+\varepsilon)}
  \qquad(\varepsilon>0).
\]
Hence, for every \(\delta>0\), there exists \(C_\delta>0\) such that
\[
  A(R)\leq C_\delta R^\delta,\qquad R\geq1.
\]
To obtain \eqref{eq:cutoff-relative-bound}, it suffices to find
integers \(j_k\to\infty\) such that
\(A(6^{j_k+1})/A(6^{j_k})\to1\).
We argue indirectly. Then there exist \(\eta>0\) and an integer \(j_0\) such that
\[
  A(6^{j+1})\geq(1+\eta)A(6^j),\qquad j\geq j_0.
\]
Hence, for \(j>j_0\),
\[
  A(6^j)\geq(1+\eta)A(6^{j-1})
  \geq\cdots\geq(1+\eta)^{j-j_0}A(6^{j_0}).
\]
Choose \(\delta>0\) such that \(6^\delta<1+\eta\).
Combining this with \(A(6^j)\leq C_\delta6^{\delta j}\), we obtain
\[
  \frac{A(6^{j_0})}{(1+\eta)^{j_0}}
  \left(\frac{1+\eta}{6^\delta}\right)^j\leq C_\delta.
\]
This contradicts
\[
  \left(\frac{1+\eta}{6^\delta}\right)^j
  \xrightarrow[j\to\infty]{}\infty.
\]
With \(R_k=2\cdot6^{j_k}\), this gives, as \(k\rightarrow\infty\),
\begin{equation}
  \label{eq:cutoff-relative-bound}
  0\leq\frac{R_k^2D_{R_k}}{A(R_k)}
  \leq C\frac{A(6^{j_k+1})-A(6^{j_k})}{A(6^{j_k})}
  \rightarrow0.
\end{equation}
Finally, to obtain \(u_k\in C_c(\mathbb Z^N\setminus\Sigma)\), we translate \(\chi_{R_k}U\) in a \(y\)-direction.
Let \(e_1\) be the first coordinate vector in \(\mathbb Z^m\)
and define
\[
  u_k(x,y)=(\chi_{R_k}U)(x,y-3R_ke_1),
  \qquad (x,y)\in\mathbb Z^n\times\mathbb Z^m.
\]
By definition, \(\operatorname{supp}u_k\subset\{|y|\geq R_k\}\),
so \(u_k\in C_c(\mathbb Z^N\setminus\Sigma)\).
Moreover, \(u_k(0,3R_ke_1)=U(0)>0\).
By a change of variables, we have
\begin{equation}
  \label{eq:cutoff-hardy}
  \begin{aligned}
    \sum_{x\in\mathbb Z^n}\sum_{y\in\mathbb Z^m\setminus\{0\}}\frac{|u_k(x,y)|^2}{|y|^2}
    &=\sum_{(x,\eta)\in\operatorname{supp}(\chi_{R_k}U)}
    \frac{|(\chi_{R_k}U)(x,\eta)|^2}{|\eta+3R_ke_1|^2}\\
    &\geq\frac1{16R_k^2}\sum_{z\in B_{R_k}}U(z)^2
    =\frac{A(R_k)}{16R_k^2}.
  \end{aligned}
\end{equation}
Here \(\chi_{R_k}=1\) and \(|\eta+3R_ke_1|\leq4R_k\) on \(B_{R_k}\).
By the definition of \(u_k\),
\begin{equation}
  \label{eq:cutoff-translation}
  \sum_{z\in\mathbb Z^N}|Du_k(z)|^2
  =\sum_{z\in\mathbb Z^N}|D(\chi_{R_k}U)(z)|^2,\qquad
  \sum_{z\in\mathbb Z^N}|u_k(z)|^{2^*}
  =\sum_{z\in\mathbb Z^N}|\chi_{R_k}(z)U(z)|^{2^*}.
\end{equation}
By \eqref{eq:cutoff-product}, \eqref{eq:cutoff-jensen}, \eqref{eq:cutoff-hardy}, and
\eqref{eq:cutoff-translation}, we obtain
\[
  \begin{aligned}
    &\sum_{z\in\mathbb Z^N}|Du_k(z)|^2
    -c_m\sum_{x\in\mathbb Z^n}\sum_{y\in\mathbb Z^m\setminus\{0\}}\frac{|u_k(x,y)|^2}{|y|^2}\\
    &\leq\overline S_N\left(\sum_{z\in\mathbb Z^N}|u_k(z)|^{2^*}\right)^{2/2^*}
    +D_{R_k}-\frac{c_m A(R_k)}{16R_k^2}\\
    &<\overline S_N\left(\sum_{z\in\mathbb Z^N}|u_k(z)|^{2^*}\right)^{2/2^*}
  \end{aligned}
\]
for all large \(k\), since \(R_k^2D_{R_k}/A(R_k)\rightarrow0\) and \(c_m>0\).
By the definition of \(\overline\kappa_{m,n}\), for all large \(k\),
\[
  \overline\kappa_{m,n}\leq
  \frac{\displaystyle\sum_{z\in\mathbb Z^N}|Du_k(z)|^2
  -c_m\sum_{x\in\mathbb Z^n}\sum_{y\in\mathbb Z^m\setminus\{0\}}
  |u_k(x,y)|^2/|y|^2}
  {\displaystyle\left(\sum_{z\in\mathbb Z^N}|u_k(z)|^{2^*}\right)^{2/2^*}}
  <\overline S_N.
\]
This completes the proof.
\end{proof}

\subsection[Comparison of Green functions for N=3]{Comparison of Green functions for \(N=3\)}
 By \eqref{eq:general-sobolev-upper}, for \((m,n)=(1,2)\) it remains to prove
\(\overline\kappa_{1,2}\geq\overline S_3\).
For this purpose, we compare the Green function of \(H_+\) on \(D_+=\mathbb Z^2\times\mathbb N\)
with that of \(\Delta_{\mathrm d}\) on \(\mathbb Z^3\). Let
\(G_+\) be the minimal
Green function of
\[
  H_+=\Delta_2+L_1,
  \qquad
  (L_1a)_k=\left(2-\frac1{4k^2}\right)a_k-a_{k-1}-a_{k+1},
  \quad a_0=0.
\]
The operator \(\Delta_2\) is the discrete
Laplacian on \(\mathbb Z^2\), with the convention of Section~\ref{sec:preliminaries}.
Let \(g_3\) be the Green function of \(\Delta_{\mathrm d}\) on \(\mathbb Z^3\).
To prove \(G_+((x,k),(y,l))<g_3(x-y,k-l)\), we first establish the following lemmas.

To estimate \((e^{-tL_1})_{kl}\), set
\(\beta=1/\sqrt2\), \(r_k=k+(\beta-1)/2\), and define
\(B:\ell^2(\mathbb N)\to\ell^2(\mathbb N)\) by
\[
  B_{kk}=2-\frac{2\beta^2}{(2k+\beta-2)(2k+\beta)},
  \qquad
  B_{k,k+1}=B_{k+1,k}=-a_k,
  \quad
  a_k=\frac{(2k+\beta)^2-\beta^2}
  {(2k+\beta)\sqrt{(2k+\beta)^2-1}}.
\]

For \(|k-l|\geq2\), set \(B_{kl}=0\). Let \(v=2\beta/(\beta+2)\), \(\eta=v-1/4>0\), and
\(B_*=B+\eta P_1\), where \(P_1f=f(1)\delta_1\).
Let \(G_B\) and \(G_{B_*}\) be the Green functions of
\(\Delta_2+B\) and \(\Delta_2+B_*\),
respectively, on \(\mathbb Z^2\times\mathbb N\).

For \(t>0\) and \(d\in\mathbb Z\), set
\begin{equation}
  \label{eq:bessel-series}
  p(t,d)=e^{-2t}I_{|d|}(2t),\qquad
  I_\nu(z)=\sum_{j=0}^\infty
  \frac{(z/2)^{2j+\nu}}{j!\Gamma(j+\nu+1)}.
\end{equation}

For \(\Delta_{\mathrm d}\) acting on \(\mathbb Z\), we have
\begin{equation}
  \label{eq:free-heat-fourier}
  (e^{-t\Delta_{\mathrm d}})_{kl}=p(t,k-l)
  =\frac1{2\pi}\int_{-\pi}^{\pi}e^{-2t(1-\cos\theta)}e^{\mathrm i(k-l)\theta}\,d\theta.
\end{equation}

To prove \eqref{eq:green-order-use} for \(x\neq y\), we first establish
\[
  0<G_+<G_{B_*}\leq G_B.
\]

\begin{lemma}
\label{lem:jacobi-reference}
The operator \(B\) is nonnegative with spectrum \([0,4]\), and
\[
  0<(e^{-tL_1})_{kl}<(e^{-tB_*})_{kl}\leq(e^{-tB})_{kl}
  \qquad(t>0).
\]
If \(k\neq l\), then
\[
  (-\log B)_{kl}=\frac{2\sqrt{r_kr_l}}{|k-l|(r_k+r_l)},
\]
whereas
\[
  (-\log B)_{kk}
  =\int_0^\infty e^{-r_ks}
  \frac{2\cosh(\beta s/2)-1-\cosh(s/2)}{2\sinh(s/2)}\,ds<0.
\]

Moreover, \(0<G_+<G_{B_*}\leq G_B<\infty\).
\end{lemma}
\begin{proof}
We first prove  \(\sigma(B)=[0,4]\). Let
\[
  p_j(x)=\left(\frac{2j+\beta+1}{2^{\beta+1}}\right)^{1/2}
  P_j^{(0,\beta)}(x),\qquad j\geq0,
\]
where \(P_j^{(0,\beta)}\) is the Jacobi polynomial of degree \(j\), defined by the Rodrigues formula \cite[Table~18.5.1]{DLMF},
\[
  P_j^{(0,\beta)}(x)=\frac{(-1)^j}{2^j j!}(1+x)^{-\beta}
  \frac{d^j}{dx^j}\bigl[(1-x)^j(1+x)^{j+\beta}\bigr].
\] By
\cite[Section~1]{ACL2019}, we have
\[
  \int_{-1}^1p_i(x)p_j(x)(1+x)^\beta\,dx=\delta_{ij},\qquad
  \overline{\operatorname{span}\{p_j:j\geq0\}}^{\,L^2((-1,1),(1+x)^\beta dx)}
  =L^2((-1,1),(1+x)^\beta dx).
\]
By \cite[Sections~1--2]{ACL2019}, we have
\[
  \sum_{l\geq1}B_{kl}p_{l-1}(x)=2(1-x)p_{k-1}(x),
  \qquad k\geq1,\quad -1<x<1.
\]
Thus the map \(\mathcal U f=\sum_{k\geq1}f(k)p_{k-1}\) satisfies
\[
  \|\mathcal U f\|_{L^2((-1,1),(1+x)^\beta dx)}^2
  =\sum_{k\geq1}|f(k)|^2,\qquad
  (\mathcal U Bf)(x)=2(1-x)(\mathcal U f)(x),
\]
and \(\sigma(B)=[0,4]\).

Next, we prove \eqref{eq:semigroup-comparison-explicit}.
If \(q=2k+\beta\), then
\(a_k^2=1+[4q^2(q^2-1)]^{-1}>1\). Also
\((2k+\beta-2)(2k+\beta)<4k^2\). Hence, \(B_{kk}<(L_1)_{kk}\) for \(k\geq1\).
By the definition of \(\eta\),
\[
  (B_*)_{11}=(L_1)_{11},\qquad
  (B_*)_{kk}=B_{kk}<(L_1)_{kk}\quad(k\geq2).
\]

For \(k,l\geq1\),
\[
  0\leq(2I-L_1)_{kl}\leq(2I-B_*)_{kl}\leq(2I-B)_{kl}.
\]
For a bounded operator \(A\) on \(\ell^2(\mathbb N)\), we have
\begin{equation}
  \label{eq:operator-exponential}
  (e^{-tA})_{kl}=e^{-2t}\sum_{j=0}^\infty\frac{t^j}{j!}((2I-A)^j)_{kl}.
\end{equation}
Applying \eqref{eq:operator-exponential} to \(A=L_1,B_*,B\), we obtain
\((e^{-tL_1})_{kl}\leq(e^{-tB_*})_{kl}\leq(e^{-tB})_{kl}\).
For \(k<l\), we have
\[
  ((2I-B_*)^{l-k})_{kl}=\prod_{i=k}^{l-1}a_i>1
  =((2I-L_1)^{l-k})_{kl}.
\]
For \(k=l\geq2\), we have \((2I-B_*)_{kk}>(2I-L_1)_{kk}\).
For \(k=l=1\),
\[
  ((2I-B_*)^2)_{11}-((2I-L_1)^2)_{11}=a_1^2-1>0.
\]
Since \(A_{kl}=A_{lk}\) for \(A\in\{L_1,B_*,B\}\), by
\eqref{eq:operator-exponential}, for \(t>0\) and \(k,l\geq1\), we obtain
\begin{equation}
  \label{eq:semigroup-comparison-explicit}
  0<(e^{-tL_1})_{kl}<(e^{-tB_*})_{kl}\leq(e^{-tB})_{kl}.
\end{equation}

To prove \(T_{kl}>0\) in \eqref{eq:log-positive-gap}, let \(B_b\) be the operator
obtained by replacing \(\beta\) with \(b>0\)
in \(B\). The polynomials \(p_j\), now defined with parameter \(b\), satisfy
\begin{equation}
  \label{eq:jacobi-differential}
  -\bigl((1-x)(1+x)^{b+1}p_j'(x)\bigr)'
  =j(j+b+1)(1+x)^bp_j(x).
\end{equation}

By \eqref{eq:jacobi-differential} with indices \(i>j\), we have
\[
  -\bigl((1-x)(1+x)^{b+1}(p_jp_i'-p_ip_j')\bigr)'
  =(i-j)(i+j+b+1)(1+x)^bp_ip_j.
\]
Consequently, integration by parts gives
\[
  \begin{aligned}
    &(i-j)(i+j+b+1)
    \int_{-1}^1-\log(2(1-x))p_i p_j(1+x)^b\,dx\\
    &=\int_{-1}^1(1+x)^{b+1}(p_jp_i'-p_ip_j')\,dx\\
    &=2^{b+1}p_j(1)p_i(1).
  \end{aligned}
\]

By integration by parts
and \(i>j\), we have
\[
  \begin{aligned}
    &\int_{-1}^1(1+x)^{b+1}(p_jp_i'-p_ip_j')\,dx\\
    &=2^{b+1}p_i(1)p_j(1)
    -\int_{-1}^1p_i\bigl((b+1)p_j+2(1+x)p_j'\bigr)(1+x)^b\,dx
    =2^{b+1}p_i(1)p_j(1),
  \end{aligned}
\]
since
\[
  (b+1)p_j+2(1+x)p_j'\in\operatorname{span}\{p_0,\ldots,p_j\},
  \qquad \int_{-1}^1p_i p_\ell(1+x)^b\,dx=0\quad(0\leq\ell\leq j<i).
\]
By \(p_i,p_j,p_i',p_j'\in C([-1,1])\) and \(b>0\), we also have
\[
  \begin{aligned}
    (1-x)(1+x)^{b+1}\log(2(1-x))(p_jp_i'-p_ip_j')&\rightarrow0
    &&\text{as }x\rightarrow1^-\text{ or }x\rightarrow-1^+,\\
    (1+x)^{b+1}p_i(x)p_j(x)&\rightarrow0
    &&\text{as }x\rightarrow-1^+.
  \end{aligned}
\]
Since \(P_j^{(0,b)}(1)=1\), we obtain
\begin{equation}
  \label{eq:jacobi-log-off-diagonal}
  (-\log B_b)_{i+1,j+1}
  =\frac{\sqrt{(2i+b+1)(2j+b+1)}}{(i-j)(i+j+b+1)}
  \qquad(i>j).
\end{equation}

Since \((-\log B)_{kl}=(-\log B)_{lk}\), taking \(b=\beta\), \(k=i+1\), and \(l=j+1\)
in \eqref{eq:jacobi-log-off-diagonal}, we obtain
\[
  (-\log B)_{kl}=\frac{2\sqrt{r_kr_l}}{|k-l|(r_k+r_l)},\qquad k\neq l.
\]
Since
\(\left.\partial_\alpha(1-x)^\alpha\right|_{\alpha=0}=\log(1-x)\),
we introduce \(\alpha\) to calculate \((-\log B_b)_{j+1,j+1}\).
Define
\begin{equation}
  \label{eq:jacobi-hj-definition}
  h_j:=\int_{-1}^1\bigl(P_j^{(\alpha,b)}(x)\bigr)^2(1-x)^\alpha(1+x)^b\,dx,
  \qquad P_j^{(\alpha,b)}(x)=c_jx^j+\sum_{i=0}^{j-1}c_{ji}x^i.
\end{equation}
Then
\[
  h_j=\frac{2^{\alpha+b+1}\Gamma(j+\alpha+1)\Gamma(j+b+1)}
  {(2j+\alpha+b+1)j!\Gamma(j+\alpha+b+1)},
  \qquad
  c_j=\frac{\Gamma(2j+\alpha+b+1)}{2^j j!\Gamma(j+\alpha+b+1)}.
\]

Since
\[
  \partial_\alpha P_j^{(\alpha,b)}
  -\frac{\partial_\alpha c_j}{c_j}P_j^{(\alpha,b)}
  =\sum_{i=0}^{j-1}\left(\partial_\alpha c_{ji}
  -\frac{\partial_\alpha c_j}{c_j}c_{ji}\right)x^i
  \in\operatorname{span}\{P_0^{(\alpha,b)},\ldots,P_{j-1}^{(\alpha,b)}\}
\]
and the Jacobi polynomials satisfy
\[
  \int_{-1}^1P_j^{(\alpha,b)}(x)P_i^{(\alpha,b)}(x)
  (1-x)^\alpha(1+x)^b\,dx=0,
  \qquad 0\leq i<j,
\]
we have
\[
  \int_{-1}^1P_j^{(\alpha,b)}\partial_\alpha P_j^{(\alpha,b)}
  (1-x)^\alpha(1+x)^b\,dx=\frac{\partial_\alpha c_j}{c_j}h_j.
\]
Differentiating \eqref{eq:jacobi-hj-definition} with respect to \(\alpha\), we obtain
\[
  \partial_\alpha h_j=
  \int_{-1}^1\log(1-x)\bigl(P_j^{(\alpha,b)}(x)\bigr)^2
  (1-x)^\alpha(1+x)^b\,dx
  +2\frac{\partial_\alpha c_j}{c_j}h_j.
\]
Consequently, with \(\psi=\Gamma'/\Gamma\),
\[
  \begin{aligned}
    (-\log B_b)_{j+1,j+1}
    &=-\log2-\int_{-1}^1\log(1-x)p_j(x)^2(1+x)^b\,dx\\
    &=-\log2-\left[\frac{\partial_\alpha h_j}{h_j}
    -2\frac{\partial_\alpha c_j}{c_j}\right]_{\alpha=0}\\
    &=-2\log2-\psi(j+1)-\psi(j+b+1)
    +2\psi(2j+b+1)+\frac1{2j+b+1}.
  \end{aligned}
\]

Let \(r=j+(b+1)/2\). By
\[
  2\psi(2r)=\psi(r)+\psi(r+1/2)+2\log2,\qquad
  \psi(v)-\psi(u)=\int_0^\infty\frac{e^{-us}-e^{-vs}}{1-e^{-s}}\,ds
  \quad(u,v>0),
\]
we obtain
\[
  (-\log B_b)_{j+1,j+1}
  =\int_0^\infty e^{-rs}
  \frac{2\cosh(bs/2)-1-\cosh(s/2)}{2\sinh(s/2)}\,ds.
\]

For \(|\alpha|\leq1/2\) and fixed \(j\geq0\) and \(b>0\), the differentiated integrands satisfy
\[
  \left|\partial_\alpha\left[(P_j^{(\alpha,b)}(x))^2(1-x)^\alpha(1+x)^b\right]\right|
  \leq C_{j,b}(1-x)^{-1/2}(1+x)^b(1+|\log(1-x)|)\in L^1(-1,1).
\]
Thus we may differentiate \eqref{eq:jacobi-hj-definition} with respect to \(\alpha\) by Lebesgue's dominated
convergence theorem.
For \(b=\beta=1/\sqrt2\) and \(s>0\),
\[
  2\cosh(\beta s/2)-1-\cosh(s/2)
  =\sum_{j=2}^\infty(2\beta^{2j}-1)\frac{(s/2)^{2j}}{(2j)!}<0.
\]
Hence, \((-\log B)_{kk}<0\) for every \(k\geq1\).
By
\[
  \log(1-2\rho\cos\theta+\rho^2)
  =-2\sum_{j=1}^\infty\frac{\rho^j\cos(j\theta)}j,
  \qquad0<\rho<1,
\]
and convergence in \(L^1(-\pi,\pi)\) as \(\rho\uparrow1\), we obtain
\[
  (\log\Delta_{\mathrm d})_{kl}
  =\frac1{2\pi}\int_{-\pi}^\pi e^{i(k-l)\theta}\log(2-2\cos\theta)\,d\theta
  =\begin{cases}0,&k=l,\\-1/|k-l|,&k\neq l.\end{cases}
\]
Consequently,
\begin{equation}
  \label{eq:log-positive-gap}
  \begin{aligned}
    T_{kl}&:=\int_0^\infty\frac{p(t,k-l)-(e^{-tB})_{kl}}{t}\,dt\\
    &=\begin{cases}
      \displaystyle\frac{(\sqrt{r_k}-\sqrt{r_l})^2}{|k-l|(r_k+r_l)},&k\neq l,\\
      -(-\log B)_{kk},&k=l,
    \end{cases}
    \quad>0.
  \end{aligned}
\end{equation}
By \eqref{eq:operator-exponential} and \eqref{eq:free-heat-fourier},
\(p(t,k-l)-(e^{-tB})_{kl}=O(t)\) as \(t\rightarrow0^+\).
For fixed \(k,l\),
\(\max_{-1\leq x\leq1}|p_{k-1}(x)p_{l-1}(x)(1+x)^\beta|<\infty\), so
\begin{equation}
  \label{eq:jacobi-heat-integral}
  (e^{-tB})_{kl}
  =\int_{-1}^1e^{-2t(1-x)}p_{k-1}(x)p_{l-1}(x)(1+x)^\beta\,dx
  =O(t^{-1})\quad(t\rightarrow\infty).
\end{equation}
By \eqref{eq:free-heat-fourier}, \(p(t,d)\leq C(1+t)^{-1/2}\) uniformly in
\(d\in\mathbb Z\).
Hence, the integral in \eqref{eq:log-positive-gap} converges absolutely, since
\[
  \left|\frac{p(t,k-l)-(e^{-tB})_{kl}}t\right|
  \leq C_{k,l}\begin{cases}1,&0<t\leq1,\\t^{-3/2}+t^{-2},&t\geq1.\end{cases}
\]
Finally, we compare \(G_+\), \(G_{B_*}\), and \(G_B\).
For \(A\in\{L_1,B_*,B\}\), let \(G_{L_1}=G_+\). By
\eqref{eq:semigroup-comparison-explicit} and \eqref{eq:jacobi-heat-integral},
\[
  \begin{aligned}
    G_A((X,k),(0,l))
    &=\int_0^\infty p(t,X_1)p(t,X_2)(e^{-tA})_{kl}\,dt\\
    &\leq C_{k,l}\int_0^\infty(1+t)^{-2}\,dt<\infty,
    \qquad X\in\mathbb Z^2.
  \end{aligned}
\]
Integrating \eqref{eq:semigroup-comparison-explicit}, we obtain
\(0<G_+<G_{B_*}\leq G_B<\infty\).
This completes the proof.
\end{proof}

To prove \(G_B((X,k),(0,l))<g_3(X,k-l)\) for \(X\ne0\),
we improve \eqref{eq:log-positive-gap} by proving
\[
  T_{kl}(s):=\int_s^\infty
  \frac{p(t,k-l)-(e^{-tB})_{kl}}{t}\,dt>0,
  \qquad s\geq0,\quad k,l\geq1.
\]

\begin{lemma}
\label{lem:heat-crossing}
For each \(k,l\geq1\), \((e^{-tB})_{kl}-p(t,k-l)\) changes sign exactly
once for \(t>0\), from positive to negative. For every \(s\geq0\),
\[
  T_{kl}(s):=\int_s^\infty
  \frac{p(t,k-l)-(e^{-tB})_{kl}}{t}\,dt>0.
\]
For every \(X\in\mathbb Z^2\setminus\{0\}\),
\(G_B((X,k),(0,l))<g_3(X,k-l)\).
\end{lemma}
\begin{proof}
To determine the sign of \((e^{-tB})_{kl}-p(t,k-l)\), we first study the coefficients of
\(e^{4t}[(e^{-tB})_{kl}-p(t,k-l)]\).
Since \(B_{kl}=B_{lk}\) and \(p(t,k-l)=p(t,l-k)\), we may assume that \(k\leq l\). Set \(\mu=k-1\), \(d=l-k\),
\(\nu=l-1\),
\(b=1+\beta\) and \(S=k+l+\beta-1\). With \(M_B=4I-B\) and
\(M_0=4I-\Delta_{\mathrm d}\), set
\[
  F_j=(M_0^j)_{kl}=\binom{2j}{j-d},\qquad
  c_j=\frac{(M_B^j)_{kl}}{F_j}\quad(j\geq d).
\]

By
\[
  (2+z+z^{-1})^j=z^{-j}(1+z)^{2j}
  =\sum_{a=0}^{2j}\binom{2j}{a}z^{a-j},
\]
we have \(F_j=\binom{2j}{j-d}\).
If \(x_0=k\), \(x_j=l\), \(|x_i-x_{i-1}|\leq1\), and \(x_a=0\), then
\(a\geq k\), \(j-a\geq l\), and hence \(j\geq k+l\).
Thus, for \(\max\{d,1\}\leq j<k+l\),
\[
  \begin{aligned}
    (M_B^j)_{kl}
    &=\sum_{x_1,\ldots,x_{j-1}\geq1}
    \prod_{i=1}^j(M_B)_{x_{i-1},x_i}
    \geq(M_0^j)_{kl},\\
    (M_B^{d+1})_{kl}&>(M_0^{d+1})_{kl},
  \end{aligned}
\]
by \((M_B)_{ii}>2\) and \((M_B)_{i,i+1}=a_i>1\).
When \(d=0\), we also have \(c_0=1\). Hence, \(c_j\geq1\) for \(d\leq j<k+l\) and \(c_{d+1}>1\).
For \(j<d\), \((M_B^j)_{kl}=(M_0^j)_{kl}=0\).

We next prove  \(c_j\rightarrow0\). By
\(\mathcal U M_B f=2(1+x)\mathcal U f\), with \(y=(1+x)/2\), we have
\begin{equation}
  \label{eq:jacobi-moments}
  (M_B^j)_{kl}
  =4^j\sqrt{(2\mu+b)(2\nu+b)}
  \int_0^1y^{j+b-1}P_\mu^{(0,b-1)}(2y-1)P_\nu^{(0,b-1)}(2y-1)\,dy.
\end{equation}

Since \(\max_{0\leq y\leq1}|P_\mu^{(0,b-1)}(2y-1)P_\nu^{(0,b-1)}(2y-1)|<\infty\),
\[
  0<(M_B^j)_{kl}\leq C_{k,l}4^j\int_0^1y^{j+b-1}\,dy
  =\frac{C_{k,l}4^j}{j+b}.
\]
Together with \(F_j\geq c_d4^j(j+1)^{-1/2}\) for all large \(j\), this gives
\(c_j=O(j^{-1/2})\to0\).
To prove that \(c_j-1\) changes sign once, we first derive
\eqref{eq:moment-recurrence} from the Jacobi differential equation
\eqref{eq:jacobi-differential}. Set
\[
  U(y)=P_\mu^{(0,b-1)}(2y-1),\qquad
  V(y)=P_\nu^{(0,b-1)}(2y-1),
\]
and
\[
  \lambda_j=j(j+b),\qquad L=\lambda_\mu+\lambda_\nu,\qquad
  \delta=\lambda_\nu-\lambda_\mu.
\]
By \eqref{eq:jacobi-differential}, we have
\begin{equation}
  \label{eq:jacobi-pair-ode}
  \begin{aligned}
    y(1-y)U''+[b-(b+1)y]U'+\lambda_\mu U&=0,\\
    y(1-y)V''+[b-(b+1)y]V'+\lambda_\nu V&=0,
  \end{aligned}
  \qquad 0<y<1.
\end{equation}

For integers \(r\geq0\), define
\[
  M_r=\int_0^1y^{r+b-1}UV\,dy,\qquad
  J_r=\int_0^1y^{r+b}(1-y)U'V'\,dy.
\]
By \eqref{eq:jacobi-pair-ode} and integration by parts, for \(r\geq1\) we obtain
\begin{equation}
  \label{eq:jacobi-pair-integrals}
  \begin{aligned}
    2J_r&=r(r+b-1)M_{r-1}+[L-r(r+b)]M_r,\\
    \int_0^1y^{r+b}(1-y)(U'V-UV')\,dy
    &=-\frac{\delta}{r+1}M_{r+1}.
  \end{aligned}
\end{equation}
By \eqref{eq:jacobi-pair-ode}, we also have
\begin{equation}
  \label{eq:jacobi-derivative-product}
  2y(1-y)(U'V')'+4[b-(b+1)y]U'V'+L(UV)'
  =-\delta(U'V-UV').
\end{equation}
Multiplying \eqref{eq:jacobi-derivative-product} by \(y^{r+b}(1-y)\),
integrating by parts, and using \eqref{eq:jacobi-pair-integrals}, we obtain
\begin{equation}
  \label{eq:jacobi-derivative-moments}
  2(b-1-r)(J_r-J_{r+1})-L(r+b)M_r+L(r+b+1)M_{r+1}
  =\frac{\delta^2}{r+1}M_{r+1}.
\end{equation}

By \eqref{eq:jacobi-moments}, we have
\[
  c_r=\frac{4^r\sqrt{(2\mu+b)(2\nu+b)}}{F_r}M_r,\qquad r\geq d.
\]
For \(r\geq k+l-1\), by the definition of \(F_r\), we have
\begin{equation}
  \label{eq:free-moment-ratios}
  \frac{F_{r+1}}{4F_r}
  =\frac{(r+1)(2r+1)}{2((r+1)^2-d^2)},\qquad
  \frac{4F_{r-1}}{F_r}
  =\frac{2(r^2-d^2)}{r(2r-1)}.
\end{equation}
Substituting \eqref{eq:jacobi-pair-integrals} and \eqref{eq:free-moment-ratios}
into \eqref{eq:jacobi-derivative-moments}, with
\[
  L=\frac{S^2+d^2-b^2}{2},\qquad \delta=dS,
\]
we obtain
\begin{equation}
  \label{eq:moment-recurrence}
  A_r(c_{r+1}-c_r)-C_r(c_r-c_{r-1})+E_rc_r=0,
\end{equation}
where
\[
  A_r=\frac{(r+1)^2-S^2}{4},\quad
  C_r=\frac{(r^2-d^2)(r^2-\beta^2)}{4r^2-1},\quad
  E_r=\frac{r^2+(4\beta^2-1)d^2-\beta^2}{4(4r^2-1)}.
\]
For integers \(r\geq R=k+l-1\),
since \(S=R+\beta<R+1\), \(r>d\),
and \(\beta^2=1/2\), we have \(A_r>0\), \(C_r>0\), and \(E_r>0\).
By \eqref{eq:moment-recurrence}, we have
\[
  c_{r+1}-c_r=\frac{C_r}{A_r}(c_r-c_{r-1})-\frac{E_r}{A_r}c_r.
\]
Since \(c_r>0\),
\[
  c_r-c_{r-1}\leq0\quad\Longrightarrow\quad
  c_{r+1}-c_r<0\quad\Longrightarrow\quad c_{s+1}-c_s<0\quad(s\geq r\geq R).
\]
Such an \(r\) exists because \(c_r\to0\).
Together with \(c_r\geq1\)
for \(d\leq r\leq R\), this gives an integer
\(J=\max\{r:c_r>1\}\) such that
\[
  c_r\geq1\quad(d\leq r\leq J),\qquad c_r\leq1\quad(r>J).
\]

To prove that \((e^{-tB})_{kl}-p(t,d)\) changes sign exactly once, set
\begin{equation}
  \label{eq:heat-difference-series}
 H(t)=\sum_{r=d}^\infty\frac{F_r(c_r-1)}{r!}t^{r-J},\qquad
 (e^{-tB})_{kl}-p(t,d)=e^{-4t}t^JH(t).
\end{equation}
For \(A=\max\{1,\|M_B\|,\|M_0\|\}\), we have \(|F_r(c_r-1)|\leq2A^r\), and hence
\[
 \sum_{r=d}^\infty\frac{|F_r(c_r-1)|}{r!}t^{r-J}
 \leq2t^{-J}e^{At}<\infty,\qquad t>0.
\]
Choose \(K>J\) with \(c_K<1\). Since \(c_r\geq1\) for \(d\leq r\leq J\) and \(c_r\leq1\) for \(r>J\), for \(0<t_1<t_2\) we obtain
\[
 \begin{aligned}
 H(t_2)-H(t_1)
 &=\sum_{r=d}^\infty\frac{F_r(c_r-1)}{r!}
   \bigl(t_2^{r-J}-t_1^{r-J}\bigr)\\
 &\leq\frac{F_K(c_K-1)}{K!}
   \bigl(t_2^{K-J}-t_1^{K-J}\bigr)<0.
 \end{aligned}
\]
Let \(j_0=\min\{r:c_r>1\}\). Then
\[
 \begin{aligned}
 H(t)&\sim\frac{F_{j_0}(c_{j_0}-1)}{j_0!}t^{j_0-J}>0
 &&\text{as }t\to0^+,\\
 H(t)&\leq\sum_{r=d}^J\frac{F_r(c_r-1)}{r!}
 +\frac{F_K(c_K-1)}{K!}t^{K-J}\to-\infty
 &&\text{as }t\to\infty.
 \end{aligned}
\]
By \eqref{eq:heat-difference-series}, \(H\) is continuous on \((0,\infty)\). Since \(H\) is strictly decreasing, there is a unique \(\tau>0\) with \(H(\tau)=0\). Thus, for \(D(t)=p(t,d)-(e^{-tB})_{kl}\),
\[
 D(t)<0\quad(0<t<\tau),\qquad D(t)>0\quad(t>\tau).
\]
By \eqref{eq:log-positive-gap}, for \(0\leq s\leq\tau\),
\[
  T_{kl}(s)=T_{kl}(0)-\int_0^s\frac{D(t)}t\,dt
  \geq T_{kl}(0)>0.
\]
For \(s>\tau\),
\(T_{kl}(s)=\int_s^\infty D(t)t^{-1}\,dt>0\).

Finally, to compare \(G_B((X,k),(0,l))\) with \(g_3(X,k-l)\)
using \(T_{kl}(t)>0\), set
\(M_X(t)=t p(t,X_1)p(t,X_2)\) for
\(X=(X_1,X_2)\in\mathbb Z^2\setminus\{0\}\).
We will integrate by parts, so we first prove  \(M_X'(t)>0\)
and verify that \(M_X(t)T_{kl}(t)\to0\) as \(t\to0^+\) and \(t\to\infty\).

To prove \(M_X'(t)>0\), we first compare \(zI_j'(z)/I_j(z)\) for integers \(j\geq0\). For integers \(j>i\geq0\), by \cite[(10.25.1)]{DLMF}, we have
\[
 \frac{d}{dz}\left[z\bigl(I_j'(z)I_i(z)-I_j(z)I_i'(z)\bigr)\right]
 =\frac{j^2-i^2}{z}I_j(z)I_i(z)>0,\qquad z>0.
\]
By \eqref{eq:bessel-series},
\[
 z\bigl(I_j'(z)I_i(z)-I_j(z)I_i'(z)\bigr)\to0
 \quad\text{as }z\to0^+.
\]
Hence
\[
 \frac{zI_j'(z)}{I_j(z)}>\frac{zI_i'(z)}{I_i(z)},\qquad z>0.
\]
By \cite[(10.29.1)]{DLMF},
\begin{equation}
  \label{eq:bessel-recurrences}
  I_0'=I_1,\qquad I_1'=I_0-I_1/z,\qquad
  I_{j-1}(z)-I_{j+1}(z)=2jI_j(z)/z.
\end{equation}
Let \(z=2t\) and \(r=I_1(z)/I_0(z)\). By \eqref{eq:free-heat-fourier},
\[
  I_0(z)-I_1(z)=\frac1\pi\int_0^\pi e^{z\cos\theta}(1-\cos\theta)\,d\theta>0.
\]
Thus \(0<r<1\). Since \(X\ne0\), at least one of \(|X_1|,|X_2|\) is at least one. Therefore,
\begin{equation}
  \label{eq:MX-increasing}
 \begin{aligned}
 t(\log M_X)'&=1-2z+\frac{zI_{|X_1|}'(z)}{I_{|X_1|}(z)}
 +\frac{zI_{|X_2|}'(z)}{I_{|X_2|}(z)}\\
 &\geq1-2z+\frac{zI_1'(z)}{I_1(z)}+\frac{zI_0'(z)}{I_0(z)}\\
 &=z(r+r^{-1}-2)>0.
 \end{aligned}
\end{equation}
By \(p(t,d)\leq C(1+t)^{-1/2}\) and \eqref{eq:bessel-series},
\[
  0<M_X(t)\leq\frac{Ct}{1+t},\qquad \lim_{t\to0^+}M_X(t)=0.
\]
By \(T_{kl}'(t)=-D(t)/t\), integration by parts
gives
\[
  g_3(X,k-l)-G_B((X,k),(0,l))
  =\int_0^\infty M_X'(t)T_{kl}(t)\,dt>0.
\]

Indeed, by \eqref{eq:jacobi-heat-integral} and
\(p(t,d)\leq C(1+t)^{-1/2}\), we have
\[
  |T_{kl}(t)|\leq C_{k,l}\int_t^\infty(s^{-3/2}+s^{-2})\,ds
  =O(t^{-1/2})\qquad(t\to\infty).
\]
Since \(T_{kl}(t)\to T_{kl}(0)\) as \(t\to0^+\), we obtain
\[
  \lim_{t\rightarrow0^+}M_X(t)T_{kl}(t)=0
  \quad\text{and}\quad
  \lim_{t\rightarrow\infty}M_X(t)T_{kl}(t)=0.
\]
This completes the proof.
\end{proof}

Lemma~\ref{lem:heat-crossing} gives
\(G_+((X,k),(0,l))<g_3(X,k-l)\) for \(X\ne0\).
For \(X=0\), we shall use the following estimates for \(p(t,0)\) and \(g_3(0)\).

Set \(c=1/(4\pi)\),
\(M_0(t)=t p(t,0)^2\), \(M_e(t)=t p(t,0)p(t,1)\), and
\(g_d=g_3(0,0,d)\). We shall use
\begin{equation}
  \label{eq:scalar-green-bounds}
  M_0(t)<c\left(1+\frac1{4t}\right),\qquad
  M_0(t)>c\ (t\geq1/4),\qquad
  \frac29<g_0<\frac4{15}.
\end{equation}
When used as a vertex of \(\mathbb Z^2\), \(e\) denotes \((1,0)\).

We
first prove \(M_0(t)<c(1+1/(4t))\). Define
\(W_j=\Gamma(j+1/2)/(\sqrt\pi\Gamma(j+1))\).

By \cite{ChenQi2002},
\begin{equation}
  \label{eq:wallis-bound}
  W_j^2<\frac1{\pi(j+1/4)},\qquad j\geq1.
\end{equation}
By \eqref{eq:bessel-series} and
\(\sum_{a=0}^j\binom ja^2=\binom{2j}j\), we have
\begin{equation}
  \label{eq:bessel-square-series}
  I_0(z)^2=\sum_{j=0}^\infty\left(\sum_{a=0}^j\frac1{(a!)^2((j-a)!)^2}\right)(z/2)^{2j}
  =\sum_{j=0}^\infty(2j+1)W_j^2\frac{(2z)^{2j}}{(2j+1)!}.
\end{equation}
For \(j=1\), \(3W_1^2=3/4<2/\pi+1/8\), since \(\pi<16/5\).
For \(j\geq2\),
by \eqref{eq:wallis-bound},
\begin{equation}
  \label{eq:wallis-coefficients}
  (2j+1)W_j^2<\frac{2j+1}{\pi(j+1/4)}
  =\frac2\pi+\frac1{2\pi(j+1/4)}
  <\frac2\pi+\frac1{4(j+1)}.
\end{equation}
By \eqref{eq:bessel-square-series}, \eqref{eq:wallis-coefficients}, and \(3W_1^2<2/\pi+1/8\), we obtain
\[
  I_0(z)^2<\frac2\pi\frac{\sinh(2z)}{2z}
  +\frac{\cosh(2z)-1}{8z^2}+\frac34-\frac2\pi.
\]
It follows that
\[
  2\pi ze^{-2z}I_0(z)^2
  <1+\frac\pi{8z}+\left(\frac{3\pi}2-4\right)ze^{-2z}.
\]
Since \(z^2e^{-2z}\leq e^{-2}\), \(\pi<22/7\) and \(e>8/3\) imply
\[
  \frac\pi8+(3\pi/2-4)e^{-2}<\frac{221}{448}<\frac12.
\]
Substitute \(z=2t\) to obtain the first estimate in
\eqref{eq:scalar-green-bounds}.

To prove \(M_0(t)>c\) for \(t\geq1/4\), set \(q(z)=\sqrt z e^{-z}I_0(z)\).
By \(z^2I_0''+zI_0'-z^2I_0=0\),
\[
  q''+2q'+\frac{q}{4z^2}=0,\qquad (e^{2z}q')'=-\frac{e^{2z}q}{4z^2}<0.
\]
By \eqref{eq:bessel-series} and \eqref{eq:bessel-recurrences},
\[
  q'(z)\sim\frac1{2\sqrt z}>0\quad(z\to0^+),\qquad
  2I_1(1)=\sum_{j=0}^\infty\frac1{4^jj!(j+1)!}
  <\sum_{j=0}^\infty\frac1{4^j(j!)^2}=I_0(1),
\]
so \(q'(1)<0\).
Thus there is a unique \(z_0\in(0,1)\) such that
\(q'>0\) on \((0,z_0)\) and \(q'<0\) on \((z_0,\infty)\).
By \cite[(10.40.1)]{DLMF} with \(\nu=0\), we have
\[
  q(z)=\sqrt z\,e^{-z}I_0(z)\longrightarrow\frac1{\sqrt{2\pi}}
  \qquad(z\to\infty).
\]
Since \(q(1/2)^2>1/(2e)>1/(2\pi)\), the
signs of \(q'\) give
\(q(z)^2>1/(2\pi)\) for \(z\geq1/2\). Hence
\(M_0(t)=q(2t)^2/2>1/(4\pi)\) for \(t\geq1/4\).

It remains to prove \(2/9<g_0<4/15\). By \cite[Section~1.2]{Guttmann2010}, with
\[
  6g_0=\frac1{(2\pi)^3}\int_{[-\pi,\pi]^3}
  \frac{d\theta}{1-(\cos\theta_1+\cos\theta_2+\cos\theta_3)/3},
\]
we have
\[
  6g_0=\frac{1-9q^2}{D}
  \left(\frac{2K(\sqrt\mu)}{\pi}\right)^2,
  \quad q=1-\frac{2\sqrt2}{3},\quad\xi=\sqrt q,
  \quad D=(1-\xi)^3(1+3\xi),\quad\mu=\frac{16\xi^3}{D}.
\]

Here \(K(a):=\int_0^{\pi/2}(1-a^2\sin^2\theta)^{-1/2}\,d\theta\).
By \(113/80<\sqrt2<17/12\), we have
\(1/18<q<7/120<1/16\). Also
\(D(D-16\xi^3)=(1-q)^3(1-9q)>0\), so \(0<\mu<1\).

By the definition of \(K\),
\[
  \frac{2K(\sqrt\mu)}\pi=\sum_{j=0}^\infty W_j^2\mu^j,
  \qquad W_j^2\leq\frac{(1/4)_j}{j!},
\]
where \(W_j^2\leq(1/4)_j/j!\) follows from
\((j+1/2)^2\leq(j+1/4)(j+1)\). Hence,
\begin{equation}
  \label{eq:elliptic-K-bounds}
  1+\frac\mu4<\frac{2K(\sqrt\mu)}\pi<(1-\mu)^{-1/4}.
\end{equation}

By \(2K(\sqrt\mu)/\pi<(1-\mu)^{-1/4}\) in \eqref{eq:elliptic-K-bounds}, we obtain
\[
  (6g_0)^2<\frac{(1-9q^2)^2}{(1-q)^3(1-9q)}<\frac{64}{25},
\]
since \((1-q)^3(1-9q)>(15/16)^3(19/40)>25/64\).

To prove \(g_0>2/9\), we use \(D<1\) and
\[
  4(1-9q^2)-5D=-1+30q-40q\sqrt q-21q^2
  >-1+20q-21q^2>0.
\]
Thus \((1-9q^2)/D>5/4\); combining this with
\(\xi^3>1/81\) and \(\mu>16/81\), we obtain
\(6g_0>(5/4)(1+\mu/2)>445/324>4/3\).
This proves all three estimates in \eqref{eq:scalar-green-bounds}.

To estimate \(G_B-G_{B_*}\) in Proposition~\ref{prop:green-order}, we also need the following integrals:
\begin{equation}
  \label{eq:inverse-moments}
  \int_0^\infty\frac{p(t,j)}t\,dt=\frac1j,
  \quad
  \int_0^\infty\frac{p(t,j)}{t^2}\,dt=\frac{2}{j(j^2-1)},
  \quad
  \int_0^\infty\frac{p(t,j)}{t^3}\,dt
  =\frac{12}{j(j^2-1)(j^2-4)},
\end{equation}
valid for \(j>0\), \(j>1\), and \(j>2\), respectively.
To prove \eqref{eq:inverse-moments}, define \(D\) on \(\ell^2(\mathbb N)\) by
\[
  (Df)_j=2f_j-f_{j-1}-f_{j+1},\qquad j\geq1,\qquad f_0=0.
\]
For \(a>0\), direct substitution gives
\[
  \bigl((D+aI)^{-1}\bigr)_{j1}=q(a)^j,
  \qquad q(a)=\frac{a+2-\sqrt{a(a+4)}}2.
\]

By \eqref{eq:free-heat-fourier},
\begin{equation}
  \label{eq:half-line-heat}
  (e^{-tD})_{kl}=p(t,k-l)-p(t,k+l),\qquad k,l\geq1.
\end{equation}
Taking \(l=1\) and using \eqref{eq:bessel-recurrences}, we obtain
\begin{equation}
  \label{eq:inverse-moment-transform}
  \int_0^\infty e^{-at}\frac{p(t,j)}t\,dt
  =\frac1j\int_0^\infty e^{-at}
  \bigl(p(t,j-1)-p(t,j+1)\bigr)\,dt
  =\frac{q(a)^j}{j}.
\end{equation}

Since \(e^{-at}p(t,j)/t\uparrow p(t,j)/t\) for every \(t>0\)
and \(q(a)\to1\) as \(a\downarrow0\), by Lebesgue's monotone convergence theorem, we obtain
\(\int_0^\infty p(t,j)t^{-1}\,dt=1/j\).
By \eqref{eq:inverse-moment-transform},
\(\int_0^\infty e^{-at}\,da=t^{-1}\), and
\(\int_0^\infty ae^{-at}\,da=t^{-2}\), Tonelli's theorem gives
\[
  \begin{aligned}
    \int_0^\infty\frac{p(t,j)}{t^2}\,dt
    &=\frac1j\int_0^\infty q(a)^j\,da
    =\frac1j\int_0^1(q^{j-2}-q^j)\,dq
    =\frac2{j(j^2-1)},\\
    \int_0^\infty\frac{p(t,j)}{t^3}\,dt
    &=\frac1j\int_0^\infty aq(a)^j\,da
    =\frac1j\int_0^1q^{j-3}(1-q)^3(1+q)\,dq\\
    &=\frac{12}{j(j^2-1)(j^2-4)},
  \end{aligned}
\]
where the changes of variable use \(a=q+q^{-1}-2\).
The conditions \(j>1\) and \(j>2\), respectively, ensure
integrability at \(q=0\).

To treat \(k=1\) in Proposition~\ref{prop:green-order},
we first compare \(g_j\) and \(g_{j-1}\).

\begin{lemma}
\label{lem:free-green-ratio}
For every integer \(j\geq1\), we have
\[
  g_j>\frac14g_{j-1}.
\]
\end{lemma}
\begin{proof}
Set \(\omega(\theta_1,\theta_2)=4\sin^2(\theta_1/2)+4\sin^2(\theta_2/2)\).
By \eqref{eq:free-heat-fourier}, integrating with respect to \(\theta_3\), we obtain
\begin{equation}
  \label{eq:free-green-axis-integral}
  \begin{aligned}
    g_j&=\frac1{(2\pi)^3}\int_{[-\pi,\pi]^3}
    \frac{\cos(j\theta_3)}{\omega(\theta_1,\theta_2)+2-2\cos\theta_3}
    \,d\theta_1\,d\theta_2\,d\theta_3\\
    &=\frac1{(2\pi)^2}\int_{[-\pi,\pi]^2}
    \frac{q(\omega(\theta_1,\theta_2))^j}
    {\sqrt{\omega(\theta_1,\theta_2)(\omega(\theta_1,\theta_2)+4)}}
    \,d\theta_1\,d\theta_2.
  \end{aligned}
\end{equation}
By the Cauchy--Schwarz inequality applied to \eqref{eq:free-green-axis-integral}, we obtain
\[
  g_j^2\leq g_{j-1}g_{j+1},\qquad j\geq1.
\]
Hence,
\[
  \frac{g_j}{g_{j-1}}\geq\frac{g_1}{g_0},\qquad j\geq1.
\]
By \((\Delta_{\mathrm d}g_3)(0)=1\) and \(g_3(\pm e_i)=g_1\) for \(i=1,2,3\),
we have \(6g_0-6g_1=1\).
Together with \(g_0>2/9\) in \eqref{eq:scalar-green-bounds}, this gives
\[
  \frac{g_j}{g_{j-1}}\geq1-\frac1{6g_0}>\frac14,\qquad j\geq1.
\]
This completes the proof.
\end{proof}

\begin{proposition}
\label{prop:green-order}
With \(G_+\) and \(g_3\) as defined above,
\begin{equation}
  \label{eq:green-order-use}
  0<G_+((x,k),(y,l))<g_3(x-y,k-l)
  \qquad(x,y\in\mathbb Z^2,\ k,l\geq1).
\end{equation}
\end{proposition}
\begin{proof}
By Lemmas~\ref{lem:jacobi-reference} and~\ref{lem:heat-crossing}, it remains to prove
\[
  G_+((0,k),(0,l))<g_3(0,k-l),\qquad k,l\geq1.
\]
For \(k,l\geq2\), we first estimate
\(g_3(e,k-l)-G_+((e,k),(0,l))\), where \(e=(1,0)\).
Set \(d=|k-l|\) and \(s=k+l\).
The estimate \(M_e(\tau)>1/16\) will be used to prove
\eqref{eq:adjacent-unpinned}.
For this purpose, we first show that
\(\tau>1\), where \(\tau\)
is the unique positive solution of
\[
  (e^{-\tau B})_{kl}=p(\tau,|k-l|),
\]
as in the proof of Lemma~\ref{lem:heat-crossing}. Since \(G_+(z,w)=G_+(w,z)\) and \(g_3(z)=g_3(-z)\),
we may assume that \(k\leq l\).

By \eqref{eq:half-line-heat},
\((e^{-tD})_{kl}=p(t,d)-p(t,s)\).
Let \(A=2I-D\) and \(V_{ij}=\delta_{ij}/(4j^2)\). Since \(2I-L_1=A+V\),
\[
  \bigl((A+V)^{d+1}-A^{d+1}\bigr)_{kl}
  \geq\sum_{a=0}^d(A^aVA^{d-a})_{kl}
  =\sum_{j=k}^l\frac1{4j^2}.
\]
By \eqref{eq:operator-exponential},
\[
  (e^{-tL_1})_{kl}-p(t,d)
  \geq e^{-2t}\left\{\frac{t^{d+1}}{(d+1)!}
  \sum_{j=k}^l\frac1{4j^2}-I_s(2t)\right\}.
\]

Since \((x^{-2})''=6x^{-4}>0\),
\[
  \sum_{j=k}^l\frac1{j^2}\geq\int_k^l\frac{dx}{x^2}
  +\frac12\left(\frac1{k^2}+\frac1{l^2}\right)
  \geq\frac{d+1}{kl}.
\]
By \eqref{eq:bessel-series}, \(I_s(2t)\leq t^se^{t^2/(s+1)}/s!\).

For \(R_{k,d}:=4k(k+d)d!/(d+2k)!\),
\[
  \frac{R_{k,d+1}}{R_{k,d}}
  =\frac{(k+d+1)(d+1)}{(k+d)(d+2k+1)}<1,
  \qquad
  \frac{R_{k+1,0}}{R_{k,0}}=\frac{k+1}{2k^2(2k+1)}<1.
\]
Hence \(R_{k,d}\leq R_{2,0}=2/3\), and, for \(0<t\leq1\),
\[
  \frac{I_s(2t)}{t^{d+1}\bigl(\sum_{j=k}^l(4j^2)^{-1}\bigr)/(d+1)!}
  \leq R_{k,d}t^{2k-1}e^{t^2/(s+1)}
  \leq\frac23e^{1/5}\leq\frac56<1.
\]
By \eqref{eq:semigroup-comparison-explicit},
\((e^{-tB})_{kl}>(e^{-tL_1})_{kl}>p(t,d)\) for \(0<t\leq1\).
Thus \(\tau>1\).

For \(d>0\), \eqref{eq:jacobi-log-off-diagonal} gives
\[
  T_{kl}(0)=\frac{d}{S(S+2\sqrt{r_kr_l})}>\frac{d}{2s^2},
  \qquad S=r_k+r_l<s.
\]
It is strictly positive also when \(d=0\).

By \eqref{eq:MX-increasing}, \(M_e'(t)>0\). Since
\(D(t)<0\) for \(t<\tau\) and \(D(t)>0\) for \(t>\tau\),
\[
  \begin{aligned}
    g_3(e,k-l)-G_B((e,k),(0,l))-M_e(\tau)T_{kl}(0)
    &=\int_0^\infty(M_e(t)-M_e(\tau))\frac{D(t)}t\,dt\geq0.
  \end{aligned}
\]
By \eqref{eq:bessel-series}, \(I_0(2)>9/4\) and \(I_1(2)>19/12\).
Together with \(e^4<(68/25)^4<57\), these give
\(M_e(\tau)>M_e(1)=e^{-4}I_0(2)I_1(2)>1/16\).
It follows that
\begin{equation}
  \label{eq:adjacent-unpinned}
  g_3(e,k-l)-G_B((e,k),(0,l))\geq\frac{d}{32s^2}.
\end{equation}

Since the bound in \eqref{eq:adjacent-unpinned} is zero when \(d=0\),
we also estimate \(G_B-G_{B_*}\). For this purpose, define
\(\Pi_j(X)=g_3(X,j-1)-g_3(X,j+1)\).

By \eqref{eq:bessel-recurrences}, we have
\[
  \Pi_j(X)=j\int_0^\infty
  p(t,X_1)p(t,X_2)p(t,j)\,\frac{dt}{t}>0.
\]

By \eqref{eq:inverse-moment-transform},
\[
  \sum_{X\in\mathbb Z^2}\Pi_j(X)e^{-\mathrm iX\cdot\theta}=q(\omega)^j,
  \qquad \omega=4\sin^2(\theta_1/2)+4\sin^2(\theta_2/2).
\]
Moreover, by \(\sum_{a\in\mathbb Z}p(t,a)=1\) and \eqref{eq:inverse-moments},
\[
  \sum_{X\in\mathbb Z^2}\Pi_j(X)
  =j\int_0^\infty\frac{p(t,j)}t\,dt=1,
  \qquad
  \sum_{Y\in\mathbb Z^2}\Pi_i(X-Y)\Pi_j(Y)=\Pi_{i+j}(X).
\]

We first prove the following bound:
\begin{equation}
  \label{eq:poisson-lower}
  \Pi_j(e)\geq\frac1{2\pi}L(j),\qquad
  L(x)=\frac{x^2-25/4}{(x^2-1)(x^2-4)},\quad j\geq4.
\end{equation}

By \eqref{eq:free-heat-fourier}, integration by parts and
\(u=2\sqrt t\sin(\theta/2)\) give
\[
  \begin{aligned}
    p(t,1)&=\frac2{\pi\sqrt t}\int_0^{2\sqrt t}
    u^2e^{-u^2}\sqrt{1-u^2/(4t)}\,du\\
    &\geq\frac2{\pi\sqrt t}\int_0^\infty
    u^2e^{-u^2}\left(1-\frac{u^2}{4t}\right)\,du
    =(4\pi t)^{-1/2}\left(1-\frac3{8t}\right).
  \end{aligned}
\]
Here \(\sqrt{1-z}\geq1-z\) for \(0\leq z\leq1\), and
\(u^2e^{-u^2}(1-u^2/(4t))\leq0\) for \(u\geq2\sqrt t\).
For \(t\geq3/8\),
by \eqref{eq:scalar-green-bounds},
\(M_e(t)\geq c(1-3/(8t))\). For \(0<t<3/8\), we also have \(M_e(t)\geq c(1-3/(8t))\), since
\(M_e(t)>0>c(1-3/(8t))\).
By \eqref{eq:inverse-moments}, for \(j\geq4\),
\[
  \Pi_j(e)=j\int_0^\infty M_e(t)\frac{p(t,j)}{t^2}\,dt
  \geq\frac{j}{4\pi}\left(\frac{2}{j(j^2-1)}
  -\frac{3}{8}\frac{12}{j(j^2-1)(j^2-4)}\right)
  =\frac{L(j)}{2\pi}.
\]

For \(0\leq u<1\) and \(f\in C_c(\mathbb N)\),
\[
  \langle f,(D-uP_1)f\rangle
  =\sum_{j\geq1}|f_{j+1}-f_j|^2+(1-u)|f_1|^2\geq0.
\]
For \(a>0\), substitution into \((D-uP_1+aI)f=\delta_1\) gives
\[
  \bigl((D-uP_1+aI)^{-1}\bigr)_{k1}
  =\frac{q(a)^k}{1-uq(a)}=\sum_{r=0}^\infty u^rq(a)^{k+r}.
\]
Thus the Green function of \(\Delta_2+D-uP_1\)
satisfies
\[
  K_u(X;k):=\int_0^\infty p(t,X_1)p(t,X_2)
  (e^{-t(D-uP_1)})_{k1}\,dt
  =\sum_{r=0}^\infty u^r\Pi_{k+r}(X),\qquad
  \sum_X K_u(X;k)=\frac1{1-u}.
\]
By the definitions of \(B,B_*\),
\[
  0\leq(2I-D+vP_1)_{kl}\leq(2I-B)_{kl},\qquad
  0\leq(2I-D+P_1/4)_{kl}\leq(2I-B_*)_{kl}.
\]
By \eqref{eq:operator-exponential},
\begin{equation}
  \label{eq:pinned-heat-comparison}
  (e^{-t(D-vP_1)})_{kl}\leq(e^{-tB})_{kl},\qquad
  (e^{-t(D-P_1/4)})_{kl}\leq(e^{-tB_*})_{kl}.
\end{equation}

To estimate
\(G_B((X,k),(0,l))-G_{B_*}((X,k),(0,l))\),
we use \(B_*-B=\eta P_1\).
Define \(P_1f(X,k)=\mathbf1_{\{k=1\}}f(X,1)\) on \(\ell^2(D_+)\)
and consider the inverses
\[
  R_a=(\Delta_2+B+aI)^{-1},\qquad
  R_{*,a}=(\Delta_2+B_*+aI)^{-1},\quad a>0.
\]
Since \(B_*-B=\eta P_1\),
\[
  R_a-R_{*,a}=\eta R_aP_1R_{*,a}.
\]

For \(C\in\{B,B_*\}\),
\[
  e^{-at}p(t,X_1)p(t,X_2)(e^{-tC})_{kl}
  \uparrow p(t,X_1)p(t,X_2)(e^{-tC})_{kl}\geq0
  \quad\text{as }a\downarrow0.
\]
By \eqref{eq:pinned-heat-comparison} and Lebesgue's monotone convergence theorem,
letting \(a\downarrow0\), we obtain
\begin{align*}
  G_B((X,k),(0,l))-G_{B_*}((X,k),(0,l))
  &\geq\eta\sum_YK_v(X-Y;k)K_{1/4}(Y;l)\\
  &=\eta\sum_{a,b\geq0}v^a4^{-b}\Pi_{s+a+b}(X).
\end{align*}

By \(\sum_{X\in\mathbb Z^2}\Pi_j(X)=1\), we have
\[
  \sum_{X\in\mathbb Z^2}\sum_{a,b\geq0}v^a4^{-b}\Pi_{s+a+b}(X)
  =\frac4{3(1-v)}<\infty.
\]
Since \(v>1/2\) and \(\eta>29/108\), retaining the terms with \(a+b=0,1,2,3\) and using
\eqref{eq:poisson-lower} gives
\[
  G_B((e,k),(0,l))-G_{B_*}((e,k),(0,l))
  \geq\frac\eta{2\pi}\sum_{j=0}^3h_jL(s+j),
  \quad(h_0,h_1,h_2,h_3)=(1,3/4,7/16,15/64).
\]

For fixed \(j\geq0\)
and \(s\geq4\),
\[
  s^2L(s+j)=\left(\frac{s}{s+j}\right)^2(s+j)^2L(s+j),\qquad
  \frac{d}{ds}\left(\frac{s}{s+j}\right)^2=\frac{2sj}{(s+j)^3}\geq0.
\]
Also,
\[
  x^2L(x)=1-\frac{(5/4)x^2+4}{x^4-5x^2+4},
\]
and, for \(z=x^2\geq16\),
\[
  \frac{d}{dz}\left(\frac{5z/4+4}{z^2-5z+4}\right)
  =\frac{-5z^2/4-8z+25}{(z^2-5z+4)^2}<0.
\]
Hence \(s^2L(s+j)\geq16L(4+j)\), and
\[
  s^2\sum_{j=0}^3h_jL(s+j)
  >\frac{25}{16}.
\]
Since \(\eta/(2\pi)>203/4752\), this proves
\(G_B((e,k),(0,l))-G_{B_*}((e,k),(0,l))>1/(15s^2)\). Define
\(D_e(k,l)=g_3(e,k-l)-G_+((e,k),(0,l))\) for \(k,l\geq1\). Combining \(G_B((e,k),(0,l))-G_{B_*}((e,k),(0,l))>1/(15s^2)\) with
\eqref{eq:adjacent-unpinned} and \(G_+<G_{B_*}\), we obtain
\begin{equation}
  \label{eq:adjacent-full-gap}
  4D_e(k,l)>\frac{d/8+4/15}{s^2}.
\end{equation}

Set \(D_0(k,l)=g_{|k-l|}-G_+((0,k),(0,l))\) for \(k,l\geq1\),
and \(V_k=1/(4k^2)\).
The Green equations give
\begin{equation}
  \label{eq:D0-green-equation}
  (6-V_k)D_0(k,l)-D_0(k-1,l)-D_0(k+1,l)
  =4D_e(k,l)-V_kg_{|k-l|},\qquad k,l\geq2.
\end{equation}
To prove \(D_0(k,l)>0\), we first show that
\(4D_e(k,l)-V_kg_{|k-l|}>0\).
First, \(g_d<1/(10d)\) for \(d\geq1\). When \(d=1\),
this is \(g_1=g_0-1/6<1/10\), whereas for \(d\geq2\),
\eqref{eq:scalar-green-bounds} and \eqref{eq:inverse-moments} give
\[
  g_d<\frac{1}{4\pi d}\left(1+\frac1{2(d^2-1)}\right)
  \leq\frac7{24\pi d}<\frac1{10d}.
\]

For \(d=0\), we have \(s^2V_kg_0=g_0<4/15\).
For \(d\geq1\), since \(k\geq2\) and \(l\leq k+d\),
\[
  s^2V_kg_d<\frac{(1+d/4)^2}{10d}
  <\frac d8+\frac4{15}.
\]
Hence, by \eqref{eq:adjacent-full-gap}, \(4D_e(k,l)>V_kg_d\) in every case.

Fix \(l\geq2\). To prove \(D_0(k,l)>0\), we first show that
\(D_0(\cdot,l)\in\ell^\infty(\mathbb N)\).
By \(B\geq0\) and \eqref{eq:jacobi-heat-integral}, we have
\[
  0\leq(e^{-tB})_{kk}\leq1,\qquad
  0\leq(e^{-tB})_{ll}\leq C_l(1+t)^{-1},\qquad k\geq1,\ t>0.
\]
By the Cauchy--Schwarz inequality, we obtain
\[
  |(e^{-tB})_{kl}|^2\leq(e^{-tB})_{kk}(e^{-tB})_{ll}
  \leq C_l(1+t)^{-1}.
\]
By \eqref{eq:free-heat-fourier} and Lemma~\ref{lem:jacobi-reference}, we obtain
\[
  \begin{aligned}
    0<G_+((X,k),(0,l))&\leq G_B((X,k),(0,l))\\
    &\leq C_l\int_0^\infty(1+t)^{-3/2}\,dt<\infty,
    \qquad X\in\mathbb Z^2,\ k\geq1.
  \end{aligned}
\]
Since \(0<g_3(z)\leq g_3(0)\) for \(z\in\mathbb Z^3\), we have
\[
  \sup_{k\geq1}|D_0(k,l)|\leq g_3(0)+C_l<\infty.
\]

By \eqref{eq:D0-green-equation} and \(4D_e(k,l)>V_kg_{|k-l|}\), we have
\[
  (6-V_k)D_0(k,l)-D_0(k-1,l)-D_0(k+1,l)>0,
  \qquad k\geq2.
\]
By Lemmas~\ref{lem:jacobi-reference} and~\ref{lem:heat-crossing}, we have \(D_e(1,l)>0\).
Hence, by Lemma~\ref{lem:free-green-ratio} and the equations for \(g_3\) and \(G_+\), we obtain
\begin{equation}
  \label{eq:D0-first-layer-equation}
  \frac{23}{4}D_0(1,l)-D_0(2,l)
  =4D_e(1,l)+g_l-\frac14g_{l-1}>0.
\end{equation}

To prove  \(D_0(k,l)\geq0\) for \(k\geq1\), set
\(w(k)=\max\{-D_0(k,l),0\}\).
By \eqref{eq:D0-green-equation}, \eqref{eq:D0-first-layer-equation}, and
\(4D_e(k,l)>V_kg_{|k-l|}\) for \(k\geq2\), we have
\[
  w(1)\leq\frac4{23}w(2),\qquad
  w(k)\leq\frac{16}{95}\bigl(w(k-1)+w(k+1)\bigr),\quad k\geq2.
\]
Since \(w\in\ell^\infty(\mathbb N)\), taking the supremum over \(k\geq1\) gives
\[
  \|w\|_{\ell^\infty}\leq\max\left\{\frac4{23},\frac{32}{95}\right\}
  \|w\|_{\ell^\infty}=\frac{32}{95}\|w\|_{\ell^\infty}.
\]
Hence \(w=0\), so \(D_0(k,l)\geq0\) for every \(k\geq1\).
By \eqref{eq:D0-green-equation} and \eqref{eq:D0-first-layer-equation}, we further obtain
\[
  D_0(k,l)>0,\qquad k\geq1,\quad l\geq2.
\]

For \(k\geq2\), symmetry gives
\(D_0(k,1)=D_0(1,k)>0\).
Finally, by Lemmas~\ref{lem:jacobi-reference}--\ref{lem:free-green-ratio}
and \(D_0(2,1)>0\), the Green equations give
\[
  \frac{23}{4}D_0(1,1)
  =D_0(2,1)+4D_e(1,1)+g_1-\frac14g_0>0.
\]
Thus \(D_0(k,l)>0\) for all \(k,l\geq1\).
Together with Lemmas~\ref{lem:jacobi-reference} and~\ref{lem:heat-crossing},
this proves \eqref{eq:green-order-use}.
This completes the proof.
\end{proof}

\subsection[The duality argument for (m,n)=(1,2)]{The duality argument for \((m,n)=(1,2)\)}
To prove the Hardy--Sobolev--Maz'ya inequality on \(D_+\) using
\(G_+(z,w)<g_3(z-w)\) in Proposition~\ref{prop:green-order},
we introduce the following operators.
For \(u\in C_c(D_+)\), extended by zero to \(\mathbb Z^3\), define
\[
  Q_+(u)=\sum_{z\in\mathbb Z^3}|Du(z)|^2
  -\frac14\sum_{x\in\mathbb Z^2}\sum_{k\geq1}\frac{|u(x,k)|^2}{k^2}.
\]
Let \(\mathcal H_+\) be the completion of \(C_c(D_+)\) under the norm
\(Q_+^{1/2}\). We also denote by \(Q_+(u,v)\) the continuous bilinear form
on \(\mathcal H_+\).

By \eqref{eq:HSM-sufficiency-bound} with \(m=1\) and \(n=2\),
\begin{equation}
  \label{eq:halfspace-HSM-bound}
  \|u\|_{\ell^6(D_+)}^2\leq C Q_+(u),\qquad u\in\mathcal H_+.
\end{equation}
Thus \(\mathcal H_+\) embeds continuously into \(\ell^6(D_+)\).
For \(f\in\ell^{6/5}(D_+)\), H\"older's inequality gives
\[
  \left|\sum_{z\in D_+}f(z)u(z)\right|
  \leq C\|f\|_{\ell^{6/5}(D_+)}Q_+(u)^{1/2},
  \qquad u\in\mathcal H_+.
\]
The Riesz representation theorem therefore gives a unique
\(R_+f\in\mathcal H_+\) satisfying
\begin{equation}
  \label{eq:Rplus-definition}
  Q_+(R_+f,u)=\sum_{z\in D_+}f(z)u(z),\qquad u\in\mathcal H_+.
\end{equation}

To compare \(R_+f\) with the solution on \(\mathbb Z^3\), we extend
\(f\) by zero to \(\mathbb Z^3\). The Sobolev inequality and the Riesz
representation theorem give a unique \(R_0f\in D^{1,2}(\mathbb Z^3)\) satisfying
\begin{equation}
  \label{eq:R0-definition}
  \sum_{j=1}^3\sum_{z\in\mathbb Z^3}
  D_j(R_0f)(z)D_jv(z)=\sum_{z\in\mathbb Z^3}f(z)v(z),
  \qquad v\in D^{1,2}(\mathbb Z^3).
\end{equation}
By H\"older's inequality and the Sobolev inequality on \(D^{1,2}(\mathbb Z^3)\),
we obtain
\[
  \begin{aligned}
    \sum_{z\in\mathbb Z^3}|D(R_0f)(z)|^2
    &=\langle f,R_0f\rangle\\
    &\leq\|f\|_{6/5}\|R_0f\|_6\\
    &\leq\overline S_3^{-1/2}\|f\|_{6/5}
    \left(\sum_{z\in\mathbb Z^3}|D(R_0f)(z)|^2\right)^{1/2}.
  \end{aligned}
\]
Hence,
\begin{equation}
  \label{eq:R0-dual-bound}
  \langle f,R_0f\rangle\leq\overline S_3^{-1}\|f\|_{\ell^{6/5}(D_+)}^2.
\end{equation}
By \eqref{eq:Rplus-definition}, \eqref{eq:R0-definition}, H\"older's inequality,
and the corresponding Sobolev inequalities, we obtain
\begin{equation}
  \label{eq:resolvents-bounded}
  \|R_+f\|_{\ell^6(D_+)}+\|R_0f\|_{\ell^6(\mathbb Z^3)}
  \leq C\|f\|_{\ell^{6/5}(D_+)}.
\end{equation}

To use the comparison in Proposition~\ref{prop:green-order}, we first establish the following representations
for every \(f\in C_c(D_+)\):
\begin{equation}
  \label{eq:resolvent-green-representations}
  \begin{aligned}
    R_+f(z)&=\sum_{w\in D_+}G_+(z,w)f(w),&&z\in D_+,\\
    R_0f(z)&=\sum_{w\in D_+}g_3(z-w)f(w),&&z\in\mathbb Z^3.
  \end{aligned}
\end{equation}
Choose increasing finite sets \(D_j\)
such that
\(\operatorname{supp}f\subset D_1\) and
\(\bigcup_{j\geq1}D_j=D_+\).
Let \(u_j\) solve \(H_+u_j=f\) in \(D_j\), with
\(u_j=0\) on \(D_+\setminus D_j\).
By the equations for \(u_j\) and \(R_+f\), we have
\[
  Q_+(R_+f-u_j,\varphi)=0,
  \qquad \varphi\in C_c(D_j).
\]

Since \(\bigcup_j C_c(D_j)=C_c(D_+)\) is dense in \(\mathcal H_+\), we obtain
\[
  Q_+(u_j-R_+f)
  =\inf_{\varphi\in C_c(D_j)}Q_+(\varphi-R_+f)
  \xrightarrow[j\rightarrow\infty]{}0.
\]

By \eqref{eq:halfspace-HSM-bound}, \(u_j(z)\to R_+f(z)\) for each \(z\in D_+\).

For each \(w\in D_j\), the Green function \(G_j(\cdot,w)\) satisfies
\[
  \bigl(H_+G_j(\cdot,w)\bigr)(z)=\delta_w(z)\quad(z\in D_j),
  \qquad G_j(z,w)=0\quad(z\in D_+\setminus D_j).
\]
Since \(G_j(z,w)\uparrow G_+(z,w)\) as \(j\rightarrow\infty\)
and \(f\) has finite support, we have
\[
  u_j(z)
  =\sum_{w\in\operatorname{supp}f}G_j(z,w)f(w)
  \xrightarrow[j\rightarrow\infty]{}
  \sum_{w\in D_+}G_+(z,w)f(w).
\]

This proves the formula for \(R_+f\) in \eqref{eq:resolvent-green-representations}.

To obtain the representation of \(R_0f\) in \eqref{eq:resolvent-green-representations},
recall that \(f(z)=0\) for \(z\in\mathbb Z^3\setminus D_+\), and define
\[
  v(z):=\sum_{w\in\operatorname{supp}f}g_3(z-w)f(w),
  \qquad z\in\mathbb Z^3.
\]
By \(\Delta_{\mathrm d}g_3=\delta_0\) and \eqref{eq:R0-definition}
with test function \(\delta_z\), we have
\[
  \Delta_{\mathrm d}v(z)=f(z)=\Delta_{\mathrm d}(R_0f)(z),
  \qquad z\in\mathbb Z^3.
\]
Since \(f\) has finite support, the estimate \(g_3(z)\leq C(1+|z|)^{-1}\)
in the proof of Lemma~\ref{lem:free-summability} gives \(v(z)\to0\) as \(|z|\to\infty\).
By \eqref{eq:resolvents-bounded}, \(R_0f\in\ell^6(\mathbb Z^3)\), so
\(v-R_0f\) also tends to zero at infinity.
Applying the maximum principle on finite balls to the harmonic function
\(v-R_0f\), we obtain \(v=R_0f\) on \(\mathbb Z^3\).

Combining Proposition~\ref{prop:green-order} with
\eqref{eq:resolvent-green-representations} and \eqref{eq:R0-dual-bound},
we obtain the following strict estimate.
\begin{lemma}
\label{lem:strict-dual-three}
For every \(0\ne f\in\ell^{6/5}(D_+)\),
\[
  \langle f,R_+f\rangle
  <\overline S_3^{-1}\|f\|_{\ell^{6/5}(D_+)}^2.
\]
\end{lemma}
\begin{proof}
Let \(0\ne f\in\ell^{6/5}(D_+)\), and choose \(z_0\in D_+\) such that
\(f(z_0)\ne0\). By Proposition~\ref{prop:green-order},
\[
  \varepsilon_f:=\bigl(g_3(0)-G_+(z_0,z_0)\bigr)|f(z_0)|^2>0.
\]
Choose finite sets \(K_j\subset K_{j+1}\) such that
\(z_0\in K_1\) and \(\bigcup_{j\geq1}K_j=D_+\), and set
\(f_j=f\mathbf1_{K_j}\).

By \eqref{eq:resolvent-green-representations} and \(G_+(z,w)>0\), we have
\[
  \langle f_j,R_+f_j\rangle\leq\langle |f_j|,R_+|f_j|\rangle.
\]
Moreover, Proposition~\ref{prop:green-order} gives
\[
  \begin{aligned}
    \langle |f_j|,R_0|f_j|\rangle-\langle |f_j|,R_+|f_j|\rangle
    &=\sum_{z,w\in K_j}\bigl(g_3(z-w)-G_+(z,w)\bigr)|f(z)||f(w)|\\
    &\geq\varepsilon_f.
  \end{aligned}
\]
Combining these inequalities with \eqref{eq:R0-dual-bound}, we obtain
\[
  \langle f_j,R_+f_j\rangle
  \leq\overline S_3^{-1}\|f_j\|_{\ell^{6/5}(D_+)}^2-\varepsilon_f.
\]
Since \(f_j\to f\) in \(\ell^{6/5}(D_+)\), by
\eqref{eq:resolvents-bounded} and H\"older's inequality,
letting \(j\to\infty\) gives
\begin{equation}
  \langle f,R_+f\rangle
  \leq\overline S_3^{-1}\|f\|_{\ell^{6/5}(D_+)}^2-\varepsilon_f
  <\overline S_3^{-1}\|f\|_{\ell^{6/5}(D_+)}^2.
\end{equation}
This completes the proof.
\end{proof}

We next prove that
\begin{equation}
  \label{eq:halfspace-strict-HSM}
  Q_+(u)>\overline S_3\|u\|_{\ell^6(D_+)}^2,
  \qquad 0\ne u\in\mathcal H_+.
\end{equation}
For \(0\ne u\in\mathcal H_+\), set
\(f=|u|^4u/\|u\|_{\ell^6(D_+)}^5\).
Then \(\|f\|_{\ell^{6/5}(D_+)}=1\) and
\(\langle f,u\rangle=\|u\|_{\ell^6(D_+)}\).
By \eqref{eq:Rplus-definition}, the Cauchy--Schwarz inequality,
and Lemma~\ref{lem:strict-dual-three}, we obtain
\[
  \|u\|_{\ell^6(D_+)}^2=|Q_+(R_+f,u)|^2
  \leq\langle f,R_+f\rangle Q_+(u)<\overline S_3^{-1}Q_+(u).
\]
This proves \eqref{eq:halfspace-strict-HSM}.
To apply \eqref{eq:halfspace-strict-HSM} to functions on
\(\mathbb Z^3\setminus\Sigma\), set \(D_-=\mathbb Z^2\times(-\mathbb N)\).
For \(v\in C_c(D_-)\), define
\begin{equation}
  \label{eq:negative-halfspace-reflection}
  \widetilde v(x,k)=v(x,-k),\qquad Q_-(v)=Q_+(\widetilde v).
\end{equation}
For \(u\in C_c(\mathbb Z^3\setminus\Sigma)\), set
\begin{equation}
  \label{eq:halfspace-components}
  u_\pm=u\mathbf1_{D_\pm}.
\end{equation}
Since no edge joins \(D_+\) to \(D_-\), we have
\[
  \sum_{z\in\mathbb Z^3}Du_+(z)\cdot Du_-(z)
  -\frac14\sum_{z\in\mathbb Z^3\setminus\Sigma}
  \frac{u_+(z)u_-(z)}{d_\Sigma(z)^2}=0.
\]
By the definitions of \(u_+\) and \(u_-\), we have
\begin{equation}
  \label{eq:halfspace-decomposition}
  \begin{aligned}
    u&=u_++u_-,\qquad
    \|u\|_{\mathcal H_{1,2}}^2=Q_+(u_+)+Q_-(u_-),\\
    \|u\|_{\ell^6(\mathbb Z^3)}^6
    &=\|u_+\|_{\ell^6(D_+)}^6+\|u_-\|_{\ell^6(D_-)}^6.
  \end{aligned}
\end{equation}
Since \(Q_\pm(u_\pm)\leq\|u\|_{\mathcal H_{1,2}}^2\), the maps
\(u\mapsto u_\pm\) extend continuously to the completions of
\(C_c(D_\pm)\) under \(Q_\pm^{1/2}\).
Consequently, \eqref{eq:halfspace-decomposition} holds for every
\(u\in\mathcal H_{1,2}\).

For every nonzero \(v\) in the completion of \(C_c(D_-)\) under \(Q_-^{1/2}\), set
\[
  \widetilde v(x,k)=v(x,-k),\qquad (x,k)\in D_+.
\]
Then \(\widetilde v\in\mathcal H_+\). By \eqref{eq:negative-halfspace-reflection} and \eqref{eq:halfspace-strict-HSM}, we obtain
\begin{equation}
  \label{eq:negative-halfspace-strict}
  Q_-(v)=Q_+(\widetilde v)
  >\overline S_3\|\widetilde v\|_{\ell^6(D_+)}^2
  =\overline S_3\|v\|_{\ell^6(D_-)}^2.
\end{equation}
Let \(0\ne u\in\mathcal H_{1,2}\), and let \(u_\pm\) be defined by
\eqref{eq:halfspace-components}. Since \(u=u_++u_-\ne0\), at least one of
\(u_+\) and \(u_-\) is nonzero.
By \eqref{eq:halfspace-strict-HSM}, \eqref{eq:halfspace-decomposition}, and
\eqref{eq:negative-halfspace-strict}, we obtain
\begin{equation}
  \label{eq:fullspace-strict-HSM}
  \begin{aligned}
    \|u\|_{\mathcal H_{1,2}}^2
    &=Q_+(u_+)+Q_-(u_-)\\
    &>\overline S_3\bigl(\|u_+\|_6^2+\|u_-\|_6^2\bigr)\\
    &\geq\overline S_3\bigl(\|u_+\|_6^6+\|u_-\|_6^6\bigr)^{1/3}
    =\overline S_3\|u\|_6^2.
  \end{aligned}
\end{equation}
By \eqref{eq:fullspace-strict-HSM} and Lemma~\ref{lem:general-sobolev-upper},
we obtain
\begin{equation}
  \label{eq:kappa-three-equality}
  \overline\kappa_{1,2}=\overline S_3.
\end{equation}

We are now in a position to prove Theorem~\ref{thm:sobolev-comparison}.\par\nopagebreak
\begin{proof}[Proof of Theorem~\ref{thm:sobolev-comparison}]
Proposition~\ref{prop:m2} proves the assertion for \(m=2\),
and \eqref{eq:kappa-three-equality} proves the assertion for \((m,n)=(1,2)\). If \(m\neq2\), the Hardy constant
\(c_m\) is positive: for \(m=1\) it is \(1/4\), and for \(m\geq3\) this was
proved in Section~\ref{sec:positivity-criterion}. Lemma~\ref{lem:strict-sobolev} applies
when \(N\geq4\). Nonnegativity follows from the Hardy inequality.
This completes the proof.
\end{proof}

\section[Nonattainment of kappa(2,n) and kappa(1,2)]{Nonattainment of \(\overline\kappa_{2,n}\) and \(\overline\kappa_{1,2}\)}
\label{sec:attainment}

In this section, we prove Corollary~\ref{cor:nonattainment-m2} and
Theorem~\ref{thm:attainment}.
Since \(\overline\kappa_{2,n}=\overline S_{n+2}\), the proof of
Corollary~\ref{cor:nonattainment-m2} reduces to showing that no Sobolev
extremal function on \(\mathbb Z^{n+2}\) satisfies \(u|_\Sigma=0\).
For this purpose, we use the Euler--Lagrange equation
\eqref{eq:sobolev-extremal-equation}.
Theorem~\ref{thm:attainment} follows from \eqref{eq:fullspace-strict-HSM}
and \(\overline\kappa_{1,2}=\overline S_3\) in Theorem~\ref{thm:sobolev-comparison}.

\subsection[The case m=2]{The case \(m=2\)}
Throughout this subsection, \(N=n+2\geq3\),
and \(\Sigma=\mathbb Z^n\times\{0\}\).

Since \(c_2=0\), we have
\(\|u\|_{\mathcal H_{2,n}}=\|u\|_{D^{1,2}(\mathbb Z^N)}\)
for every \(u\in C_c(\mathbb Z^N\setminus\Sigma)\).
Hence \(\mathcal H_{2,n}\) is the closure of \(C_c(\mathbb Z^N\setminus\Sigma)\)
in \(D^{1,2}(\mathbb Z^N)\).
By Lemma~\ref{lem:Hmn-embedding}, every \(u\in\mathcal H_{2,n}\)
satisfies \(u|_\Sigma=0\).

We are now in a position to prove Corollary~\ref{cor:nonattainment-m2}.\par\nopagebreak

\begin{proof}[ Proof of Corollary~\ref{cor:nonattainment-m2}]
We argue by contradiction. Suppose that a nonzero
\(u\in\mathcal H_{2,n}\) attains \(\overline\kappa_{2,n}\).

Then \(u\in D^{1,2}(\mathbb Z^N)\), and
\[
  \frac{\sum_{z\in\mathbb Z^N}|Du(z)|^2}{\|u\|_{2^*}^2}
  =\overline\kappa_{2,n}=\overline S_N.
\]
Since \(|u|\in\ell^{2^*}\) and \(\sum_{z\in\mathbb Z^N}|D(|u|)(z)|^2\leq \sum_{z\in\mathbb
Z^N}|Du(z)|^2<\infty\),
\eqref{eq:full-space-identification} gives
\(|u|\in D^{1,2}(\mathbb Z^N)\). The definition of
\(\overline S_N\) then shows that \(|u|\) is also an extremal
function for the Sobolev inequality on \(\mathbb Z^N\).

Replacing \(u\) by \(|u|/\|u\|_{2^*}\), we may assume that \(u\geq0\),
\(\|u\|_{2^*}=1\), and \(u|_\Sigma=0\).
The Euler--Lagrange equation on \(\mathbb Z^N\) is
\begin{equation}
  \label{eq:sobolev-extremal-equation}
  \Delta_{\mathrm d} u(z)=\overline S_Nu(z)^{2^*-1},
  \qquad z\in\mathbb Z^N.
\end{equation}
Since \(u\) attains \(\overline S_N\) in \(D^{1,2}(\mathbb Z^N)\),
\eqref{eq:sobolev-extremal-equation} also holds for \(z\in\Sigma\).

Choose \(z_0\in\Sigma\). Since \(u(z_0)=0\),
by \eqref{eq:sobolev-extremal-equation}, we have
\[
  0=\Delta_{\mathrm d} u(z_0)=-\sum_{z\sim z_0}u(z).
\]

Since \(u\geq0\), we have
\[
  u(z)=0\ \Longrightarrow\ \sum_{w\sim z}u(w)=0
  \ \Longrightarrow\ u(w)=0\quad(w\sim z).
\]
Starting from \(z_0\) and applying this at each neighboring vertex, we obtain
\(u\equiv0\) on \(\mathbb Z^N\). This contradicts
\(\|u\|_{2^*}=1\).
This completes the proof.
\end{proof}

\subsection[The case (m,n)=(1,2)]{The case \((m,n)=(1,2)\)}

We are now in a position to prove Theorem~\ref{thm:attainment}.\par\nopagebreak
\begin{proof}[Proof of Theorem~\ref{thm:attainment}]
Since \(\overline\kappa_{1,2}=\overline S_3\) by
Theorem~\ref{thm:sobolev-comparison}, \eqref{eq:fullspace-strict-HSM} shows that
\(\overline\kappa_{1,2}\) is not attained in \(\mathcal H_{1,2}\).
This completes the proof.
\end{proof}

\appendix

\section[Estimates for c9]{Estimates for \(c_9\)}
\label{sec:nine-estimates}

To complete the proof of Theorem~\ref{thm:main} for \(m=9\),
we need to compare \(c_9\) with
\(49/4-(n/(n+9))^2\) in
\eqref{eq:positivity-criterion}.
For this purpose, we prove the following bounds.

\begin{proposition}
\label{prop:nine-bounds}
The Hardy constant \(c_9\) satisfies
\[
  \frac{25137}{2116}<\frac{297}{25}\leq c_9<\frac{1440}{121}.
\]
\end{proposition}

We prove the lower bound in Subsection~\ref{subsec:nine-lower-bound}
and the upper bound in Subsection~\ref{subsec:nine-upper-bound}.
Throughout this appendix,
sums over \(x\) are over \(\mathbb Z^9\),
and sums over \(x\neq0\) are over \(\mathbb Z^9\setminus\{0\}\).
We use test functions \(u\in C_c(\mathbb Z^9)\) with \(u(0)=0\).

\subsection{The lower bound}\label{subsec:nine-lower-bound}

\begin{proof}[Proof of the lower bound in Proposition~\ref{prop:nine-bounds}]
Set \(\lambda=297/25\). For every oriented edge with nonzero
endpoints, choose \(a_{xy}>0\) with \(a_{xy}a_{yx}=1\).
The identity
\[
  |u(x)-u(y)|^2-(1-a_{xy})|u(x)|^2-(1-a_{yx})|u(y)|^2
  =\left|\sqrt{a_{xy}}u(x)-a_{xy}^{-1/2}u(y)\right|^2
\]
gives, for \(u\in C_c(\mathbb Z^9)\) with \(u(0)=0\),
\begin{equation}
  \label{eq:nine-edge-weight}
  \sum_x |Du(x)|^2\geq\sum_{x\neq0}W(x)|u(x)|^2,
  \qquad W(x)=18-\sum_{y\sim x,\ y\neq0}a_{xy}.
\end{equation}
It suffices to choose \(a_{xy}\) so that \(W(x)\geq\lambda/|x|^2\) for
\(x\neq0\). We first prove \(W(x)\geq\lambda/|x|^2\) for \(|x|^2\geq9\),
and then consider \(0<|x|^2<9\).
\subsubsection[Construction of axy for max radii squared >= 9]{Construction of \(a_{xy}\) for
\(\max\{|x|^2,|y|^2\}\geq9\)}
To obtain $a_{xy}a_{yx}=1$, we
define
$a_{xy}=(1-F_{xy})/(1+F_{xy})$ with $F_{yx}=-F_{xy}$.

To define \(F_{xy}\) in \eqref{nine:quadratic-field},
we first construct a polynomial \(Y\) satisfying \(\Delta_{\mathrm d}Y=0\). Set $s=|x|^2$ and $P=\sum_i
x_i^4$. The identities
$-\Delta_{\mathrm d}P=12s+18$, $-\Delta_{\mathrm d}s^2=44s+18$, and
$-\Delta_{\mathrm d}s=18$ show that the polynomial
\[
  Y(x)=P-\frac{3s^2+8s}{11}
\]
satisfies $\Delta_{\mathrm d}Y=0$. For adjacent vertices, set
$S=|x|^2+|y|^2$ and $\delta=|y|^2-|x|^2$. Since
$P(y)-P(x)=(\delta^3+\delta)/2$, one has
\begin{equation*}
  Y(y)-Y(x)=\frac{\delta}{22}(11\delta^2-6S-5).
\end{equation*}

Let $C=\{x\in\mathbb Z^9:0<|x|^2<9\}$. On every edge with nonzero endpoints not both in $C$, define
\begin{equation}
  \label{nine:quadratic-field}
  {\quad
  F_{xy}=\frac{\displaystyle\delta\left(\frac{115}{56}S-\frac52
  +\frac{2(11\delta^2-6S-5)}{3S}\right)}
  {\displaystyle S^2-\frac{185}{42}S-\frac49+\frac{68}{9}\delta^2},
  \qquad a_{xy}=\frac{1-F_{xy}}{1+F_{xy}}.\quad}
\end{equation}

Here \(2\delta(11\delta^2-6S-5)/(3S)=44[Y(y)-Y(x)]/(3S)\).

We first prove that $|F_{xy}|<1$.
Since $\max(|x|,|y|)\ge3$ and
$\big||x|-|y|\big|\le1$, we have $S\ge13$.
If $y=x\pm e_i$, then
\[
  \delta^2=(2x_i\pm1)^2\leq4|x|^2\pm4x_i+1=2S-1.
\]
Set
\[
  A=\frac{115}{56}S-\frac52+\frac{2(11\delta^2-6S-5)}{3S},\qquad
  D=S^2-\frac{185}{42}S-\frac49+\frac{68}{9}\delta^2.
\]
By \(S\geq13\), \(\delta^2\leq2S-1\), and \(-185S/42-4/9>-9S/2\), we obtain
\[
  0<2S-7<A<\frac{21}{10}S+9<3S,\qquad
  D>\frac{S^2+9\delta^2}{2}\ge3S|\delta|.
\]

Here \((S^2+9\delta^2)/2-3S|\delta|=(S-3|\delta|)^2/2\geq0\). Hence $D>|\delta|A$ and $|F_{xy}|<1$.
Since $F_{yx}=-F_{xy}$,
we have
$a_{xy}>0$ and $a_{xy}a_{yx}=1$.

\subsubsection{Reduction to one variable}
We now prove that
\begin{equation}
  \label{nine:quadratic-exterior}
  W(x)\ge\frac{297}{25|x|^2},\qquad |x|^2\ge9.
\end{equation}

For this purpose, we will prove \eqref{nine:harmonic-middle-average}, namely,
\[
  sW(x)=\sum_{i:t_i>0}t_iQ_s(t_i),\qquad
  t_i=x_i^2/s,\quad\sum_{i=1}^9t_i=1,
\]
where \(Q_s\) is defined below.
For $s\ge9$, substitute $\delta=1+2X$ and $S=2s+1+2X$
into \eqref{nine:quadratic-field} and express $F/(1+F)$ as
$N_s(X)/H_s(X)$, where $N_s(X)=\sum_{j=0}^3n_jX^j$ and
$H_s(X)=\sum_{j=0}^3h_jX^j$
have the following coefficients:
\[
  \begin{array}{c|l|l}
    j&n_j&h_j\\\hline
    0&4140s^2-2412s-225&4032s^3+1308s^2-1100s+1643\\
    1&8280s^2+3456s+11922&20376s^2+28256s+28466\\
    2&16560s+43668&59120s+86532\\
    3&37848&72344
  \end{array}
\]
For $t=X^2/s\in[0,1]$, define
\begin{align*}
  B&=(h_0+h_2st)^2-st(h_1+h_3st)^2,\\
  J&=(n_0+n_2st)(h_0+h_2st)-st(n_1+n_3st)(h_1+h_3st),\\
  \phi_s(t)&=\frac{4J}{B},\qquad T=900sn_0-297h_0,\\
  \mathcal P(s,t)&=\frac13\left[100s\frac{h_0J-n_0B}{t}+TB\right]
  =\sum_{j=0}^3C_j(s)t^j.
\end{align*}

Since \((h_0J-n_0B)|_{t=0}=0\), $\mathcal P(s,t)$ is a
polynomial also at $t=0$. For $s\ge9$ and $|X|\le\sqrt s$, the bounds
$S\ge13$ and $\delta^2\le2S-1$ give $B=H_s(X)H_s(-X)>0$ and $h_0>0$.

By the definitions of \(F_{xy}\), \(N_s\), and \(H_s\), we have
\[
  \phi_s(x_i^2/s)
  =\frac{2F_{x,x+e_i}}{1+F_{x,x+e_i}}
  +\frac{2F_{x,x-e_i}}{1+F_{x,x-e_i}}.
\]
Set
\[
  Q_s(t)=s\left(9\phi_s(0)+\frac{\phi_s(t)-\phi_s(0)}t\right),
  \qquad 0<t\le1.
\]
By the definitions of \(\phi_s,T\), and \(\mathcal P\),
\begin{equation}
  \label{nine:harmonic-scalar-identity}
  \begin{aligned}
    Q_s(t)-\lambda
    &=\frac{100s(h_0J-n_0B)+tTB}{25th_0B}\\
    &=\frac{3\mathcal P(s,t)}{25h_0B}.
  \end{aligned}
\end{equation}
If $t_i=x_i^2/s$, then $\sum_i t_i=1$ and
\begin{equation}
  \label{nine:harmonic-middle-average}
  sW(x)=\sum_{i:t_i>0}t_iQ_s(t_i).
\end{equation}
Indeed, the definition of $Q_s$ and $\sum_i t_i=1$ give
\[
  \sum_{i:t_i>0}t_iQ_s(t_i)
  =s\left(9\phi_s(0)+\sum_{i:t_i>0}
  \bigl(\phi_s(t_i)-\phi_s(0)\bigr)\right)
  =s\sum_{i=1}^9\phi_s(t_i)=sW(x).
\]
Since $h_0B>0$, \eqref{nine:harmonic-scalar-identity} shows
that it suffices to prove $\mathcal P(s,t)\geq0$ for the
values of $s,t$ under consideration.

\subsubsection[The case s >= 17]{The case $s\ge17$}
For $s\ge17$, we express $\mathcal P(s,t)$
as a sum involving
$(4t-1)^2$, $1-t$, and $t$,
then prove that its coefficients are positive.
By the definitions of $n_j,h_j,B,J,T$, $C_j/4^j$ is an integer polynomial.
Set $q_j=C_j/(4^j10^8)$ and
\[
  \begin{aligned}
    U_2(s,t)&=q_2+(2+4t)q_3,\qquad U_0(s)=q_0-q_2-2q_3,\\
    U_1(s)&=q_0+4q_1+7q_2+10q_3.
  \end{aligned}
\]
Then
\begin{equation}
  \label{nine:harmonic-positive-cubic}
  \begin{aligned}
    \frac{\mathcal P(s,t)}{10^8}
    &=q_0+4q_1t+16q_2t^2+64q_3t^3\\
    &=U_2(s,t)(4t-1)^2+(1-t)U_0(s)+tU_1(s).
  \end{aligned}
\end{equation}
To prove $U_0(s),U_1(s),U_2(s,t)>0$, we first compare $U_2(s,t)$ with $U_2(s,1)$. By the definitions of $n_0,h_0$, and $T$, for $s\ge17$ we have
\[
 n_0>3000s^2,\qquad h_0<5000s^3,\qquad
 T=900sn_0-297h_0>0.
\]
By the definitions of $\mathcal P$ and $q_3$, we obtain
\[
 C_3=-\frac{s^3h_3^2T}{3}<0,\qquad q_3<0.
\]
Hence
\[
 U_2(s,t)=q_2+(2+4t)q_3\ge q_2+6q_3=U_2(s,1),\qquad 0\le t\le1.
\]

To prove $U_0(s),U_1(s),U_2(s,1)>0$, set $z=s/17-1$ and
\[
 V_i(z)=17^{-9}U_i(17(1+z))\quad(i=0,1),\qquad
 V_2(z)=17^{-9}U_2(17(1+z),1).
\]
Let
\[
 V_i(z)=\sum_{j=0}^{d_i}b_{i,j}z^j,\qquad (d_0,d_1,d_2)=(9,9,7).
\]
By substituting $s=17(1+z)$ into the definitions of $q_j,U_0,U_1$, and $U_2$, we obtain
\[
 b_{0,j}>3\ (0\le j\le9),\qquad
 b_{1,j}>7\ (0\le j\le9),\qquad
 b_{2,j}>64\ (0\le j\le7).
\]
Since $z\ge0$, we have $V_0(z),V_1(z),V_2(z)>0$. Therefore, $U_0(s),U_1(s),U_2(s,t)>0$ for $s\ge17$ and $0\le t\le1$. By \eqref{nine:harmonic-positive-cubic}, we obtain $\mathcal P(s,t)>0$. Combining this with \eqref{nine:harmonic-scalar-identity} and \eqref{nine:harmonic-middle-average} gives
\[
 W(x)\ge\frac{297}{25|x|^2},\qquad |x|^2\ge17.
\]

\subsubsection[The case 9 <= s <= 16]{The case $9\le s\le16$}
It remains to prove \eqref{nine:quadratic-exterior} for
$s=9,\ldots,16$. For $1\le k\le\lfloor\sqrt s\rfloor$,
set $Q_{s,k}=Q_s(k^2/s)$.
Substitution into \eqref{nine:quadratic-field} gives the following lower bounds.
In each row, $k$ ranges over the indicated admissible integers.
\begin{center}
\begin{tabular}{ccc}
  \toprule
  $s$&\shortstack{Lower bound for\\$\min_{k\ne2}Q_{s,k}$}&\shortstack{Lower bound for\\$Q_{s,2}$}\\\midrule
  $9,10,11,12,15,16$&$\lambda$&$\lambda$\\
  $13,14$&$\lambda+1/20$&$\lambda-1/300$\\\bottomrule
\end{tabular}
\end{center}

For \(s\in\{9,10,11,12,15,16\}\), we have \(Q_{s,k}\geq297/25\).
For \(s=13,14\), at most three coordinates have absolute value \(2\).
Their total weight in \eqref{nine:harmonic-middle-average} is at most $12/s\le12/13$.

For \(k\neq2\), we have \(Q_{s,k}\geq\lambda+1/20\). Therefore
\[
  sW(x)\ge\lambda+\frac1{13}\frac1{20}-\frac{12}{13}\frac1{300}
  =\lambda+\frac1{1300}>\lambda.
\]

This proves \eqref{nine:quadratic-exterior} for \(9\leq|x|^2\leq16\).

\subsubsection[The vertices with 0 < |x| squared < 9]{The vertices with \(0<|x|^2<9\)}
It remains to prove \(W(x)\geq\lambda/|x|^2\) for
\(x\in C=\{0<|x|^2<9\}\). We retain \eqref{nine:quadratic-field} when \(x\in C\) and \(y\notin C\),
and determine \(a_{xy}\) when \(x,y\in C\). Under coordinate permutations and sign
changes, \(C\) has the following fourteen orbits:
\[
  \Omega_j=[1^j]\quad(1\leq j\leq8),\qquad
  A_k=[2,1^k]\quad(0\leq k\leq4),\qquad B=[2,2],
\]
where omitted coordinates are zero. On the edges not contained
in \(C\), we keep \eqref{nine:quadratic-field} and denote
\(a_{xy}\) by \(a(s,t)\), where \(s=|x|^2\), \(t=|y|^2\).
For \(x\in C\), set
\[
  d(x)=18-\frac\lambda{|x|^2}
  -\sum_{y\sim x,\ y\notin C,\ y\neq0}a(|x|^2,|y|^2).
\]
Counting the neighbors of each orbit gives
\begin{equation}
  \label{eq:core-diagonals}
  \begin{aligned}
    d_{\Omega_j}&=18-\lambda/j
    -\mathbf1_{\{j>5\}}j\,a(j,j+3)
    -2\mathbf1_{\{j=8\}}a(j,j+1),\\
    d_{A_k}&=18-\lambda/s-a(s,s+5)
    -(k-\mathbf1_{\{k=1\}})a(s,s+3)
    -8\mathbf1_{\{k=4\}}a(s,s+1),\quad s=k+4,\\
    d_B&=18-\lambda/8-2a(8,13)-14a(8,9).
  \end{aligned}
\end{equation}
The orbit \(\Omega_j\) has \(j\) neighbors in \(\Omega_{j-1}\),
\(2(9-j)\) in \(\Omega_{j+1}\), and \(j\) in \(A_{j-1}\)
whenever the named orbit lies in \(C\). The orbit \(A_{j-1}\) has one
neighbor in \(\Omega_j\), \(j-1\) in \(A_{j-2}\), and
\(2(9-j)\) in \(A_j\). There is also the edge between
\(A_1\) and \(B\), with multiplicities one and two.

To obtain \(W(x)\geq\lambda/|x|^2\) on \(C\), choose
\(v:C\to(0,\infty)\) to be constant on each of \(\Omega_j,A_k\), and \(B\),
with the values in the following table.
\begin{center}
\begin{tabular}{cr}
\toprule
\(K\)&\(v(x)\), \(x\in K\)\\\midrule
\(\Omega_1\)&141223\\
\(\Omega_2\)&52640\\
\(\Omega_3\)&23528\\
\(\Omega_4\)&12641\\
\(\Omega_5\)&7739\\
\(\Omega_6\)&5123\\
\(\Omega_7\)&3623\\
\(\Omega_8\)&2691\\
\(A_0\)&22044\\
\(A_1\)&11500\\
\(A_2\)&6907\\
\(A_3\)&4623\\
\(A_4\)&3345\\
\(B\)&4685\\\bottomrule
\end{tabular}
\end{center}
By \eqref{eq:core-diagonals} and the numbers of neighbors listed above, we obtain
\[
  d(x)v(x)-\sum_{y\sim x,\,y\in C}v(y)>\frac1{10},\qquad x\in C.
\]
Set
\[
  a_{xy}=\frac{v(y)}{v(x)},\qquad x,y\in C,\quad x\sim y.
\]
For \(x,y\in\mathbb Z^9\setminus\{0\}\) satisfying
\[
  x\sim y,\qquad \max\{|x|^2,|y|^2\}\ge9,
\]
we retain the definition of \(a_{xy}\) in \eqref{nine:quadratic-field}.
Then \(a_{xy}>0\), \(a_{xy}a_{yx}=1\), and
\[
  W(x)-\frac\lambda{|x|^2}
  =\frac{d(x)v(x)-\sum_{y\sim x,\,y\in C}v(y)}{v(x)}>0,
  \qquad x\in C.
\]
Combining this with \eqref{nine:quadratic-exterior}, we obtain
\(W(x)\ge\lambda/|x|^2\) for every \(x\ne0\).
By \eqref{eq:nine-edge-weight}, we deduce that \(c_9\ge297/25\).
This completes the proof.
\end{proof}

\subsection{The upper bound}
\label{subsec:nine-upper-bound}
 To prove \(c_9<1440/121\), we construct a function \(u\) satisfying
\[
  \sum_{x\in\mathbb Z^9}|Du(x)|^2
  <\frac{1440}{121}\sum_{x\in\mathbb Z^9\setminus\{0\}}
  \frac{|u(x)|^2}{|x|^2}.
\]
 To apply the definition of \(c_9\), we then approximate \(u\) by finitely supported functions
that vanish at the origin.

\begin{proof}[Proof of the upper bound in Proposition~\ref{prop:nine-bounds}]
Set $\tau=1/12$ and $\lambda=1440/121$. For $x\ne0$, set $s=|x|^2$ and $P=\sum_i x_i^4$, and set
\begin{equation}
  \label{nine:trial-function}
  \begin{split}
    \mathcal B(x)&=\beta_\tau(s):=
    \frac{1-(1+\tau s)e^{-\tau s}}{s^2}
    =\int_0^\tau b e^{-bs}\,db,\\
    G(x)&=\frac{32}{15}e^{-4s/3}
    +\left(\frac{13}{100}+\frac{P}{550}\right)e^{-3s/8}
    +\frac1{70}e^{-s/8}+\frac{3P}{17}e^{-s},\\
    u(x)&=\mathcal B(x)+G(x).
  \end{split}
\end{equation}

We set \(\mathcal B(0)=G(0)=u(0)=0\).  To verify the strict
inequality above for the function \(u\) in \eqref{nine:trial-function}, write
\[
  Q(f,g)=\sum_x Df(x)\cdot Dg(x)-\lambda\sum_{x\neq0}\frac{f(x)g(x)}{|x|^2},\qquad Q(f)=Q(f,f).
\]
 It suffices to prove
\begin{equation}
  \label{nine:elementary-final-bounds}
  Q(u)<0.
\end{equation}

Let $G=\sum_{i=1}^4G_i$, where $G_i(x)=(c_i+d_iP)e^{-a_is}$ for $x\ne0$ and $G_i(0)=0$. The coefficients are
\begin{equation}
  \label{nine:grouped-gaussians}
  \begin{array}{c|ccc}
    i&a_i&c_i&d_i\\\hline
    1&4/3&32/15&0\\
    2&3/8&13/100&1/550\\
    3&1/8&1/70&0\\
    4&1&0&3/17
  \end{array}
\end{equation}
Also set $g_{a,\ell}(x)=P(x)^\ell e^{-a|x|^2}$ for $x\ne0$, and $g_{a,\ell}(0)=0$.
The decomposition
\[
  Q(u)=Q(G)+2Q(G,\mathcal B)+Q(\mathcal B)
\]
reduces the proof of \eqref{nine:elementary-final-bounds} to estimating \(Q(G)\), \(2Q(G,\mathcal B)\), and
\(Q(\mathcal B)\).

\subsubsection[Estimates for the Gaussian sums]{Estimates for \(\sum_x|DG(x)|^2\)}
To estimate \(\sum_x|DG(x)|^2\), we first express
\(\sum_x DG_i(x)\cdot DG_j(x)\) as products of one-dimensional sums. Set $q=e^{-1/24}$ and
\[
  (\rho_1,\rho_2,\rho_3,\rho_4)=(q^{32},q^9,q^3,q^{24}),
  \qquad f_i(k)=\rho_i^{k^2}.
\]

We use the values \(f_i(k)\), \(1\leq i\leq4\), \(0\leq k\leq K_i\), where
\[
  (K_1,K_2,K_3,K_4)=(4,6,7,4).
\]

Starting from $f_i(0)=1$ and $b_i(0)=\rho_i$,
we obtain these twenty-five values by
\[
  f_i(k+1)=f_i(k)b_i(k),\qquad b_i(k+1)=\rho_i^2b_i(k).
\]

Let $\Delta v(k)=v(k+1)-v(k)$.

By definition, we have
\begin{equation}
  \label{nine:factored-difference}
  \Delta(k^pf_i(k))=f_i(k)\{(k+1)^pb_i(k)-k^p\},
  \qquad p\in\{0,4\},\quad k\ge0.
\end{equation}

For each pair $i,j$, set $t=a_i+a_j$ and define
\[
  M_r=\sum_{k\in\mathbb Z}k^r f_i(k)f_j(k)\quad(r=0,4,8),
  \qquad \theta=M_0,
\]
\[
  D_{pr}=\sum_{k\in\mathbb Z}\Delta\bigl(k^p f_i(k)\bigr)
  \Delta\bigl(k^r f_j(k)\bigr)
  \quad(p,r\in\{0,4\}).
\]
We use the convention $k^0=1$, including at $k=0$.

For \(x=(k,x')\in\mathbb Z\times\mathbb Z^8\), set \(P'(x')=\sum_{j=1}^8(x'_j)^4\).
On $\mathbb Z^9$, set $\widetilde G_i=(c_i+d_iP)e^{-a_is}$.
For \(x=(k,x')\), we have
\[
  \widetilde G_i(k+1,x')-\widetilde G_i(k,x')=e^{-a_i|x'|^2}
  \{(c_i+d_iP')\Delta f_i+d_i\Delta(k^4f_i)\}.
\]
By the definitions of \(M_r\), \(\theta\), and \(P'\), we have
\[
  \begin{aligned}
    \sum_{x'\in\mathbb Z^8}e^{-t|x'|^2}&=\theta^8,\qquad
    \sum_{x'\in\mathbb Z^8}P'(x')e^{-t|x'|^2}=8M_4\theta^7,\\
    \sum_{x'\in\mathbb Z^8}P'(x')^2e^{-t|x'|^2}
    &=8M_8\theta^7+56M_4^2\theta^6.
  \end{aligned}
\]
Set
\begin{align*}
  X_0&=c_ic_jD_{00}+c_id_jD_{04}+d_ic_jD_{40}+d_id_jD_{44},\\
  X_1&=(c_id_j+d_ic_j)D_{00}+d_id_j(D_{04}+D_{40}),\\
  X_2&=d_id_jD_{00}{.}
\end{align*}

Since
\(\widetilde G_i(x_1,\ldots,x_9)=\widetilde G_i(x_{\sigma(1)},\ldots,x_{\sigma(9)})\)
for every permutation \(\sigma\) of \(\{1,\ldots,9\}\), we obtain
\begin{equation}
  \label{nine:grouped-energy}
  \begin{split}
    \sum_x DG_i(x)\cdot DG_j(x)&=9\theta^6\{\theta^2X_0+8\theta M_4X_1
    +(8\theta M_8+56M_4^2)X_2\}+O_{ij},\\
    \sum_x |DG(x)|^2&=\sum_{1\le i\le j\le4}(2-\delta_{ij})\sum_x DG_i(x)\cdot DG_j(x).
  \end{split}
\end{equation}
By \(G_i=\widetilde G_i-c_i\delta_0\) and
\(\sum_{x\in\mathbb Z^9}|D\delta_0(x)|^2=18\), we have
\[
  O_{ij}=-c_i\sum_{x\in\mathbb Z^9}D\widetilde G_j(x)\cdot D\delta_0(x)
  -c_j\sum_{x\in\mathbb Z^9}D\widetilde G_i(x)\cdot D\delta_0(x)+18c_ic_j.
\]
Since
\[
  \sum_{x\in\mathbb Z^9}D\widetilde G_i(x)\cdot D\delta_0(x)
  =(\Delta_{\mathrm d}\widetilde G_i)(0)
  =18\{c_i-(c_i+d_i)e^{-a_i}\},
\]
we obtain
\[
  O_{ij}=18\{c_ic_j(e^{-a_i}+e^{-a_j}-1)+c_id_je^{-a_j}+d_ic_je^{-a_i}\}.
\]

\paragraph[Estimates for Mr and Dpr.]{Estimates for \(M_r\) and \(D_{pr}\).}
To estimate \(M_r\) and \(D_{pr}\) using finite sums, for each pair \(i,j\), choose
\[
  K=\begin{cases}
    7,&t=1/4,\\
    6,&1/2\le t<1,\\
    4,&t\ge1.
  \end{cases}
\]
These cases cover all pairings used here. Since \(k^pf_i(k)\) is even, we have
\(\Delta(k^pf_i(k))|_{k=-j-1}=-\Delta(j^pf_i(j))\), so
\(\Delta(k^pf_i(k))\Delta(k^rf_j(k))\) is invariant under \(k\mapsto-k-1\). Set
\[
  \widehat M_r=\mathbf1_{\{r=0\}}
  +2\sum_{k=1}^{K-1}k^r f_i(k)f_j(k),
\]
\[
  \widehat D_{pr}=2\sum_{k=0}^{K-1}f_i(k)f_j(k)
  \{(k+1)^pb_i(k)-k^p\}\{(k+1)^rb_j(k)-k^r\},
  \qquad \varepsilon_r(t)=4K^re^{-K^2t}.
\]
Terms whose coefficients vanish in \eqref{nine:grouped-energy} are omitted.
Then
\begin{equation}
  \label{nine:elementary-gaussian-sum}
  \widehat M_r\le M_r\le\widehat M_r+\varepsilon_r(t),\qquad
  \widehat D_{pr}\le D_{pr}\le\widehat D_{pr}+\varepsilon_{p+r}(t).
\end{equation}

If the term \(\Delta(k^4f_i(k))\) occurs, then \(d_i\neq0\) and \(a_i\geq3/8\). For $k\ge4$,
\begin{equation}
  \label{nine:quartic-gaussian-ratio}
  \frac{(k+1)^4e^{-a_i(k+1)^2}}{k^4e^{-a_ik^2}}
  \le\left(\frac54\right)^4e^{-27/8}<1.
\end{equation}
For $p=0$, the functions $k^pf_i(k)$ decrease for $k\geq0$. Hence, for \(k\geq K\),
\[
  -k^pf_i(k)\leq\Delta(k^pf_i(k))\leq0,\qquad
  0\leq\Delta(k^pf_i(k))\Delta(k^rf_j(k))\leq k^{p+r}e^{-tk^2}.
\]
For \(0\leq r\leq8\) and \(k\geq K\), we have
\begin{equation}
  \label{nine:gaussian-tail-ratio}
  \frac{(k+1)^re^{-t(k+1)^2}}{k^re^{-tk^2}}
  \leq\left(\frac{K+1}K\right)^8e^{-(2K+1)t}<\frac12.
\end{equation}
To prove \eqref{nine:gaussian-tail-ratio}, it suffices to take \((K,t)=(7,1/4),(6,1/2),(4,1)\). Equations
\eqref{nine:quartic-gaussian-ratio} and \eqref{nine:gaussian-tail-ratio} follow from \(e^x>1+x+x^2/2\) for
\(x>0\).
By \eqref{nine:gaussian-tail-ratio}, we obtain
\[
  2\sum_{k\geq K}k^re^{-tk^2}
  \leq2K^re^{-K^2t}\sum_{j=0}^\infty2^{-j}
  =4K^re^{-K^2t}.
\]
Here \(e^{-K^2t}=f_i(K)f_j(K)\), with \(f_i(K)=q^{24a_iK^2}\).

Finally, we use
\begin{equation}
  \label{nine:elementary-common-q}
  \sum_{j=0}^{5}\frac{(-1/24)^j}{j!}<q
  <\sum_{j=0}^{6}\frac{(-1/24)^j}{j!}.
\end{equation}

By \eqref{nine:elementary-common-q}, \eqref{nine:factored-difference}, and
\eqref{nine:elementary-gaussian-sum}, we obtain rational intervals for
\(M_r\) and \(D_{pr}\), and hence for \(\sum_x DG_i(x)\cdot DG_j(x)\) by
\eqref{nine:grouped-energy}. Lower endpoints are rounded down and upper endpoints up,
including when \(D_{04}<0\) or \(D_{40}<0\).

Next, we estimate \(\sum_{x\neq0}G_i(x)G_j(x)/|x|^2\).
By the definition of \(G_i\), we have
\begin{equation}
  \label{nine:gaussian-hardy-products}
  \sum_{x\neq0}\frac{G_i(x)G_j(x)}{|x|^2}=c_ic_jF_0(t)
  +(c_id_j+d_ic_j)F_1(t)+d_id_jF_2(t),
\end{equation}
where $F_\ell(t)=\sum_{x\ne0}P(x)^\ell e^{-t|x|^2}/|x|^2$.

Since \(c_ic_j\geq0\), \(c_id_j+d_ic_j\geq0\), and \(d_id_j\geq0\), we may substitute lower bounds for
\(F_\ell(t)\) in \eqref{nine:gaussian-hardy-products}.

\subsubsection[Estimates for Fell(t)]{Estimates for \(F_\ell(t)\)}
We first estimate \(F_\ell(t)\) from below by retaining the vertices with
$1\leq|x|^2\leq14$. On these vertices, each coordinate has absolute
value at most $3$, and at most one has absolute value $3$.
If $i,j,k$ coordinates have absolute values $1,2,3$, respectively, then
\[
  |x|^2=i+4j+9k,\qquad P(x)=i+16j+81k.
\]
Hence, for \(1\leq n\leq14\),
\begin{equation*}
  r_\ell(n):=\sum_{\substack{x\in\mathbb Z^9\\|x|^2=n}}P(x)^\ell=\sum_{\substack{i+4j+9k=n\\i,j,k\ge0,\
  i+j+k\le9}}
  2^{i+j+k}\binom9k\binom{9-k}j\binom{9-k-j}i
  (i+16j+81k)^\ell.
\end{equation*}
For each $n\le14$, at most four triples occur. For example, at $n=14$ they are
$(6,2,0)$, $(2,3,0)$, $(5,0,1)$, and $(1,1,1)$, giving
\[
  r_\ell(14)=64512\cdot38^\ell+40320\cdot50^\ell
  +32256\cdot86^\ell+4032\cdot98^\ell.
\]

Set
\[
  F_\ell^{[14]}(t)=\sum_{n=1}^{14}\frac{r_\ell(n)}n e^{-nt}.
\]

For \(t\geq1\), we use \(F_\ell(t)\geq F_\ell^{[14]}(t)\), since \(P(x)^\ell e^{-t|x|^2}/|x|^2\geq0\).
For the values of \(t\) used in \(F_\ell^{[14]}(t)\), we have \(24t\in\mathbb N\).
To evaluate \(F_\ell^{[14]}(t)\) as a polynomial in \(z=q^{24t}\), set \(p_{15}=0\) and
\begin{equation}
  \label{nine:finite-hardy-recurrence}
  p_n=z\left(\frac{r_\ell(n)}n+p_{n+1}\right)
  \quad(n=14,\ldots,1),\qquad F_\ell^{[14]}(t)=p_1.
\end{equation}

For \(0<t<1\), we obtain a lower bound for \(F_\ell(t)\)
by estimating \(-F_\ell'(t)\) and integrating over \([t,1]\).
For this purpose, set \(T(t)=\sum_{k\in\mathbb Z}e^{-tk^2}\) and
\(M_j(t)=\sum_{k\in\mathbb Z}k^j e^{-tk^2}\).
By \cite[(1.8.16)]{DLMF}, we have
\begin{equation}
  \label{nine:theta-lower}
  T(t)=\sqrt{\frac\pi t}\left(1+2\sum_{k=1}^\infty e^{-\pi^2k^2/t}\right)
  \geq\sqrt{\frac\pi t}.
\end{equation}
Differentiating with respect to \(t\), we obtain
\[
  M_4(t)=\frac{d^2T}{dt^2}(t),\qquad
  M_8(t)=\frac{d^4T}{dt^4}(t).
\]
For \(k\geq1\), set \(z=\pi^2k^2/t>9\).
We have
\[
  \begin{aligned}
    \frac{d^2}{dt^2}(t^{-1/2}e^{-\pi^2k^2/t})&=t^{-5/2}e^{-z}q_4(z),\\
    \frac{d^4}{dt^4}(t^{-1/2}e^{-\pi^2k^2/t})&=t^{-9/2}e^{-z}q_8(z),
  \end{aligned}
\]
where
\[
  q_4(z)=z^2-3z+\frac34>0,\qquad
  q_8(z)=z^4-14z^3+\frac{105}{2}z^2-\frac{105}{2}z+\frac{105}{16}.
\]
For \(z=9+w\), we obtain
\[
  q_8(9+w)=w^4+22w^3+\frac{321}{2}w^2
  +\frac{813}{2}w+\frac{2265}{16}>0
  \quad(w\ge0).
\]
Hence, for \(0<t\leq1\),
\begin{equation}
  \label{nine:moment-lower}
  M_4(t)\ge\frac3{4t^2}\sqrt{\frac\pi t},\qquad
  M_8(t)\ge\frac{105}{16t^4}\sqrt{\frac\pi t}.
\end{equation}

Set $\eta_\ell=7/2+2\ell$, $C_0'=1$, $C_1'=27/4$, and $C_2'=1593/16$.
Combining \eqref{nine:theta-lower} and \eqref{nine:moment-lower} with
\[
  \begin{aligned}
    -F_0'(t)&=T(t)^9-1,\\
    -F_1'(t)&=9M_4(t)T(t)^8,\\
    -F_2'(t)&=9M_8(t)T(t)^8+72M_4(t)^2T(t)^7,
  \end{aligned}
\]
and integrating over $[t,1]$, we obtain
\begin{equation}
  \label{nine:hardy-tail}
  F_\ell(t)\ge F_\ell^{[14]}(1)
  +\frac{C_\ell'\pi^{9/2}}{\eta_\ell}(t^{-\eta_\ell}-1)
  +\mathbf1_{\{\ell=0\}}(t-1)
  \quad(0<t<1).
\end{equation}
This holds for $\ell=0,1,2$. Since the coefficients of $G$ are nonnegative,
summing over all pairs gives a lower bound for
$\sum_{x\neq0}|G(x)|^2/|x|^2$.
For $a=1/8$ or $3/8$,
we also have
\[
  \sum_{x\neq0}\frac{g_{a,\ell}(x)\mathcal B(x)}{|x|^2}=\int_0^\tau bF_\ell(a+b)\,db.
\]
Since $a+\tau<1$, \eqref{nine:hardy-tail} gives
\[
  \begin{aligned}
    \sum_{x\neq0}\frac{g_{a,\ell}(x)\mathcal B(x)}{|x|^2}
    &\geq\frac{\tau^2}{2}F_\ell^{[14]}(1)
    +\frac{C_\ell'\pi^{9/2}}{\eta_\ell}
    \left(\int_0^\tau b(a+b)^{-\eta_\ell}\,db-\frac{\tau^2}{2}\right)\\
    &\quad+\mathbf1_{\{\ell=0\}}
    \left(\frac{(a-1)\tau^2}{2}+\frac{\tau^3}{3}\right).
  \end{aligned}
\]
We evaluate \(\int_0^\tau b(a+b)^{-\eta_\ell}\,db\) using \eqref{nine:elementary-stable-integral} and
\eqref{nine:positive-integral-recurrence} below.

\subsubsection[Estimates for Q(B)]{Estimates for \(Q(\mathcal B)\)}
We first evaluate the integrals needed to estimate \(\sum_x|D\mathcal B(x)|^2\). For $k=1,2,3$, define
\[
  J_k=\int_0^\tau\int_0^\tau
  \frac{a^{k+1}b^{k+1}}{(a+b)^{k+9/2}}\,da\,db.
\]

Since the integrand is unchanged under \(a\leftrightarrow b\), we integrate over
\(0<b<a<\tau\) and multiply by two. With \(b=av\) and then
\(w=(1+v)^{-1/2}\), we obtain
\begin{equation}
  \label{nine:radial-Jk-integral}
  \begin{aligned}
    J_k&=2\int_0^\tau a^{k-3/2}\,da
    \int_0^1\frac{v^{k+1}}{(1+v)^{k+9/2}}\,dv\\
    &=\frac{4\tau^{k-1/2}}{k-1/2}
    \int_{1/\sqrt2}^{1}w^4(1-w^2)^{k+1}\,dw.
  \end{aligned}
\end{equation}

Taking \(k=1,2\) in \eqref{nine:radial-Jk-integral}, we obtain
\[
  \begin{aligned}
    J_1&=\left(\frac{64}{315}-\frac{107\sqrt2}{1260}\right)\tau^{1/2},\\
    J_2&=\left(\frac{128}{3465}-\frac{533\sqrt2}{27720}\right)\tau^{3/2}.
  \end{aligned}
\]

By \eqref{nine:radial-Jk-integral} and \(1-w^2\leq1/2\) for \(1/\sqrt2\leq w\leq1\), we obtain
\(J_3\leq(3\tau/10)J_2\).

To estimate \(\sum_{x\in\mathbb Z^9}|D\mathcal B(x)|^2\),
we use \eqref{nine:gaussian-nine-dimensional-error} to compare
\[
  \sum_{x\in\mathbb Z^9}e^{-t|x+\rho|^2}
  \quad\text{with}\quad
  \int_{\mathbb R^9}e^{-t|X+\rho|^2}\,dX
  =\left(\frac\pi t\right)^{9/2},
\]
for \(0<t\leq2\tau\) and \(\rho\in\mathbb R^9\).
For \(a,b>0\), direct integration gives
\[
  \begin{aligned}
    &\sum_{i=1}^9\int_{\mathbb R^9}
    (e^{-a|x+e_i|^2}-e^{-a|x|^2})(e^{-b|x+e_i|^2}-e^{-b|x|^2})\,dx\\
    &=18\left(\frac\pi{a+b}\right)^{9/2}
    \left(1-e^{-ab/(a+b)}\right).
  \end{aligned}
\]
Since $1-e^{-z}\leq z-z^2/2+z^3/6$ for $z\geq0$,
\begin{equation}
  \label{nine:radial-continuous-energy}
  \begin{aligned}
    &18\pi^{9/2}\int_0^\tau\int_0^\tau
    \frac{ab}{(a+b)^{9/2}}\left(1-e^{-ab/(a+b)}\right)\,da\,db\\
    &\leq18\pi^{9/2}\left(J_1-\frac12J_2+\frac16J_3\right)\\
    &\leq18\pi^{9/2}\left(J_1-\left(\frac12-\frac\tau{20}\right)J_2\right)
    =18\pi^{9/2}\left(J_1-\frac{119}{240}J_2\right).
  \end{aligned}
\end{equation}

We shall prove
\begin{equation}
  \label{nine:tail-energy-upper}
  \begin{split}
    \sum_x |D(\mathcal B)(x)|^2\le U_B:={}&
    18\pi^{9/2}\left(J_1-\frac{119}{240}J_2\right)\\
    &+18\left(\tau^2\beta_\tau(1)-\frac{\tau^4}{4}\right)
    +360\tau^4\left(\frac\pi{2\tau}\right)^{9/2}e^{-\pi^2/(2\tau)}.
  \end{split}
\end{equation}
The double integral in \eqref{nine:radial-continuous-energy} is finite, since
\[
  0\leq\frac{ab}{(a+b)^{9/2}}(1-e^{-ab/(a+b)})
  \leq\frac{a^2b^2}{(a+b)^{11/2}},\qquad J_1<\infty.
\]
For \(t=a+b\leq2\tau\), by \cite[(1.8.16)]{DLMF}, we have
\begin{equation}
  \label{nine:gaussian-one-dimensional-error}
  \left|\sum_{k\in\mathbb Z}e^{-t(k+\theta)^2}-\left(\frac\pi t\right)^{1/2}\right|
  \leq\varepsilon\left(\frac\pi t\right)^{1/2},\qquad
  \varepsilon:=4e^{-\pi^2/t}<\frac1{1000},\quad\theta\in\mathbb R.
\end{equation}
Applying \eqref{nine:gaussian-one-dimensional-error} in each of the nine coordinates, for \(\rho\in\mathbb
R^9\), we obtain
\begin{equation}
  \label{nine:gaussian-nine-dimensional-error}
  \left|\sum_{x\in\mathbb Z^9}e^{-t|x+\rho|^2}
  -\left(\frac\pi t\right)^{9/2}\right|
  \leq((1+\varepsilon)^9-1)\left(\frac\pi t\right)^{9/2}
  \leq10\varepsilon\left(\frac\pi t\right)^{9/2}.
\end{equation}
Applying \eqref{nine:gaussian-nine-dimensional-error} with \(\rho=0\) and \(\rho=(b/t)e_i\), and using
\(e^{-ab/t}\leq1\), gives
\[
  \begin{aligned}
    \Bigg|&\sum_{i=1}^9\sum_{x\in\mathbb Z^9}
    (e^{-a|x+e_i|^2}-e^{-a|x|^2})(e^{-b|x+e_i|^2}-e^{-b|x|^2})\\
    &-18\left(\frac\pi t\right)^{9/2}(1-e^{-ab/t})\Bigg|
    \leq1440\left(\frac\pi t\right)^{9/2}e^{-\pi^2/t}.
  \end{aligned}
\]
Since
\begin{equation}
  \label{eq:gaussian-tail-derivative}
  \frac{d}{dt}(t^{-\nu}e^{-\pi^2/t})
  =t^{-\nu-2}e^{-\pi^2/t}(\pi^2-\nu t)>0
  \qquad(0<t<\pi^2/\nu),
\end{equation}
with \(\nu=9/2\), we obtain
\begin{equation}
  \label{nine:radial-discrete-error}
  1440\int_0^\tau\int_0^\tau ab\left(\frac\pi{a+b}\right)^{9/2}
  e^{-\pi^2/(a+b)}\,da\,db
  \leq360\tau^4\left(\frac\pi{2\tau}\right)^{9/2}e^{-\pi^2/(2\tau)}.
\end{equation}
Finally, set \(\widetilde{\mathcal B}(x)=\int_0^\tau b e^{-b|x|^2}\,db\) for \(x\in\mathbb Z^9\). Since
\(\mathcal B=\widetilde{\mathcal B}-(\tau^2/2)\delta_0\), we have
\begin{equation}
  \label{nine:radial-origin-correction}
  \begin{aligned}
    \sum_x|D\mathcal B(x)|^2-\sum_x|D\widetilde{\mathcal B}(x)|^2
    &=18\left[\beta_\tau(1)^2-\left(\beta_\tau(1)-\frac{\tau^2}{2}\right)^2\right]\\
    &=18\left(\tau^2\beta_\tau(1)-\frac{\tau^4}{4}\right).
  \end{aligned}
\end{equation}
Combining \eqref{nine:radial-origin-correction} with \eqref{nine:radial-continuous-energy} and
\eqref{nine:radial-discrete-error} proves \eqref{nine:tail-energy-upper}.

Next, to estimate
\(\sum_{x\neq0}\frac{|\mathcal B(x)|^2}{|x|^2}\), set
\[
  C_*=F_0^{[14]}(1)-\frac{2\pi^{9/2}}7-1.
\]
Equation~\eqref{nine:hardy-tail} gives $F_0(t)\ge(2\pi^{9/2}/7)t^{-7/2}+t+C_*$
for $0<t\le2\tau$. Moreover,
\[
  \int_0^\tau\int_0^\tau ab(a+b)^{-7/2}\,da\,db
  =\frac{16-7\sqrt2}{15}\sqrt\tau.
\]

Consequently,
\begin{equation*}
  \sum_{x\neq0}\frac{|\mathcal B(x)|^2}{|x|^2}\ge L_B:=
  \frac{2\pi^{9/2}}7\frac{16-7\sqrt2}{15}\sqrt\tau
  +\frac{\tau^5}{3}+\frac{\tau^4}{4}C_*.
\end{equation*}
By \(\mathcal B(x)=\int_0^\tau b e^{-b|x|^2}\,db\) for \(x\neq0\) and Tonelli's theorem,
\[
  \sum_{x\neq0}\frac{|\mathcal B(x)|^2}{|x|^2}
  =\int_0^\tau\int_0^\tau abF_0(a+b)\,da\,db.
\]
Since \(D_jg_{a,0}(x)D_jg_{b,0}(x)\geq0\), Tonelli's theorem gives
\[
  \sum_x |D(\mathcal B)(x)|^2
  =\int_0^\tau\int_0^\tau ab\,\sum_x D(g_{a,0})(x)\cdot D(g_{b,0})(x)\,da\,db,
\]
and \eqref{nine:tail-energy-upper} gives \(\sum_x|D\mathcal B(x)|^2<\infty\).

\subsubsection[The cross terms for a=1/8 and a=3/8]{The cross terms for \(a=1/8\) and \(a=3/8\)}
For \(a=1/8\) or \(3/8\), we estimate
\(\sum_{x\in\mathbb Z^9}Dg_{a,\ell}(x)\cdot D\mathcal B(x)\)
by comparing \(\mathcal E_{\ell0}(a,b)\) with
\(\mathcal E_{\ell0}^{\rm c}(a,b)\) in \eqref{nine:continuous-cross-definition}.
For this comparison, we use \cite[(1.8.14)]{DLMF}.

For \(a=1/8\) or \(3/8\), set \(\mathcal E_{\ell0}(a,b)=\sum_x D(g_{a,\ell})(x)\cdot D(g_{b,0})(x)\).
Let $t=a+b$, $z=ab/t$, and $\nu_a=1/2-a\tau/(6(a+\tau))$.

For \(\ell=0,1\), define
\[
  n_\ell=\left(\frac t\pi\right)^{9/2}\int_{\mathbb R^9}P(x)^\ell e^{-t|x|^2}\,dx,
  \qquad
  A_\ell=\left(\frac t\pi\right)^{9/2}
  \int_{\mathbb R^9}P(x)^\ell e^{-t|x+(b/t)e_1|^2}\,dx.
\]
By integration in each coordinate, we have
\[
  n_0=A_0=1,\qquad n_1=\frac{27}{4t^2},\qquad
  A_1=n_1+\frac{3b^2}{t^3}+\frac{b^4}{t^4}.
\]
For \(\ell=0,1\), define
\begin{equation}
  \label{nine:continuous-cross-definition}
  \mathcal E_{\ell0}^{\rm c}(a,b)
  =18\left(\frac\pi t\right)^{9/2}(n_\ell-A_\ell e^{-z})
  +18\{e^{-a}+\mathbf1_{\{\ell=0\}}(e^{-b}-1)\}.
\end{equation}
Since $z\le a\tau/(a+\tau)$, Taylor's inequality gives
$1-e^{-z}\le z-\nu_a z^2$. Dropping the contribution $-b^4e^{-z}/t^4$ and using $e^{-z}\ge1-z$
in \(-3b^2e^{-z}/t^3\), we obtain
\begin{equation}
  \label{nine:cross-taylor}
  n_1-A_1e^{-z}
  \le n_1(z-\nu_a z^2)-\frac{3b^2}{t^3}(1-z).
\end{equation}

To estimate \(\int_0^\tau b\mathcal E_{\ell0}^{\rm c}(a,b)\,db\), define
\[
  \mathcal J_p(a)=\int_0^\tau\frac{b^2}{(a+b)^p}\,db,\qquad
  \mathcal K_p(a)=\int_0^\tau\frac{b^3}{(a+b)^{p+1}}\,db.
\]

By \eqref{nine:continuous-cross-definition}, \eqref{nine:cross-taylor}, and \(b^4=tb^3-ab^3\), we obtain
$\int_0^\tau b\mathcal E_{\ell0}^{\rm c}(a,b)\,db\le C_\ell(a)$, where
\begin{equation*}
  \begin{split}
    C_0(a)={}&18\pi^{9/2}
    \{a\mathcal J_{11/2}(a)-\nu_a a^2\mathcal K_{11/2}(a)\}+O_0(a),\\
    C_1(a)={}&18\pi^{9/2}\Bigl\{
    \frac{27a}{4}\mathcal J_{15/2}(a)
    -a^2\Bigl(\frac{27\nu_a}{4}+3\Bigr)\mathcal K_{15/2}(a)\\
    &\hspace{26mm}-3(1-a)\mathcal K_{13/2}(a)\Bigr\}+O_1(a).
  \end{split}
\end{equation*}
The corrections at the origin are
\[
  O_0(a)=18\left\{(e^{-a}-1)\frac{\tau^2}{2}
  +1-(1+\tau)e^{-\tau}\right\},\qquad
  O_1(a)=9\tau^2e^{-a}.
\]

We next estimate \(|\mathcal E_{\ell0}(a,b)-\mathcal E_{\ell0}^{\rm c}(a,b)|\). Since $t=a+b\le11/24$,
it suffices to work with $0<t\le1/2$.
Set $w=\pi^2/t>18$. For $n\ge1$, the number $r_0(n)$ of points with $|k|^2=n$ satisfies
\[
  r_0(1)=18,\qquad r_0(n)\le2^n\binom{n+8}{8}.
\]
Since \(e^{-w}<10^{-4}\), for \(n\geq2\)
and \(0\leq j\leq4\), we have
\[
  \frac{2^{n+1}\binom{n+9}{8}(n+1)^je^{-w(n+1)}}
  {2^n\binom{n+8}{8}n^je^{-wn}}
  =2e^{-w}\frac{n+9}{n+1}\left(\frac{n+1}{n}\right)^j
  \leq2\cdot10^{-4}\frac{11}{3}\left(\frac32\right)^4<\frac1{100}.
\]
Consequently,
\begin{equation}
  \label{nine:small-scale-frequencies}
  \sum_{k\ne0}|k|^{2j}e^{-w|k|^2}
  \le e^{-w}\left(18+\frac{2880\cdot10^{-4}}{1-1/100}\right)
  <20e^{-w},\qquad 0\le j\le4.
\end{equation}

Indeed,
\[
  \#\{(r_1,\ldots,r_9)\in\mathbb N_0^9:r_1+\cdots+r_9=n\}
  =\binom{n+8}{8},
\]
and each choice of the squares $x_i^2=r_i$ allows at most $2^n$ choices of signs.

To apply \eqref{nine:small-scale-frequencies} to the terms containing \(P\),
set \(\rho=b/t\leq1\) and \(z_i=-\rho\delta_{ij}+\mathrm i\pi k_i/t\).

Integration in each coordinate gives
\[
  \begin{aligned}
    &\left(\frac t\pi\right)^{9/2}
    \int_{\mathbb R^9}P(x)e^{-t|x+\rho e_j|^2}e^{2\pi\mathrm i k\cdot x}\,dx\\
    &=e^{-\pi^2|k|^2/t-2\pi\mathrm i\rho k_j}
    \sum_{i=1}^9\left(z_i^4+\frac{3z_i^2}{t}+\frac3{4t^2}\right).
  \end{aligned}
\]
For \(k\neq0\), we have
\begin{equation}
  \label{nine:shifted-fourier-bound}
  \left|\sum_{i=1}^9\left(z_i^4+\frac{3z_i^2}{t}+\frac3{4t^2}\right)\right|
  \leq 1+\frac3t+\frac{27}{4t^2}
  +\left(2+\frac3t\right)\frac{\pi^2|k|^2}{t^2}
  +\frac{\pi^4|k|^4}{t^4}
  \le\frac54\frac{\pi^4|k|^4}{t^4}.
\end{equation}
Indeed, $t\le1/2$, $\pi^2>9$, and $|k|\ge1$ give
\[
  \frac{t^4+3t^3+27t^2/4}{\pi^4|k|^4}\leq\frac{17}{648},\qquad
  \frac{2t^2+3t}{\pi^2|k|^2}\leq\frac29,\qquad
  1+\frac29+\frac{17}{648}<\frac54.
\]

Equation \eqref{nine:shifted-fourier-bound} also holds for \(\rho=0\). By \eqref{nine:small-scale-frequencies}
and \eqref{nine:shifted-fourier-bound}, we obtain
\[
  |\mathcal E_{\ell0}(a,b)-\mathcal E_{\ell0}^{\rm c}(a,b)|
  \le K_\ell\left(\frac\pi t\right)^{9/2+4\ell}e^{-\pi^2/t},
  \qquad K_0=720,\quad K_1=900.
\]

The correction \(18\{e^{-a}+\mathbf1_{\{\ell=0\}}(e^{-b}-1)\}\) occurs in both \(\mathcal E_{\ell0}(a,b)\) and
\eqref{nine:continuous-cross-definition}, so it cancels in their difference.
By \eqref{eq:gaussian-tail-derivative}, with \(\nu=9/2+4\ell\) and
\(0<t\leq11/24\), we obtain
\[
  \left|\int_0^\tau b\bigl(\mathcal E_{\ell0}(a,b)-\mathcal E_{\ell0}^{\rm c}(a,b)\bigr)\,db\right|
  \leq D_\ell:=\frac{K_\ell\tau^2}{2}
  \left(\frac{24\pi}{11}\right)^{9/2+4\ell}
  e^{-24\pi^2/11}.
\]

The elementary bounds $157/50<\pi<22/7$, $e>19/7$, and $e^{1/2}>3/2$ give
\[
  D_\ell<\frac{K_\ell}{288}
  \left(\frac{48}{7}\right)^{4+4\ell}
  \frac83\frac23\left(\frac7{19}\right)^{21},
  \qquad D_0<\frac1{100000},\quad D_1<\frac1{40}.
\]

Since \(2\sum_{i=2}^3(c_iD_0+d_iD_1)<1/10000\), we obtain
\begin{equation*}
  2\sum_x D(G_2+G_3)(x)\cdot D(\mathcal B)(x)
  \le2\sum_{i=2}^3\{c_iC_0(a_i)+d_iC_1(a_i)\}+\frac1{10000}.
\end{equation*}

\subsubsection[Estimates for the integrals]{Estimates for the integrals}
We now evaluate \(\mathcal J_p(a)\), \(\mathcal K_p(a)\), and \(\int_0^\tau b(a+b)^{-\eta_\ell}\,db\). First,
integration by parts gives
\[
  \tau^3(a+\tau)^{-p}
  =\bigl[b^3(a+b)^{-p}\bigr]_0^\tau
  =3\mathcal J_p(a)-p\mathcal K_p(a).
\]
Hence
\begin{equation}
  \label{nine:integration-by-parts}
  \mathcal K_p(a)=\frac{3\mathcal J_p(a)-\tau^3(a+\tau)^{-p}}p.
\end{equation}
It remains to evaluate $\mathcal J_p(a)$ for $p=11/2$ when $a=1/8$,
and for $p=11/2,13/2,15/2$ when $a=3/8$.
To
evaluate \(\mathcal J_p(a)\) and \(\int_0^\tau b(a+b)^{-\eta_\ell}\,db\), set
\[
  I_{q,p}(a)=\int_0^\tau b^q(a+b)^{-p}\,db,\qquad
  r=\sqrt{\frac a{a+\tau}},\qquad
  v=1-r=\frac{\tau}{(a+\tau)(1+r)}.
\]
With $w=\sqrt{a/(a+b)}$, we obtain
\begin{equation}
  \label{nine:elementary-stable-integral}
  \begin{aligned}
    I_{1,7/2}(a)&=\frac{2}{a^{3/2}}\int_r^1 w^2(1-w^2)\,dw\\
    &=\frac{2v^2}{15a^{3/2}}(2+4r+6r^2+3r^3),\\
    I_{2,11/2}(a)&=\frac{2}{a^{5/2}}\int_r^1 w^4(1-w^2)^2\,dw\\
    &=\frac{2v^3}{315a^{5/2}}
    (8+24r+48r^2+80r^3+120r^4+105r^5+35r^6).
  \end{aligned}
\end{equation}

Likewise,
\[
  \begin{aligned}
    \tau^{q+1}(a+\tau)^{-p}
    &=\bigl[b^{q+1}(a+b)^{-p}\bigr]_0^\tau\\
    &=(q+1)I_{q,p}(a)-p\bigl(I_{q,p}(a)-aI_{q,p+1}(a)\bigr),
  \end{aligned}
\]
and therefore
\begin{equation}
  \label{nine:positive-integral-recurrence}
  apI_{q,p+1}(a)=(p-q-1)I_{q,p}(a)
  +\tau^{q+1}(a+\tau)^{-p}.
\end{equation}

By \eqref{nine:positive-integral-recurrence}, we obtain
\(I_{1,11/2}\), \(I_{2,13/2}\), and \(I_{2,15/2}\) from
\eqref{nine:elementary-stable-integral}. Together with
\(\mathcal J_p=I_{2,p}\) and \eqref{nine:integration-by-parts}, this gives
\(\mathcal J_p(a)\), \(\mathcal K_p(a)\), and \(I_{1,\eta_\ell}(a)\) for the values of \(p,\ell,a\) used
above.

\subsubsection[The cross terms for a=4/3 and a=1]{The cross terms for \(a=4/3\) and \(a=1\)}
For these values of \(a\), we estimate
\(\sum_x D(g_{a,\ell})(x)\cdot D(\mathcal B)(x)\) directly.

Let \(f=g_{4/3,0}\) or \(g_{1,1}\), and take \(x,y\in\mathbb Z^9\setminus\{0\}\)
with \(x\sim y\) and \(|x|<|y|\). Set \(s=|x|^2\) and \(s'=|y|^2\).
We first verify that \(f(x)\geq f(y)\). For \(f=g_{4/3,0}\), this follows
from its definition. For \(f=g_{1,1}\), let \(j\) be the coordinate in which
\(x\) and \(y\) differ, and set \(r=|x_j|\), so \(|y_j|=r+1\). Then
\[
  P(y)-P(x)=(r+1)^4-r^4,\qquad s'-s=2r+1,
\]
and, since \(x\ne0\),
\[
  \frac{P(y)}{P(x)}\leq
  \begin{cases}2,&r=0,\\16,&r\geq1.\end{cases}
\]
Thus \(P(y)e^{-s'}\leq P(x)e^{-s}\), which gives \(f(y)\leq f(x)\).

Consequently, the function
\[
  M_f(b)=\sum_{\substack{x\sim y\\x,y\ne0,\ |x|<|y|}}
  (f(x)-f(y))\int_{|x|^2}^{|y|^2}e^{-br}\,dr
\]
satisfies, for \(j\geq0\) and \(b\geq0\),
\[
  (-1)^jM_f^{(j)}(b)
  =\sum_{\substack{x\sim y\\x,y\ne0,\ |x|<|y|}}
  (f(x)-f(y))\int_{|x|^2}^{|y|^2}r^je^{-br}\,dr\geq0.
\]
Indeed, for \(x\sim y\) and \(0<|x|<|y|\),
\[
  0\leq(f(x)-f(y))\int_{|x|^2}^{|y|^2}r^je^{-br}\,dr
  \leq C_j(1+|x|)^{2j+5}e^{-|x|^2},
  \qquad\sum_{x\in\mathbb Z^9}(1+|x|)^{2j+5}e^{-|x|^2}<\infty.
\]
Set $f_1:=f(e_1)$ and $N_j:=\sum_{x\ne0}|x|^{2j}f(x)$ for $j=0,1,2$. Then
\begin{equation}
  \label{nine:monotone-energy}
  \sum_x Df(x)\cdot D(g_{b,0})(x)=18f_1e^{-b}+bM_f(b),\qquad
  M_f(0)=18(N_0+f_1).
\end{equation}
Summation by parts and the identities
\[
  -\Delta_{\mathrm d}s=18,\qquad -\Delta_{\mathrm d}(s^2)=44s+18,\qquad
  -\Delta_{\mathrm d}(s^3)=78s^2+78s+18,
\]
also give
\[
  M_f'(0)=-22N_1-9(N_0+f_1),\qquad
  M_f''(0)=26(N_2+N_1)+6(N_0+f_1).
\]

Here \(\sum_{|x|=1}f(x)=18f_1\), since \(f(\pm e_j)=f_1\) for \(1\leq j\leq9\).

To obtain an upper bound for \(\int_0^\tau b^2M_f(b)\,db\),
we choose a cubic polynomial \(R_f\) satisfying
\[
  R_f^{(j)}(0)=M_f^{(j)}(0)\quad(j=0,1,2),\qquad
  R_f(1/8)=M_f(1/8).
\] By the remainder formula for polynomial interpolation and \(M_f^{(4)}\geq0\), for some \(\xi\in(0,1/8)\), we have
\[
  M_f(b)-R_f(b)=\frac{M_f^{(4)}(\xi)}{24}b^3(b-1/8)\le0,
  \qquad 0\le b\le\tau.
\]
Integrating $b^2R_f(b)$ and using $1/8=3\tau/2$ therefore gives
\begin{equation}
  \label{nine:fast-cross-hermite}
  \begin{split}
    \sum_x Df(x)\cdot D(\mathcal B)(x)\le{}&18f_1\beta_\tau(1)
    +\tau^3\left\{\frac{23}{81}M_f(0)
    +\frac{19\tau}{108}M_f'(0)\right.\\
    &\left.\hspace{18mm}
    +\frac{2\tau^2}{45}M_f''(0)
    +\frac4{81}M_f(1/8)\right\}.
  \end{split}
\end{equation}
By \eqref{nine:monotone-energy},
\[
  M_f(1/8)=8\{\sum_x Df(x)\cdot D(g_{1/8,0})(x)-18f_1e^{-1/8}\}.
\]
By \eqref{nine:grouped-gaussians}, we have
\begin{equation}
  \label{nine:fast-pairing-factors}
  \begin{aligned}
    \sum_x Dg_{4/3,0}\cdot Dg_{1/8,0}&=\frac{1050}{32}\sum_x DG_1\cdot DG_3,\\
    \sum_x Dg_{1,1}\cdot Dg_{1/8,0}&=\frac{1190}{3}\sum_x DG_4\cdot DG_3.
  \end{aligned}
\end{equation}
We estimate \(\sum_x DG_1\cdot DG_3\) and \(\sum_x DG_4\cdot DG_3\) by \eqref{nine:grouped-energy}.

To evaluate \(N_0,N_1,N_2\), set $T=T(a)$ and $M_j=M_j(a)$. For $f=g_{a,0}$,
\[
  N_0=T^9-1,\qquad N_1=9M_2T^8,\qquad
  N_2=9M_4T^8+72M_2^2T^7.
\]
For $f=g_{a,1}$,
\[
  \begin{aligned}
    N_0&=9M_4T^8,\qquad N_1=9M_6T^8+72M_4M_2T^7,\\
    N_2&=9M_8T^8+72M_4^2T^7+144M_6M_2T^7+504M_4M_2^2T^6.
  \end{aligned}
\]

Here we use \(s=\sum_{i=1}^9x_i^2\) and \(s^2=\sum_{i=1}^9x_i^4+2\sum_{i<j}x_i^2x_j^2\).
Since \(a\geq1\), by \eqref{nine:gaussian-tail-ratio} with
\(K=4\) and \(t=a\), we have
\[
  0\leq M_j(a)-\sum_{|k|\leq3}k^je^{-ak^2}
  \leq4\cdot4^je^{-16a},\qquad j=0,2,4,6,8,
\]
where \(M_0=T\).
We use \(e^{-ak^2}=q^{24ak^2}\) to evaluate the finite sums.

Since \(g_{a,\ell}(x)\geq0\) and
\(\mathcal B(x)=\beta_\tau(|x|^2)>0\) for \(x\neq0\), we have
\begin{equation*}
  \sum_{x\neq0}\frac{g_{a,\ell}(x)\mathcal B(x)}{|x|^2}\ge
  \sum_{n=1}^{14}\frac{r_\ell(n)}n e^{-an}\beta_\tau(n).
\end{equation*}
For \((a,\ell)=(4/3,0)\) and \((1,1)\), we use \(\beta_\tau(n)\), \(1\leq n\leq14\), with \(e^{-\tau
n}=q^{2n}\) and \(e^{-an}=q^{24an}\).

\subsubsection[Rational bounds for Q(u)]{Rational bounds for \(Q(u)\)}

 To prove \(Q(u)<0\) using rational bounds, we use \eqref{nine:elementary-common-q}
for \(q=e^{-1/24}\) and
\[
  \frac{355}{113}-\frac1{10^6}<\pi<\frac{355}{113}
\] for \(\pi\).

By Machin's identity \(\pi=16\arctan(1/5)-4\arctan(1/239)\) and
\[
  \arctan x=\sum_{j=0}^{M-1}\frac{(-1)^jx^{2j+1}}{2j+1}+R_M(x),
  \qquad 0\leq(-1)^MR_M(x)\leq\frac{x^{2M+1}}{2M+1},
  \quad 0<x<1,
\]
with \((x,M)=(1/5,5)\) and \((1/239,2)\), we obtain
\(355/113-10^{-6}<\pi<355/113\).
For \(x\geq0\), choose an integer \(k\geq0\) such that \(z=x/2^k\leq1/8\).
Set
\[
  L(z):=\sum_{j=0}^7\frac{(-z)^j}{j!},\qquad
  U(z):=\sum_{j=0}^6\frac{(-z)^j}{j!}.
\]
Since \(0<L(z)\leq e^{-z}\leq U(z)\) for \(0\leq z\leq1/8\),
\[
  L(x/2^k)^{2^k}\leq e^{-x}\leq U(x/2^k)^{2^k}.
\]
For square roots, we use
\(0\leq a\leq b,\ a^2\leq x\leq b^2\Longrightarrow a\leq\sqrt x\leq b\).
At each step, we replace a computed interval \([l,u]\ne[0,0]\) by
\[
  [l,u]\subset
  [h\lfloor l/h\rfloor,h\lceil u/h\rceil],\qquad
  h=10^{\max\{-40,\lfloor\log_{10}M\rfloor-6\}},\quad
  M=\max\{|l|,|u|\}.
\]
Thus the resulting interval still contains the quantity being estimated.

To evaluate \(M_r\) and \(D_{pr}\), we use
\eqref{nine:factored-difference}, \(a^{2j}=(a^j)^2\) and \(a^{2j+1}=a(a^j)^2\) for \(j\geq0\),
and the positive formulas
\eqref{nine:elementary-stable-integral}--\eqref{nine:positive-integral-recurrence}
for \(I_{1,\eta_\ell}\) and \(\mathcal J_p\). The quantities $\mathcal K_p$ are eliminated
by \eqref{nine:integration-by-parts}. We use \eqref{nine:fast-pairing-factors} in
\eqref{nine:fast-cross-hermite}. Each \(F_\ell^{[14]}(t)\) is evaluated by
\eqref{nine:finite-hardy-recurrence}, with its exact parameter
$24t\in\mathbb N$ used before rounding. Substitution gives the
rational bounds below. In the table, the values for
\(\sum_x Df\cdot Dg\) are upper bounds, and those for
\(\sum_{x\neq0}fg/|x|^2\) are lower bounds. Both sums include a factor
\(2\) for \(2Q(G,\mathcal B)\). By
\[
  Q(f,g)=\sum_{x\in\mathbb Z^9}Df(x)\cdot Dg(x)
  -\lambda\sum_{x\in\mathbb Z^9\setminus\{0\}}\frac{f(x)g(x)}{|x|^2},
  \qquad\lambda=\frac{1440}{121},
\]
we obtain the upper bounds for \(Q(G)\), \(2Q(G,\mathcal B)\), and
\(Q(\mathcal B)\) below.
\begin{equation}
  \label{nine:Q-rational-data}
  \begin{array}{c|c|c|c}
    \text{Term}&\sum_x Df\cdot Dg&\sum_{x\neq0}fg/|x|^2
    &\text{Upper bound}\\\hline
    Q(G)&\dfrac{157067}{500}&\dfrac{103819}{4000}&\dfrac{3282}{625}\\[6pt]
    2Q(G,\mathcal B)&\dfrac{1719877}{50000}&\dfrac{754577}{200000}
    &-\dfrac{65643}{6250}\\[6pt]
    Q(\mathcal B)&\dfrac{1854293}{25000}&\dfrac{1447863}{250000}
    &\dfrac{52487}{10000}
  \end{array}
\end{equation}
By \eqref{nine:trial-function} and \eqref{nine:Q-rational-data}, we obtain
\[
  \begin{aligned}
    Q(u)&=Q(G)+2Q(G,\mathcal B)+Q(\mathcal B)\\
    &\leq\frac{3282}{625}-\frac{65643}{6250}
    +\frac{52487}{10000}=-\frac{149}{50000}<0,
  \end{aligned}
\]
which proves \eqref{nine:elementary-final-bounds}.

\subsubsection[Approximation by functions in Cc(Z9)]{Approximation by functions in \(C_c(\mathbb Z^9)\)}
Finally, to use the definition of \(c_9\), we approximate \(u\) by
\(u_R\in C_c(\mathbb Z^9)\) with \(u_R(0)=0\).
By \(\beta_\tau(s)=\int_0^\tau b e^{-bs}\,db\) in \eqref{nine:trial-function}, we have
\(0<\beta_\tau(s)\leq s^{-2}\) and \(|\beta_\tau'(s)|\leq2s^{-3}\).
Thus, for all sufficiently large \(|x|\),
\begin{equation}
  \label{nine:trial-decay}
  |u(x)|\leq C|x|^{-4},\qquad |D_ju(x)|\leq C|x|^{-5}\quad(1\leq j\leq9).
\end{equation}
Let \(u_R=\chi_Ru\), where \(\chi_R\) is the cutoff from
Section~\ref{sec:sobolev-comparison}, now on \(\mathbb Z^9\).
Then \(u_R\in C_c(\mathbb Z^9)\) and \(u_R(0)=0\).
By \eqref{nine:trial-decay} and \(|D_j\chi_R|\leq R^{-1}\), we obtain
\[
  \begin{aligned}
    &\sum_x|D(u-u_R)(x)|^2
    +\sum_{x\ne0}\frac{|u(x)-u_R(x)|^2}{|x|^2}\\
    &\leq C\sum_{|x|>R-1}|x|^{-10}
    +\frac C{R^2}\sum_{R-1<|x|<2R+1}|x|^{-8}\\
    &\leq\frac CR\longrightarrow0\qquad(R\to\infty).
  \end{aligned}
\]
Therefore,
\[
  c_9\leq\lim_{R\rightarrow\infty}
  \frac{\displaystyle\sum_x|Du_R(x)|^2}
  {\displaystyle\sum_{x\neq0}|u_R(x)|^2/|x|^2}
  =\frac{\displaystyle\sum_x|Du(x)|^2}
  {\displaystyle\sum_{x\neq0}|u(x)|^2/|x|^2}<\frac{1440}{121}.
\]
This completes the proof.
\end{proof}

\noindent
\textbf{Conflict of interest.}
On behalf of all authors, the corresponding author states that there is no conflict of interest.

\medskip

\noindent
\textbf{Ethics approval.}
 Not applicable.

\medskip

\noindent
\textbf{Data Availability Statements.}
Data sharing is not applicable to this article as no datasets were generated or analysed during the current study.

\medskip

\noindent
\textbf{AI Declaration.}
The authors used ChatGPT (OpenAI, GPT-6) to improve the English and the clarity of the exposition in the preparation of this manuscript. All suggestions were independently
reviewed and verified by the authors, who take full responsibility for the content of the paper.

\medskip
\noindent
\textbf{Acknowledgements.}
C. Ji is supported by the National Natural Science Foundation of China (No.~12571117).

\end{document}